\documentclass[12pt]{amsart}
\usepackage{blindtext}
\usepackage{url}
\usepackage{tabularx}
\usepackage[utf8]{inputenc}
\usepackage[margin=2.5cm]{geometry}
\usepackage[utf8]{inputenc}
\usepackage{amsmath}
\usepackage{amsfonts}
\usepackage{color}
\usepackage{slashbox}
\usepackage{amssymb}
\usepackage{amsthm}
\usepackage{array}
\usepackage{rotating}
\usepackage{physics}
\usepackage{tikz}
\usepackage{braids}
\usepackage{tikz-cd}
\usepackage[bottom]{footmisc}
\tikzset{%
	symbol/.style={%
		draw=none,
		every to/.append style={%
			edge node={node [sloped, allow upside down, auto=false]{$#1$}}}
	}
}
\usepackage[all]{xy}
\usepackage{cite}
\usepackage{mathdots}
\usepackage{pstricks}
\usepackage{mathrsfs}
\usepackage{setspace}
\usepackage{hyperref}
\usepackage{fancyhdr}
\theoremstyle{theorem}
\newtheorem{theorem}{Theorem}
\numberwithin{theorem}{subsection}
\newtheorem{lemma}[theorem]{Lemma}
\newtheorem{proposition}[theorem]{Proposition}
\newtheorem{corollary}[theorem]{Corollary}
\theoremstyle{conjecture}

\theoremstyle{definition}
\newtheorem{definition}[theorem]{Definition}
\newtheorem{example}[theorem]{Example}
\theoremstyle{remark}
\newtheorem{remark}[theorem]{Remark}
\theoremstyle{diagram}

\theoremstyle{notation}
\newtheorem{notation}[theorem]{Notation}
\title[Strictifying operational coherences]{Strictifying operational coherences and weak functor classifiers in low dimensions}
\author{Adrian Miranda}
\thanks{This material is based upon work done for my PhD thesis. I am grateful for having been supported by the MQRES PhD Scholarship, 20192497, and by EPSRC under grant EP/V002325/2. I thank Steve Lack for his guidance while I was conducting this research, Nicola Gambino for his guidance while I was preparing this paper, and Calum Hughes for proof reading Section 2.}
\address{Department of Mathematics, University of Manchester, 
	United Kingdom}
\email{adrian.miranda@manchester.ac.uk}
\begin{document}

\maketitle

\begin{abstract}
	\noindent Weak structures abound in higher category theory, but are often suitably equivalent to stricter structures that are easier to understand. We extend strictification for tricategories and trihomomorphisms to trinatural transformations, trimodifications and perturbations. Along the way we distinguish between the operational coherences, which are possible to strictify, and the coherences on globular inputs, which remain weak. We introduce generalised path objects for $\mathbf{Gray}$-categories, which help reduce proofs in the three-dimensional setting to known results. Upon closing the resulting semi-strict trinatural transformations under composition, we state the hom-triequivalences of what we expect to be a `semi-strictification tetra-adjunction'.
\end{abstract}

\onehalfspacing

\section{Introduction}

\subsection{Motivation and Goals}
\noindent Many examples of higher categories of interest are \emph{weak} in the sense that they satisfy familiar axioms up to higher dimensional data, rather than on the nose. The same is true for higher dimensional functors, natural transformations and so on. Such structures can be difficult to work with in practice, but can often be understood via simpler structures to which they are suitably equivalent. These simpler structures are \emph{semi-strict}, in that they satisfy some axioms familiar from lower dimensional settings on the nose. Indeed, structures which satisfy all axioms from lower dimensional settings on the nose are typically too strict to model their fully weak counterparts up to weak equivalence \cite{GPS tricategory}. Identifying which kinds of coherence data can be made strict, and which must remain weak to avoid loss of generality, is a major open problem in higher category theory.
\\
\\
\noindent This paper studies particular well-behaved strictification constructions in low dimensions $n \in \{2, 3\}$. Specifically, any bicategory (resp. tricategory) is biequivalent (resp. triequivalent) to a $2$-category (resp. $\mathbf{Gray}$-category) whose underlying category (resp. sesquicategory) is free on a graph (resp. $2$-computad). We examine how the resulting freeness in codimension one of the strictified structure informs how the strictification construction extends to weak maps between higher categories, in which the usual functoriality or naturality axioms hold up to suitable equivalences. 
\\
\\
\noindent For $n=2$, these results about strictification are captured by properties of the left adjoint to the inclusion $I_{2}: \mathbf{Gray} \rightarrow \mathbf{Bicat}$, and of the $2$-adjunction and a triadjunction it underlies \cite{Campbell Strictification}. We present an alternative way of strictifying bicategories that is motivated by the universal property of the unit $\eta_{\mathcal{A}}: \mathcal{A} \rightsquigarrow \mathbf{st}_{2}\left(\mathcal{A}\right)$ as a universal pseudofunctor from $\mathcal{A}$ into a $2$-category. For $n=3$, a strictification construction for tricategories and trihomomorphisms has been studied in \cite{Gurski Coherence in Three Dimensional Category Theory}. One of the goals of extending three-dimensional strictification to higher dimensional weak maps between trihomomorphisms is to provide an analogous universal property.

\subsection{Main results and contributions}

\begin{itemize}
	\item Theorem \ref{2-globular set modelling weak higher cells via 3 pseudofunctors} shows how trinatural transformations, trimodifications and perturbations can be modelled via generalised path objects for $\mathbf{Gray}$-categories using partially weak versions of trihomomorphisms that we introduce in Definition \ref{Trihomomorphism special cases}.
	\item Lemma \ref{Gurski Lemma semi strict} extends Lemma 15.4 of \cite{Gurski Coherence in Three Dimensional Category Theory} to trimodifications and perturbations. As well as being a key tool for the rest of the paper, it is likely to be of interest in its own right.
	\item Theorem \ref{Gr on weak higher cells between tricategories} extends strictification of tricategories and trihomomorphisms to trinatural transformations, trimodifications and perturbations.
	\item Proposition \ref{tricategory of tricategories} combines our results with those of \cite{Low dimensional structures formed by tricategories} to give a tricategorical statement of three-dimensional strictification. The two-dimensional analogue is reviewed in Remark \ref{Two Dimensional Statement of Two Dimensional Strictification}.
	\item Proposition \ref{alternative composition if domain is cofibrant resolves semi-strict's failure to be closed under composition} shows that composition of semi-strict trinatural transformations out of cofibrant $\mathbf{Gray}$-categories can be redefined in such a way that their alternative composite is semi-strict. The newly defined composition extends to a $\mathbf{Gray}$-category structure that is triequivalent to the standard one. This corrects an error in Chapter 15 of \cite{Gurski Coherence in Three Dimensional Category Theory}.
	\item Definition \ref{coherent hom def} describes an internal hom ${[\mathfrak{A}, \mathfrak{B}]}_\text{ssg}$ of $\mathbf{Gray}$-categories, making no assumptions on the domain and codomain. In Theorem \ref{trihomomorphism classifier} we use this hom to describe triequivalences of the form depicted below, which we conjecture may feature in a tetra-adjunction between tricategories and $\mathbf{Gray}$-categories.
	
	$$\begin{tikzcd}
		{[}\mathbf{st}_{3}\left(\mathfrak{A}\right){,} \mathfrak{B}{]}_\text{ssg}
		\arrow[r, "I_{3}"]
		& {\mathbf{Tricat}\left(I_{3}\mathbf{st}_{3}\left(\mathfrak{A}\right){,} I_{3}\left(\mathfrak{B}\right)\right)}
		\arrow[r, "\eta_{\mathfrak{A}}^{*}"]
		& {\mathbf{Tricat}\left(\mathfrak{A}{,} I_{3}\left(\mathfrak{B}\right)\right)}
	\end{tikzcd}$$
	
	\noindent In forthcoming papers we will show that $[\mathfrak{A}, \mathfrak{B}]_\text{ssg}$ equips $\mathbf{Gray}$-$\mathbf{Cat}$ with a closed structure \cite{Miranda semi-strictly generated closed structure on Gray-Cat}, and pursue applications to four dimensional coherence via enrichment over a related base.
\end{itemize}

\subsection{Key ideas and techniques}

\subsubsection{Freeness in codimension one}

\noindent Strict functors out of free structures are uniquely determined via their values on generating data. Weak functors $F$ out of free structures also have values on such data, and thus determine a strict functor $\overline{F}$ with the same domain and codomain. However, weak functors may only preserve the operations, such as identities or composition, via which general data in their domain are freely generated, up to suitable equivalences. These coherence data themselves comprise a suitable equivalence between the original $F$ and the strict functor $\overline{F}$. Since strictifications of weak categories can be chosen to be free in codimension one, the construction also extends to weak functors turning them into suitably equivalent strict functors. These ideas are made precise in Remark \ref{pseudofunctors out of cofibrant 2-categories} and Proposition \ref{trihomomorphism out of cofibrant Gray category}. Although homotopy theoretic aspects are not the focus of this paper, the property of being free in codimension one characterises cofibrant objects in the relevant model structures on $2$-$\mathbf{Cat}$ \cite{Quillen 2-cat} and $\mathbf{Gray}$-$\mathbf{Cat}$ \cite{Quillen Gray Cat}.

\subsubsection{Operational versus free coherences}
We distinguish between two different types of coherences in weak maps for higher categories. The first type, which we call \emph{free coherences}, mediate respect for the globular structure of the domain $n$-dimensional category. These coherences categorify naturality-like conditions that hold on the nose for their lower dimensional analogues. The second type of coherences are the \emph{operational coherences} of the title. These mediate respect for the various operations of $\mathscr{A}$, such as composition, identities and whiskering. Our motivation for making this distinction between different types of coherences is that when strictification is extended to weak functors, weak transformations and weak modifications, operational coherences become identities while free coherences remain weak. Distinguishing between different types of coherences in this way clarifies properties expected in semi-strict structures.

\subsubsection{Generalised path objects for $\mathbf{Gray}$-categories}

\noindent In Subsection \ref{Path Objects Section} we introduce generalised path objects for $\mathbf{Gray}$-categories. These allow us to leverage the strictification trihomomorphisms discussed in the previous paragraph to prove analogous results about higher dimensional weak maps between trihomomorphisms. In particular, there is a sense in which trinatural transformations between $\mathbf{Gray}$-functors are like pseudofunctors, despite seeming to be much more complicated. This is made precise in Proposition \ref{trinatural transformations as 2-Cat pseudofunctors}.

\subsubsection{Semi-strictly decomposable trinatural transformations}

\noindent In dimension two, part of the interest in semi-strict maps between $2$-categories is that they form an internal hom $\mathbf{Gray}\left(\mathcal{A}, \mathcal{B}\right)$ which features in a closed structure on $2$-$\mathbf{Cat}$ and captures semi-strict categories in dimension three via enrichment. As observed in \cite{crans tensor of gray categories}, semi-strict trinatural transformations fail to be closed under composition due to failure of middle four interchange in the codomain $\mathbf{Gray}$-category, posing challenges in constructing a similar internal hom between $\mathbf{Gray}$-categories. We overcome this challenge by closing semi-strict trinatural transformations under composition, and describe a new semi-strictly generated internal hom of $\mathbf{Gray}$-categories $[\mathfrak{A}, \mathfrak{B}]_\text{ssg}$.

\subsection{Outline}
\noindent The plan for this paper is as follows.

\begin{itemize}
	\item Section \ref{Section operational coherences in dimension two} reviews strictification of bicategories, pseudofunctors, pseudonatural transformations and modifications. We emphasise the new perspective provided by the distinction between operational and free coherences.
	\item Section \ref{Section Definitions of 3, k transfors} reviews the definitions of trihomomorphisms, and of trinatural transformations, trimodifications and perturbations between $\mathbf{Gray}$-functors. We again emphasise the new distinction between their operational and free coherence data, in light of which one may correctly anticipate that trinatural transformations can be partially strictified. This section also catalogues the structure of the full sub-$\mathbf{Gray}$-categories of $\mathbf{Tricat}\left(\mathfrak{A}, \mathfrak{B}\right)$ on $\mathbf{Gray}$-functors, which will be used throughout the rest of the paper.
	\item Section \ref{Section Extending strictification to 3, k transfors} extends Gurski's strictification construction for tricategories and trihomomorphisms to higher $\left(3, k\right)$-transfors. This requires developing the notions of generalised path objects for $\mathbf{Gray}$-categories.
	\item Section \ref{Section Solutios to the Failure of Semi-strict trinatural transformations to be closed under composition} explores solutions to the failure of semi-strict trinatural transformations to be closed under composition.
\end{itemize}

\begin{notation}
	\begin{itemize}
		\item 
		Two-dimensional structures, such as bicategories or $2$-categories, will be written in calligraphic font $\mathcal{A}$, $\mathcal{B}$, $\mathcal{C}$. Meanwhile, three-dimensional structures, such as tricategories or $\mathbf{Gray}$-categories, will be written in fraktur font $\mathfrak{A}$, $\mathfrak{B}$, $\mathfrak{C}$.
		\item Objects will typically be written in upper-case roman letters $X$, $Y$, $Z$, morphisms will be typically written in lower-case roman letters $f$, $g$, $h$, $2$-cells will be typically written in lower-case Greek letters $\phi$, $\psi$, and $3$-cells will be typically written in upper-case Greek letters $\Gamma$, $\Omega$. Data in mapping spaces such as $\mathbf{Bicat}\left(\mathcal{A}, \mathcal{B}\right)$ or $\mathbf{Tricat}\left(\mathfrak{A}, \mathfrak{B}\right)$ will be written with the same script and case conventions, according to their dimension.
	\end{itemize}
\end{notation}

\section{Operational coherences of $\left(2, k\right)$-transfors and a review of the strictification triadjunction}\label{Section operational coherences in dimension two}

\noindent In a bicategory \cite{Benabou Bicategories}, the usual associativity and left and right unit laws for categories hold up to invertible $2$-cells which are subject to further axioms. Pseudofunctors between bicategories obey the usual unit and composition laws for functors between categories up to invertible $2$-cells subject to further axioms. Similarly, pseudonatural transformations only obey the usual naturality condition up to suitably coherent invertible $2$-cells. Just as the morphisms of a category allow natural transformations between functors to be defined, the $2$-cells of a $2$-category allow higher transfors known as \emph{modifications} to be defined between pseudonatural transformations. These fit into a larger picture of higher maps between higher categories as the $n =2$ column of the table below. Such higher maps are often referred to as $\left(n, k\right)$-transfors, following \cite{Crans Localisations of Transfors}, \cite{crans tensor of gray categories}.
\\
\\
\emph{`Periodic Table' of pseudo $\left(n, k\right)$ transfors for $0 \leq n \leq 3$ and $0 \leq k \leq 4$.}
	\\
	\\
	\begin{tabular}{|>{\raggedright}p{0.07\linewidth}|>{\raggedright}p{0.2\linewidth}|>{\raggedright}p{0.2\linewidth}|>{\raggedright}p{0.2\linewidth}|>{\raggedright\arraybackslash}p{0.2\linewidth}|}
		\hline
		$n = $ & $0$& $1$ & $2$ & $3$  \\
		\hline
		$k = 0$
		& functions
		& functors
		& pseudofunctors
		& trihomomorphisms \\
		\hline 
		$k = 1$
		& equalities of functions
		& natural transformations
		& pseudonatural transformations 
		& trinatural transformations
		\\
		\hline
		$k = 2$
		& $-$
		& equalities of natural transformations
		& modifications
		& trimodifications
		\\
		\hline
		$k = 3$
		& $-$
		& $-$
		& equalities of modifications
		& perturbations 
		\\
		\hline
		$k=4$
		& $-$
		&$-$
		& $-$
		& equalities of perturbations
		\\
		\hline
	\end{tabular}

\noindent The data in the table above are of the `pseudo' variant, where $\left(n, k\right)$-transfors satisfy the axioms for an $\left(n-1, k\right)$-transfor up to appropriate $n$-cell equivalences. Subsection \ref{subsection definitions of 2 k transfors} reviews the definitions for data in the $n=2$ column while Section \ref{Section Definitions of 3, k transfors} reviews definitions for data in the column $n =3$. There are many pieces of data involved in $\left(n, k\right)$-transfors, and the goal of our presentation of this information will be to make distinctions between different kinds of data in a way that relates to strictification. In particular, certain data can be understood at the level of underlying globular sets, and will be referred to as \emph{globular data}. Meanwhile, other data mediate the axioms of an $\left(n-1, k\right)$-transfor. When $n=2$ these axioms are the composition and identity preservation aspects of functoriality, and the naturality condition. Data of this kind are referred to as \emph{coherences}. Coherence data will be further split into two different classes. 

\begin{itemize}
	\item \emph{free} coherences, which will mediate `naturality like' conditions, in which $m$-cells in the domain vary along $\left(m+1\right)$-cells. In our current setting, there will only be one such coherence. Coherences of this kind will remain weak after semi-strictification.
	\item \emph{operational} coherences, which will mediate respect for the operations in the underlying $\left(n-1\right)$-dimensional categorical structure of domain. In our current setting with $n =2$, these include the nullary identity operation and the binary composition operation. When the domain is freely generated, these coherences will be strictifiable.
\end{itemize}

\noindent These distinctions will be explained in greater detail in Subsection \ref{subsection free versus derivable coherences}. Analogous distinctions will be made in the three dimensional setting in Remark \ref{Remark free vs operational three dimensions}, and there will be a similar relationship to strictifiability in that setting.
\\
\\
\noindent Throughout Subsections \ref{subsection definitions of 2 k transfors} and \ref{subsection free versus derivable coherences}, $\mathcal{A}$ and $\mathcal{B}$ will be $2$-categories.

\subsection{Definitions of $\left(2, k\right)$-transfors between $2$-categories}\label{subsection definitions of 2 k transfors}

\begin{definition}\label{definition pseudofunctor}
	A \emph{pseudofunctor} $\left(F, \iota, \chi\right): \mathcal{A} \rightarrow \mathcal{B}$ consists of \begin{itemize}
		\item A function $F_{0}: \mathcal{A}_{0} \rightarrow \mathcal{B}_{0}$,
		\item For every $X, Y \in \mathcal{A}$, a functor $F_{X, Y}: \mathcal{A}\left(X, Y\right) \rightarrow \mathcal{B}\left(FX, FY\right)$. Hence in particular for every $f: X \rightarrow Y$, there is a morphism $Ff: FX \rightarrow FY$, and for every $2$-cell $\phi: f \Rightarrow g: X \rightarrow Y$, there is a $2$-cell $F\phi: Ff \Rightarrow Fg$,
		\item For every identity arrow $1_X \in \mathcal{A}$, an invertible $2$-cell $\iota_{X}: 1_{FX} \Rightarrow F1_{X}$ called the \emph{unitor},
		\item For every composable pair \begin{tikzcd}
			X \arrow[r, "f"] & Y \arrow[r, "g"] & Z
		\end{tikzcd}, an invertible $2$-cell $\chi_{g, f}: Fg.Ff \Rightarrow F\left(gf\right)$, called the \emph{compositor}.
	\end{itemize}
	\noindent These data are subject to the condition that $\chi_{g, f}$ is natural as either $g$ or $f$ vary along $2$-cells, as in the equations below.
	
	$$\begin{tikzcd}[ font=\fontsize{9}{6}]
		Fg.Ff
		\arrow[dd, "1_{Fg}.F\phi"']
		\arrow[rr, "\chi_{{}_{g{,}f}}"]
		&&F\left(gf\right)
		\arrow[dd, "F\left(1_g .\phi\right)"]
		\\
		\\
		Fg.Ff'
		\arrow[rr, "\chi_{{}_{g{,}f'}}"']
		&&F\left(gf'\right)&{}
	\end{tikzcd}\begin{tikzcd}[ font=\fontsize{9}{6}]
		Fg.Ff
		\arrow[dd, "F\psi.1_{Ff}"']
		\arrow[rr, "\chi_{{}_{g{,}f}}"]
		&&F\left(gf\right)
		\arrow[dd, "F\left(\psi .1_f\right)"]
		\\
		\\
		Fg'.Ff
		\arrow[rr, "\chi_{{}_{g'{,}f}}"']
		&&
		F\left(g'f\right)
	\end{tikzcd}$$
	\noindent Moreover, the following three diagrams must commute. These conditions will be referred to as the \emph{left unit law}, \emph{associativity law}, and \emph{right unit law} respectively.
	
	$$ \begin{tikzcd}[ font=\fontsize{9}{6}]
		1_{FY}.Ff
		\arrow[dd, equal]
		\arrow[rr, "\iota_{Y}.1_{Ff}"]
		&& F1_{Y}.Ff
		\arrow[dd, "\chi_{{}_{1_{Y}{,}f}}"]
		\\
		\\
		Ff\arrow[rr, equal]
		&&F\left(1_{Y}.f\right)
	\end{tikzcd} \begin{tikzcd}[ font=\fontsize{9}{6}]
		Fh.Fg.Ff
		\arrow[dd, "1_{Fh}.\chi_{{}_{g{,}f}}"']
		\arrow[rr, "\chi_{{}_{h{,} g}}.1_{Ff}"] && F\left(hg\right).Ff
		\arrow[dd, "\chi_{{}_{hg{,}f}}"]
		\\
		\\
		Fh.F\left(gf\right)
		\arrow[rr, "\chi_{{}_{h{,}gf}}"']
		&& F\left(hgf\right)
	\end{tikzcd} \begin{tikzcd}[ font=\fontsize{9}{6}]
		Ff.1_{FX}
		\arrow[dd, equal]
		\arrow[rr, "1_{Ff}.\iota_{X}"]
		&& Ff.F1_{X}
		\arrow[dd, "\chi_{{}_{f{,}1_{X}}}"]
		\\
		\\
		Ff\arrow[rr, equal]
		&&F\left(f.1_{X}\right)
	\end{tikzcd} $$
	
	\noindent If $\iota$ is the identity, $F$ will be called \emph{normal}, and if $\chi$ is the identity then $F$ will be called \emph{compositional}. When the unitors and compositors are clear from context, $\left(F, \iota, \chi\right)$ will just be denoted as $F$.
\end{definition}

\noindent Pseudofunctors have an underlying morphism of $2$-globular sets $F_{n}: \mathcal{A}_{n} \rightarrow \mathcal{B}_{n}$, but also have coherences given by their compositors and unitors. These coherences mediate the usual functoriality laws, which specify respect for the underlying category structure of $\mathcal{A}$. Observe that Definition \ref{definition pseudofunctor} reduces to the notion of a $2$-functor if $\chi$ and $\iota$ are taken to be identities.

\begin{definition}\label{definition pseudonatural transformation}
	Let $\left(F, \chi, \iota\right)$ be as above and let $\left(G, \eta, \mu\right)$ be another pseudofunctor from $\mathcal{A}$ to $\mathcal{B}$. A \emph{pseudonatural transformation} $p: F \Rightarrow G$ consists of the following data subject to the following axioms.
	\\
	\\
	\noindent DATA
	\\
	\begin{itemize}
		\item For every $X \in \mathcal{A}$, a morphism $p_{X}: FX \rightarrow GX$ in $\mathcal{B}$, which will be referred to as the \emph{$1$-cell component} of $p$ on $X$.
		\item For every $f: X \rightarrow Y$ in $\mathcal{A}$, an invertible $2$-cell $p_{f}: Gf.p_{X} \Rightarrow p_{Y}.Ff$, which will be referred to as the \emph{$2$-cell component}, or the \emph{pseudonaturality constraint}, of $p$ on $f$.
	\end{itemize}
	
	\noindent AXIOMS
	\\
	\begin{itemize}
		\item \emph{(Local naturality)} The assignment $f \mapsto p_{f}$ varies naturally as $f$ varies along $\phi: f \Rightarrow g$. This is depicted in the equality shown below left.
		\item \emph{(Unit law)} The equality of pasting diagrams depicted below right holds.
	\end{itemize}
	
	$$
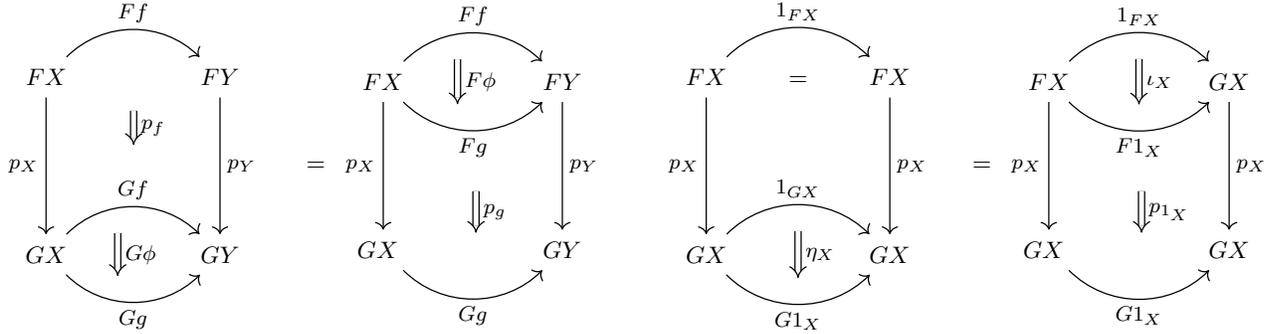
\begin{tikzcd}[column sep = 15, row sep = 20,  font=\fontsize{9}{6}]
		FX
		\arrow[rr,bend left = 45, "Ff"]
		\arrow[dd,"p_X"']
		&
		{}
		\arrow[d,Rightarrow, shorten = 5,"p_f"]
		&
		FY
		\arrow[dd,"p_Y"]
		\\
		&{}&&=
		\\
		GX
		\arrow[rr,bend left = 45,"Gf" {name = C}]
		\arrow[rr, bend right = 45, "Gg"' {name = D}]
		&
		{}
		&
		GY
		\arrow[from =C, to =D, Rightarrow, shorten = 10, shift right = 2, "G\phi"]
	\end{tikzcd}	\begin{tikzcd}[column sep = 15, row sep = 20,  font=\fontsize{9}{6}]
		FX
		\arrow[rr,bend left = 45, "Ff" {name = A}]
		\arrow[rr, bend right = 45, "Fg"' {name = B}]
		\arrow[dd,"p_X"']
		&
		&
		FY
		\arrow[dd,"p_{Y}"]
		\\
		&{}
		\arrow[d,Rightarrow, shorten = 5,"p_g"]
		\\
		GX\
		\arrow[rr, bend right = 45, "Gg"']
		&
		{}
		&
		GY &{}
		\arrow[from =A, to =B, Rightarrow, shorten = 10, shift right = 2, "F\phi"]
	\end{tikzcd}\begin{tikzcd}[column sep = 13, row sep = 20,  font=\fontsize{9}{6}]
		FX
		\arrow[rr,bend left = 45, "1_{FX}"]
		\arrow[dd,"p_{X}"']
		&
		=
		&
		FX
		\arrow[dd,"p_{X}"]
		\\
		&{}&&=
		\\
		GX
		\arrow[rr,bend left = 45,"1_{GX}" {name = C}]
		\arrow[rr, bend right = 45, "G{1_{X}}"' {name = D}]
		&
		{}
		&
		GX
		\arrow[from =C, to =D, Rightarrow, shorten = 10, "\eta_{X}"]
	\end{tikzcd}	\begin{tikzcd}[column sep = 15, row sep = 20,  font=\fontsize{9}{6}]
		FX
		\arrow[rr,bend left = 45, "1_{FX}" {name = A}]
		\arrow[rr, bend right = 45, "F{1_{X}}"' {name = B}]
		\arrow[dd,"p_{X}"']
		&
		&
		GX
		\arrow[dd,"p_{X}"]
		\\
		&{}
		\arrow[d,Rightarrow, shorten = 5,"p_{1_{X}}"]
		\\
		GX\
		\arrow[rr, bend right = 45, "G{1_{X}}"']
		&
		{}
		&
		GX &{}
		\arrow[from =A, to =B, Rightarrow, shorten = 10, "\iota_{X}"]
	\end{tikzcd} $$

	\begin{itemize}
		\item \emph{(Composition law)} The equality of pasting diagrams depicted below holds.
		
		$$
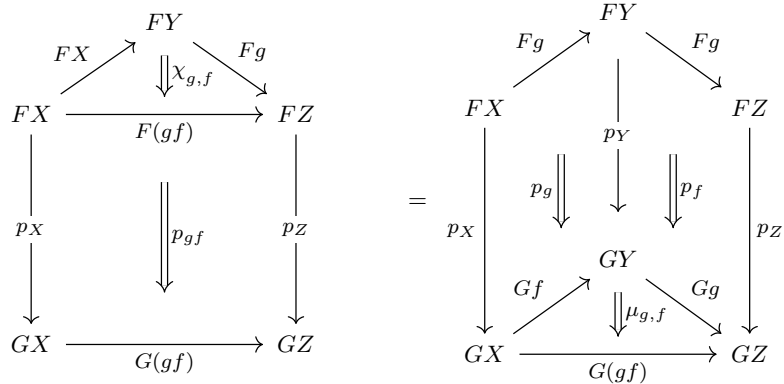
\begin{tikzcd}[row sep = 20, column sep = 25,  font=\fontsize{9}{6}]
			&FY
			\arrow[rd, "Fg"]
			\arrow[d, Rightarrow, shorten = 5, "\chi_{{}_{g{,}f}}"]
			\\
			FX
			\arrow[ddd, "p_{X}" description]
			\arrow[ru, "F{X}"]
			\arrow[rr, "F\left({gf}\right)"']
			&{}\arrow[ddd, Rightarrow, shorten = 18, "p_{gf}"]
			& FZ
			\arrow[ddd, "p_{Z}" description]
			\\
			&&&=
			\\
			\\
			GX 
			\arrow[rr, "G\left({gf}\right)"']
			&{}
			& GZ
		\end{tikzcd}\begin{tikzcd}[column sep = 25,  font=\fontsize{9}{6}]
			&FY
			\arrow[rd, "Fg"]
			\arrow[ddd, shorten = 10, "p_{Y}" description]
			\\
			FX
			\arrow[dd, Rightarrow, shorten = 10, "p_{g}"', shift left = 10]
			\arrow[ddd, "p_{X}"']
			\arrow[ru, "Fg"]
			&& FZ
			\arrow[dd, Rightarrow, shorten = 10, "p_{f}", shift right = 10]
			\arrow[ddd, "p_{Z}"]
			\\
			\\
			{}&{}GY
			\arrow[rd, "Gg"]
			\arrow[d, Rightarrow, shorten = 5, "\mu_{g{,}f}"]
			&{}
			\\
			GX
			\arrow[ru, "Gf"]
			\arrow[rr, "G\left({gf}\right)"']
			&{}
			& GZ
		\end{tikzcd}$$
	\end{itemize}
	
	\noindent Pseudonatural transformations with identity $1$-cell components will be referred to as \emph{invertible icons}, while pseudonatural transformations with identity $2$-cells will be called \emph{strict transformations}. 
\end{definition}

\noindent Both components of a pseudonatural transformation are determined on globular input data. The $2$-cell component mediates the usual naturality condition, as a $0$-cell varies along a $1$-cell, rather than respect for any of the operations in the underlying category of $\mathcal{A}$. Observe that strict transformations between $2$-functors are precisely $2$-natural transformations.

\begin{definition}\label{Definition modification}
	Let $p: F \Rightarrow G$ be as above and let $q: F \Rightarrow G$ be another pseudonatural transformation. A \emph{modification} $\sigma: p \Rrightarrow q$ consists of the data of a $2$-cell $\sigma_{X}: p_{X} \Rightarrow q_{X}$, subject to the following equation of pasting diagrams.
	
	$$\begin{tikzcd}[column sep = 18, row sep = 20,  font=\fontsize{9}{6}]
		FX
		\arrow[rr,bend left = 45, "p_{X}"]
		\arrow[dd,"Ff"']
		&
		{}
		\arrow[d,Rightarrow, shorten = 5,"p_f"]
		&
		GX
		\arrow[dd,"Gf"]
		\\
		&{}&&=
		\\
		FY
		\arrow[rr,bend left = 45,"p_{Y}" {name = C}]
		\arrow[rr, bend right = 45, "q_{Y}"' {name = D}]
		&
		{}
		&
		GY
		\arrow[from =C, to =D, Rightarrow, shorten = 10, shift right = 5, "\sigma_{Y}"]
	\end{tikzcd}	\begin{tikzcd}[column sep = 18, row sep = 20,  font=\fontsize{9}{6}]
		FX
		\arrow[rr,bend left = 45, "p_{X}" {name = A}]
		\arrow[rr, bend right = 45, "q_{X}"' {name = B}]
		\arrow[dd,"Ff"']
		&
		&
		GX
		\arrow[dd,"q_{f}"]
		\\
		&{}
		\arrow[d,Rightarrow, shorten = 5,"q_{f}"]
		\\
		FY\
		\arrow[rr, bend right = 45, "q_{Y}"']
		&
		{}
		&
		GY &{}
		\arrow[from =A, to =B, Rightarrow, shorten = 10, shift right = 5, "\sigma_{X}"]
	\end{tikzcd}$$ 
	
\end{definition}

\subsection{Free versus operational coherences}\label{subsection free versus derivable coherences}

\begin{definition}\label{globular transfors}
	Let $\mathbb{A}$ and $\mathbb{B}$ be $n$-globular sets, with source and target functions denoted $d_{1}, d_{0}: \mathbb{A}_{k+1} \rightarrow \mathbb{A}_{k}$ for $0 < k < n$ and similarly for $\mathbb{B}$. Globular $\left(n, k\right)$-transfors are defined by the following induction.
	
	\begin{itemize}
		\item Globular $\left(n, 0\right)$-transfors $F: \mathbb{A} \rightarrow \mathbb{B}$ are morphisms of globular sets.
		\item Let $S$ and $T$ be globular $\left(n, k\right)$-transfors from $\mathbb{A}$ to $\mathbb{B}$ for $0 \leq k < n$. A globular $\left(n, k+1\right)$-transfor $\Phi: S \rightarrow T$ consists of a function $\Phi_{0}: \mathbb{A}_{0} \rightarrow \mathbb{B}_{k+1}$ such that $d_{0}\Phi_{0} = T_{0}$ and $d_{1}\Phi_{0} = S_{0}$.
	\end{itemize}
\end{definition}

\noindent In general, $\left(n, k\right)$-transfors between weak-$n$-categories are expected to have underlying globular $\left(n, k\right)$-transfors. This is true in dimensions two and three, where precise definitions are known. The data in the underlying globular transfors of a $\left(2, k\right)$ transfor will be called its \emph{globular} data, while any other data will be called its \emph{coherences}. Coherence data for $\left(2, k\right)$ transfors are listed below.

\begin{itemize}
	\item $\left(k = 0\right)$ The assignments $\left(g, f\right)\mapsto \big(\chi_{g, f}: Fg.Ff \Rightarrow F\left(gf\right)\big)$ and $X \mapsto \big(\iota_{X}:1_{FX}\Rightarrow F1_{X}\big)$, which mediate respect for composition and identity laws respectively.
	\item $\left(k = 1\right)$ The assignment $f \mapsto \big(p_{f}: Gf.p_{X} \Rightarrow p_{Y}.Ff\big)$, which mediates naturality as $X$ varies along $f$.
	\item $\left(k = 2\right)$ There are no coherence data.
\end{itemize}

\noindent As anticipated, it is useful to further distinguish between two different kinds of coherence data. There are coherences which mediate respect for the operations in the underlying category of the domain (composition $\left(g, f\right)\mapsto \chi_{g, f}$ and identities $X \mapsto\iota_{X}$), and ones which mediate varying an $m$-cell in the domain along an $\left(m+1\right)$-cell (only applicable when $m=0$, and given by $f \mapsto p_{f}$).

\begin{remark}\label{free vs operational via dimensional analysis}
	As a good first approximation, the free coherences of an $\left(n, k\right)$ transfor are just those which are $\left(m+k\right)$-dimensional outputs for $m$-dimensional inputs, with these inputs being globular in shape. For example, $\chi_{g, f}$ is determined by the non-globular input data \begin{tikzcd}
		X \arrow[r, "f"] & Y \arrow[r, "g"] &Z
	\end{tikzcd}, so it is an operational coherence. Meanwhile $\iota_{X}$ is determined by a $0$-cell $X$, hence globular input data, but it is not itself $0+0$ dimensional, so it is also operational. In contrast, $p_{f}$ is a $m+k= 2$-cell component determined by a $k = 1$-transfor on a $m = 1$-cell in the domain, so it is a free coherence. In the two-dimensional setting, this `dimensional analysis' does correctly classify coherences. However, as we will discuss in Remark \ref{Remark free vs operational three dimensions}, it will require a slight caveat when we move to the three-dimensional setting to account for the fact that the analogue of $p_{f}$ will underlie an adjoint equivalence.
\end{remark}

\subsection{Strictification of bicategories as a left triadjoint to $I_{2}:\mathbf{Gray} \rightarrow \mathbf{Bicat}$}\label{Subsection strictification of bicategories as a triadjoint}

\noindent The monomorphism of categories $I_{2}: 2$-$\mathbf{Cat} \rightarrow \mathbf{Bicat}$ has a left adjoint given by strictification of bicategories $\mathcal{B} \mapsto \mathbf{st}_{2}\left(\mathcal{B}\right)$. Definition \ref{presentation of universal 2-category recipient of a pseudofunctor} describes a particular presentation for such a $2$-category $\mathbf{st}_{2}\left(\mathcal{B}\right)$ formed from a bicategory $\mathcal{B}$. This presentation is directly motivated by the universal property of the unit. See \cite{Street Limits Indexed By Category Values 2-functors} for presentations of $2$-categories via $2$-computads.

\begin{definition}\label{presentation of universal 2-category recipient of a pseudofunctor}
	The $2$-category $\mathbf{st}_{2}\left(\mathcal{B}\right)$ is defined by the following presentation.
	
	\begin{itemize}
		\item Its objects are as in $\mathcal{B}$,
		\item Its generating morphisms are morphisms $f: X \rightarrow Y$ in $\mathcal{B}$. The underlying category of $\mathbf{st}_{2}\left(\mathcal{B}\right)$ is freely generated from the graph we have thus far described. Hence \begin{itemize}
			\item the identity on an object $X$ is given by the empty path, which we denote as $\left(-\right)_{X}: X \rightarrow X$ 
			\item a general morphism is a path of the form depicted below, and will be denoted $\left(f_{m}, ..., f_{1}\right): X_{0} \rightarrow X_{m}$.
		\end{itemize} 
		
		$$\begin{tikzcd}
			X_{0} \arrow[r, "f_{1}"] & X_{1} \arrow[r, "f_{2}"] &... \arrow[r, "f_{m}"] & X_{m}
		\end{tikzcd}$$
		
		\item It has three types of generating $2$-cells: 
		\begin{enumerate}
			\item For every $2$-cell $\phi: f \Rightarrow g: X \rightarrow Y$ in $\mathcal{B}$, there is a generating $2$-cell $\phi: f \Rightarrow g$,
			\item For every object $X \in \mathcal{B}$, there are two generating $2$-cells $\iota_{X}: \left(-\right)_{X} \Rightarrow 1_{X}$ and $\iota_{X}^{-1}: 1_{X} \Rightarrow \left(-\right)_{X}$. 
			\item For every composable pair \begin{tikzcd}
				X \arrow[r, "f"]&Y\arrow[r, "g"]&Z\end{tikzcd}, there are two generating $2$-cells $\chi_{g, f}: \left(g, f\right) \Rightarrow gf$ and $\chi_{g, f}^{*}: gf\Rightarrow \left(g, f\right)$.
		\end{enumerate}
	\end{itemize}
	\noindent The free $2$-category generated by the $2$-computad thus far described is subject to the following relations in $\mathbf{st}_{2}\left(\mathcal{B}\right)$. These relations specify that the pseudofunctoriality axioms hold, and that composition and identities of data in the hom-categories $\mathcal{B}\left(X, Y\right)$ agree with the same operations performed in the categories $\mathbf{st}_{2}\left(\mathcal{B}\right)\left(X, Y\right)$.
	\begin{enumerate}
		\item For every object $X \in \mathcal{B}$, there are two relations specifying that $\iota_{X}$ and $\iota_{X}^{-1}$ are each other's inverse: $\iota_{X}^{-1}\iota_{X}= \left(-\right)_{1_{X}}$ and $\iota_{X}\iota_{X}^{-1} = \left(-\right)_{\left(-\right)_{X}}$.
		\item For every arrow $f: X \rightarrow Y$ in $\mathcal{B}$, there are three relations.
		\begin{enumerate}
			\item $\left(-\right)_{f} = 1_{f}$, where $\left(-\right)_{f}$ is the identity $2$-cell on $f$ in the free $2$-category.
			\item The following composite of $2$-cells equals $\lambda_{f}$.
			
			$$\begin{tikzcd}
				f \arrow[r, equal]
				&\left(-\right)_{Y}f
				\arrow[rr, "\left(\iota_{Y}{,}f\right)"]
				&&
				\left(1_{Y}{,}f\right)
				\arrow[rr, "\chi_{1_{Y}{,}f}"]
				&& 1_{Y}f
			\end{tikzcd}$$
			\item The following composite of $2$-cells equals $\rho_{f}$.
			
			$$\begin{tikzcd}
				f \arrow[r, equal]
				&f\left(-\right)_{X}
				\arrow[rr, "\left(f{,}\iota_{X}\right)"]
				&&
				\left(f{,}1_{X}\right)
				\arrow[rr, "\chi_{f{,}1_{X}}"]
				&& f1_{X}
			\end{tikzcd}$$
		\end{enumerate}
		\item For every $\phi: f \Rightarrow f': X \rightarrow Y$ and $\psi: g \Rightarrow g': Y \rightarrow Z$, the naturality condition for $\chi$ holds on the morphisms $\left(1_{g}, \phi\right): \left(g, f\right) \rightarrow \left(g, f'\right)$ and $\left(\psi, 1_{f}\right): \left(g, f\right) \rightarrow \left(g', f\right)$ in $\mathcal{B}\left(Y, Z\right)\times \mathcal{B}\left(X, Y\right)$.
		
		$$\begin{tikzcd}
			\left(g{,}f\right)
			\arrow[dd, "\left(1_g{,}\phi\right)"']
			\arrow[rr, "\chi_{{}_{g{,}f}}"]
			&&gf\arrow[dd, "g\phi"]
			\\
			\\
			\left(g{,}f'\right)
			\arrow[rr, "\chi_{{}_{g{,}f'}}"']
			&&gf'&{}
		\end{tikzcd}\begin{tikzcd}
			\left(g{,}f\right)
			\arrow[dd, "\left(\psi{,}1_f\right)"']
			\arrow[rr, "\chi_{{}_{g{,}f}}"]
			&&gf\arrow[dd, "\psi f"]
			\\
			\\
			\left(g'{,}f\right)
			\arrow[rr, "\chi_{{}_{g'{,}f}}"']
			&&g'f
		\end{tikzcd}$$

		\item For every composable pair \begin{tikzcd}
			f \arrow[r, "\phi"] &g \arrow[r, "\psi"] & h
		\end{tikzcd} in $\mathcal{A}\left(X, Y\right)$, denote by $\psi\phi$ their vertical composite in $\mathcal{B}$. Then the equation $\psi\phi = \left(\psi, \phi\right)$ holds, where the right hand side is the formal composite in $\mathbf{st}_{2}\left(\mathcal{B}\right)$.
		\item For every \begin{tikzcd}
			W \arrow[r, "f"]
			&X\arrow[r, "g"]
			&Y \arrow[r, "h"]
			&Z
		\end{tikzcd}, the equation depicted below holds.
		$$\begin{tikzcd}
			\left(h{,}g{,}f\right)
			\arrow[dd, equal]
			\arrow[rr, "\left(\chi_{{}_{h{,}g}}{,} f\right)"]
			&& \left(hg{,}f\right)
			\arrow[rr, "\chi_{{}_{hg{,}f}}"]
			&& \left(hg\right)f
			\arrow[dd, "{\alpha}_{{}_{h{,}g{,}f}}"]
			\\
			\\
			\left(h{,}g{,}f\right)
			\arrow[rr, "\left(h{,}\chi_{{}_{g{,}f}}\right)"']
			&&\left(h{,}gf\right)
			\arrow[rr, "\chi_{{}_{h{,}gf}}"']
			&&
			h\left(gf\right)
		\end{tikzcd}$$
	\end{enumerate}
\end{definition}

\begin{remark}\label{comparing pseudofunctor classifier to cofibrant replacement}
	 Our $\mathbf{st}_{2}\left(\mathcal{B}\right)$ is a biased version of the $2$-category described in Section 2.2.3 of \cite{Gurski Coherence in Three Dimensional Category Theory}. Instead of arbitrarily choosing a way of evaluating paths, and then defining $2$-cells in $\mathbf{st}_{2}\left(\mathcal{B}\right)$ to be $2$-cells in $\mathcal{B}$ between the resulting morphisms, we have included formal compositors $\chi_{g, f}$ and unitors $\iota_{X}$ in the presentation given in Definition \ref{presentation of universal 2-category recipient of a pseudofunctor}. It is evident from Definition \ref{presentation of universal 2-category recipient of a pseudofunctor} that there is a pseudofunctor $\eta_{\mathcal{B}}: \mathcal{B} \rightsquigarrow \mathbf{st}_{2}\left(\mathcal{B}\right)$, whose unitors and compositors are the generating $2$-cells of the form $\chi_{g, f}$ and $\iota_{X}$. Moreover, since the relations in the presentation of the $2$-category $\mathbf{st}_{2}\left(\mathcal{B}\right)$ are precisely those needed for $\eta_\mathcal{B}$ to be a pseudofunctor, it is also clear that it will be the \emph{universal} pseudofunctor from $\mathcal{B}$ to some $2$-category. An isomorphism between our $\mathbf{st}_{2}\left(\mathcal{B}\right)$ and the one of \cite{Gurski Coherence in Three Dimensional Category Theory} follows from the universal property.
\end{remark}

\noindent There is also a $2$-categorical formulation of strictification of bicategories and pseudofunctors, first given in \cite{Icons}, which we review in Remark \ref{Two Dimensional Statement of Two Dimensional Strictification}, to follow. We will build upon the work in \cite{Low dimensional structures formed by tricategories} to give an analogous tricategorical formulation of three-dimensional strictification in Proposition \ref{tricategory of tricategories}.

\begin{remark}\label{Two Dimensional Statement of Two Dimensional Strictification}
	\noindent There is a large $2$-category $\mathbf{Bicat}_\text{icon}$ whose objects are small bicategories, morphisms are pseudofunctors, and $2$-cells are invertible icons. The locally full sub-$2$-category whose objects are small $2$-categories and morphisms are $2$-functors is denoted $2$-$\mathbf{Cat}_\text{icon}$. Note that the composition of invertible icons in $\mathbf{Bicat}_\text{icon}$ is different to their usual composition as pseudonatural transformations, as described in \cite{Icons}. The inclusion $I_\text{icon}: 2$-$\mathbf{Cat}_\text{icon} \rightarrow \mathbf{Bicat}_\text{icon}$ has a left $2$-adjoint ${\left(\mathbf{st}_{2}\right)}_\text{icon}$. The right adjoint is biessentially surjective on objects, and the components of the unit are adjoint equivalences internal to $2$-$\mathbf{Cat}_\text{icon}$. On the other hand, the components of the counit are part of adjoint equivalences if and only if the $2$-category $\mathcal{A}$ is cofibrant in the Lack model structure \cite{Quillen 2-cat}. For our purposes, the most convenient characterisation of cofibrant $2$-categories is that their underlying categories are free on graphs. See \cite{Lack Homotopy Theoretic Aspects of 2-monads} for more on the connection with cofibrant replacement. The significance of this freeness in codimension one is recorded in Remark \ref{pseudofunctors out of cofibrant 2-categories}, to follow. It will be extended to the three-dimensional setting in Lemma \ref{Gurski Lemma semi strict} and Proposition \ref{trihomomorphism out of cofibrant Gray category}.
\end{remark}

\begin{remark}\label{pseudofunctors out of cofibrant 2-categories}
	Let $\left(F, \iota, \chi\right): \mathcal{A} \rightsquigarrow \mathcal{B}$ be a pseudofunctor between $2$-categories and suppose that the underlying category of $\mathcal{A}$ is free on a graph. Then there is a unique $2$-functor $\overline{F}: \mathcal{A} \rightarrow \mathcal{B}$ which agrees with $F$ on objects and generating morphisms, and an invertible icon $\overline{F} \rightarrow F$ with components described by the following induction \begin{itemize}
		\item The component on ${\left(-\right)}_{X}$ is given by the unitor $\iota_X: 1_{FX} \Rightarrow F1_{X}$,
		\item The component on a generating morphism $g: X \rightarrow Y$ is given by the identity,
		\item The component on a composite of non-identity morphisms \begin{tikzcd}
			X \arrow[r, "g"] & Y \arrow[r, "h"] & Z
		\end{tikzcd} in which $g$ is a generating morphism is given by the following pasting.
	\end{itemize}
	$$\begin{tikzcd}[row sep = 25, column sep = 30]
		FX \arrow[r, "\bar{F}g"]\arrow[d, "1_{FX}"']
		&FY \arrow[r, "\bar{F}h"] \arrow[d, "1_{FY}"]
		&FZ\arrow[d, "1_{FZ}"]
		\arrow[ld, Rightarrow, shorten = 5, "e_{h}"]
		\\
		FX \arrow[rr, bend right = 60, "F\left(hg\right)"'] \arrow[r, "Fg"']
		\arrow[r, equal, shorten = 5, shift left = 8]
		& FY \arrow[r, "Fh"']\arrow[d, Rightarrow, "\chi_{h{,}g}", shorten = 2]
		& FZ
		\\
		&{}
	\end{tikzcd}$$
	
	\noindent The $2$-functor $\overline{F}: \mathcal{A} \rightarrow \mathcal{B}$ is defined on $2$-cells $\phi: f \Rightarrow g$ via the composite depicted below.
	
	$$\begin{tikzcd}
		\overline{F}f \arrow[r, "e_{f}"] & Ff \arrow[r, "F\phi"] & Fg \arrow[r, "e_{g}^{-1}"] & \overline{F}g
	\end{tikzcd}$$ 
\end{remark}

\noindent The adjunction $\mathbf{st}_{2} \dashv I_{2}$ also underlies a triadjunction between $\mathbf{Bicat}$ and $\mathbf{Gray}$ (Theorem 3.1.3 of \cite{Campbell PhD}). This triadjunction takes into account arbitrary pseudonatural transformations, and modifications between them. It is much stricter than a typical triadjunction. In particular, it has isomorphisms, rather than just biequivalences, between hom $2$-categories. The main weak aspect of this triadjunction is the behaviour of its left triadjoint $\mathbf{st}_{2}$ between hom-bicategories. See \cite{Campbell Strictification} for an analysis of the strict aspects of the strictification triadjunction $\mathbf{st}_{2}\dashv I_{2}$. In Theorem \ref{trihomomorphism classifier}, the analogous `triequivalences between homs' will be described in the three dimensional setting. 

\begin{remark}\label{tri-fully-faithful and triessential image}
	There is a standard argument that an adjunction whose unit is an isomorphism has a fully-faithful left adjoint. This argument works equally well for triadjunctions with `internal biequivalence' in place of `isomorphism' and `tri-fully-faithful' in place of `fully-faithful'. Here `tri-fully-faithful' is taken to mean that the actions between hom-bicategories are biequivalences. We review how this argument applies to the strictification triadjunction. Firstly, the triangle identity says that the action of $\mathbf{st}_{2}$ between homs is given by the composite depicted below. But since $\eta_\mathcal{B}$ is always a biequivalence internal to $\mathbf{Bicat}$, this composite is also a biequivalence, and hence $\mathbf{st}_{2}$ is tri-fully-faithful. 
	
	$$\begin{tikzcd}
		\mathbf{Bicat}\left(\mathcal{A}{,} \mathcal{B}\right) \arrow[rrr, "\mathbf{Bicat}\left(\mathcal{A}{,} \eta_\mathcal{B}\right)"] &&& \mathbf{Bicat}\left(\mathcal{A}{,} I\mathbf{st}_{2}\left(\mathcal{B}\right)\right) \arrow[rrr, "\cong"] &&& \mathbf{Gray}\left(\mathbf{st}_{2}\left(\mathcal{A}\right){,}\mathbf{st}_{2}\left(\mathcal{B}\right)\right)
	\end{tikzcd}$$
	
	\noindent Moreover, if the underlying category of $\mathcal{A}$ is free on a graph then by Remark \ref{pseudofunctors out of cofibrant 2-categories}, the pseudofunctor $\eta_{I_{2}\mathcal{A}}$ is isomorphic to a strict $2$-functor via an invertible icon. In this case the counit of $\mathbf{st}_{2}\dashv I_{2}$ would be part of a biequivalence internal to $\mathbf{Gray}$. It therefore follows by a dual argument to the one above that $\mathbf{st}_{2}$ factors through the full-sub-tricategory $\mathbf{Gray}_c$ of $\mathbf{Gray}$ on cofibrant $2$-categories via a triequivalence \cite{Bicat Not Triequivalent To Gray}.
\end{remark}

\begin{remark}\label{dimension two strict, semi-strict, fully weak}
	The table below lists data of varying degrees of strictness in two-dimensions. 
	\\
	\begin{center}
		\begin{tabularx}{0.8\textwidth} { 
				| >{\raggedright\arraybackslash}X 
				| >{\raggedright\arraybackslash}X 
				| >{\centering\arraybackslash}X
				| >{\raggedleft\arraybackslash}X
				|}
			\hline
			Dimension & Strict & Semi-strict & Fully-weak \\
			\hline
			\hline
			$n = 0$  & $2$-categories & $2$-categories & bicategories  \\
			\hline
			$n = 1$ & $2$-functors & $2$-functors & pseudofunctors \\
			\hline
			$n = 2$ & $2$-natural transformations & pseudonatural transformations & pseudonatural transformations\\
			\hline
			$n = 3$ & modifications & modifications & modifications\\
			\hline
		\end{tabularx}
	\end{center}
	
	\noindent Observe that data appearing in the `strict' column are seen as enriched in the monoidal closed category $\left(\mathbf{Cat}, \times\right)$, and have all of their coherence data given by identities. Meanwhile, data in the `Fully weak' column can be seen as weakly enriched in the monoidal closed $2$-category $\left(\mathbf{Cat}, \times\right)$, in the sense of \cite{Garner Shulman Enriched Categories as a Free Cocompletion}. Finally, data in the `semi-strict' column may have non-identity \emph{free coherences} in the sense explained in Section \ref{subsection free versus derivable coherences}, but their \emph{operational coherences} are all identities. When strictification is extended from bicategories to $\left(2, k\right)$-transfors, fully weak transfors are mapped to semi-strict ones since $\mathbf{st}_{2}\left(\mathcal{B}\right)$ is free in codimension one.
\\
\\	
	\noindent The analogous table in dimension three is discussed in Remark \ref{dimension three strict, Gray, semi-strict, fully weak}. Once again, three-dimensional strictification will map fully weak $\left(3, k\right)$-transfors to ones which preserve operations in their domain $\mathbf{Gray}$-category on the nose, but only respect the varying of a $k$-cell in their domain along a $\left(k+1\right)$-cell up to an appropriate $\left(k+2\right)$-cell in their codomain $\mathbf{Gray}$-category.
\end{remark} 

\section{$\left(3, k\right)$-transfors}\label{Section Definitions of 3, k transfors}

\noindent We refer the reader to \cite{Gurski Coherence in Three Dimensional Category Theory} for definitions of tricategories and arbitrary $\left(3, k\right)$-transfors. Recall that categories enriched over $2$-$\mathbf{Cat}$ equipped with the $\mathbf{Gray}$-tensor product \cite{Gray Formal Category Theory} and corresponding closed structure $\mathbf{Gray}\left(\mathcal{A}, \mathcal{B}\right)$ model arbitrary tricategories up to weak equivalence \cite{GPS tricategory}. $\mathbf{Gray}$-categories satisfy the usual middle-four interchange law from two-dimensional category theory up to an invertible $3$-cell, but are otherwise strict. The aim of this section is to recall in detail the structure of the full-sub-$\mathbf{Gray}$-categories of $\mathbf{Tricat}\left(\mathfrak{A}, \mathfrak{B}\right)$ on $\mathbf{Gray}$-functors, when $\mathfrak{A}$ and $\mathfrak{B}$ are both $\mathbf{Gray}$-categories, and to distinguish between operational and free coherences of $\left(3, k\right)$-transfors.

\subsection{Trihomomorphisms between tricategories}

\begin{definition}\label{Definition trihomomorphism}
	Let $\mathfrak{A}$ and $\mathfrak{B}$ be tricategories. A \emph{trihomomorphism} $F: \mathfrak{A} \rightsquigarrow \mathfrak{B}$ consists of the following data. Note that to save space we have written, for example, $\mathfrak{B}^{3}$ in place of the product of hom-bicategories depicted below.
	
	$$\mathfrak{B}\left(FY, FZ\right)\times \mathfrak{B}\left(FX, FY\right)\times \mathfrak{B}\left(FW, FX\right)$$
	
	DATA
	\\
	\begin{itemize}
		\item A function $F_{0}: \mathfrak{A}_{0} \rightarrow \mathfrak{B}_{0}$, and a family of pseudofunctors $F_{X, Y}: \mathfrak{A}\left(X, Y\right) \rightsquigarrow \mathfrak{B}\left(X, Y\right)$.
		\item A family of pseudonatural adjoint equivalences with left adjoints $\iota^{X}: I_{FX} \rightarrow \left(F_{X, X}\right).\left(I_{X}\right)$, indexed by the objects of $\mathfrak{A}$. These will be called the \emph{unitors} of $F: \mathfrak{A} \rightsquigarrow \mathfrak{B}$.
		\item A family of pseudonatural adjoint equivalences with left adjoints as depicted below, indexed by triples of objects in $\mathfrak{A}$. These are called the \emph{compositors} of $F: \mathfrak{A} \rightsquigarrow \mathfrak{B}$.
		
		$$\chi^{X, Y, Z}: \left(\circ_{FX, FY, FZ}\right).\left(F_{Y, Z}\times F_{X, Y}\right) \rightarrow \left(F_{X, Z}\right).\circ_{X, Y, Z}$$
		
		\item A family of invertible modifications $\omega^{W, X, Y, Z}$ as depicted below, indexed by quadruplets of objects in $\mathfrak{A}$. These are called the \emph{associators} of $F: \mathfrak{A} \rightsquigarrow \mathfrak{B}$.
	\end{itemize}
	\noindent\begin{tikzcd}[column sep=18,  font=\fontsize{9}{6}]
		&
		\mathfrak{A}^3
		\arrow[ldd,"\circ \times 1"']
		\arrow[rdd,"1 \times \circ"]
		\arrow[rrr,"F^3"]
		&&
		{}
		\arrow[dd,Rightarrow,shorten=15,"1 \times \chi"]
		&
		\mathfrak{B}^3
		\arrow[rdd,"1 \times \circ"]
		&&&&
		\mathfrak{A}^3
		\arrow[ldd,"\circ \times 1"']
		\arrow[rrr,"F^3"]
		&
		{}
		\arrow[dd,Rightarrow,shorten=15,"\chi \times 1"]
		&&
		\mathfrak{B}^3
		\arrow[ldd,"\circ \times 1"']
		\arrow[rdd,"1 \times \circ"]
		\\
		\\
		\mathfrak{A}^2
		\arrow[rdd,"\circ"']
		&&
		\mathfrak{A}^2
		\arrow[ldd,"\circ"]
		\arrow[rrr,"F^2"]
		\arrow[ll,Rightarrow,shorten=15,"\alpha"']
		&
		{}
		\arrow[dd,Rightarrow,shorten=15,"\chi"]
		&&
		\mathfrak{B}^2
		\arrow[rr,Rightarrow,shorten=5,"\omega"]
		\arrow[rr,shorten=5, no head]
		\arrow[ldd,"\circ"]
		&&
		\mathfrak{A}^2
		\arrow[rdd,"\circ"']
		\arrow[rrr,"F"']
		&&
		{}
		\arrow[dd,Rightarrow,shorten=15,"\chi"]
		&
		\mathfrak{B}^2
		\arrow[rdd,"\circ"']
		&&
		\mathfrak{B}^2
		\arrow[ldd,"\circ"]
		\arrow[ll,Rightarrow,"\alpha"',shorten=15]
		\\
		\\
		&
		\mathfrak{A}
		\arrow[rrr,"F"']
		&&
		{}
		&
		\mathfrak{B}
		&&&&
		\mathfrak{A}
		\arrow[rrr,"F"']
		&
		{}
		&&
		\mathfrak{B}
	\end{tikzcd}
	\begin{itemize}
		\item A family of invertible modifications $\gamma^{X, Y}$, indexed by pairs of objects in $\mathfrak{A}$. The source of $\gamma^{X, Y}$ is the pseudonatural transformation depicted below left, while the target is $\lambda.F$. These invertible modifications are called the \emph{left unitors} of $F: \mathfrak{A} \rightsquigarrow \mathfrak{B}$.
		\item A family of invertible modifications $\delta^{X, Y}$, indexed by pairs of objects in $\mathfrak{A}$. The source of $\delta^{X, Y}$ is the pseudonatural transformation depicted below right, while the target is $\rho.F$. These invertible modifications are called the \emph{right unitors} of $F: \mathfrak{A} \rightsquigarrow \mathfrak{B}$.
	\end{itemize}
	
	\begin{tikzcd}[ font=\fontsize{9}{6}]
		1 \times \mathfrak{A}
		\arrow[rr,"I \times 1"]
		\arrow[rrrr,bend left=40,"\left(I{,}F\right)"{name=A}]
		\arrow[rrdd,bend right=15,"1"'{name=B}]
		&
		{}
		\arrow[to=B,Rightarrow,shorten=5,"\lambda", near start]
		&
		\mathfrak{A}^2
		\arrow[rr,"F^2"]
		\arrow[dd,"\circ"']
		\arrow[to=A,Leftarrow,shorten=5,"\iota"',near start]
		&
		{}
		\arrow[dd,Rightarrow,shorten=15,"\chi"]
		&
		\mathfrak{B}^2
		\arrow[dd,"\circ"]
		\\
		\\
		&&
		\mathfrak{A}
		\arrow[rr,"F"']
		&
		{}
		&
		\mathfrak{B}
	\end{tikzcd}\begin{tikzcd}[ font=\fontsize{9}{6}]
		\mathfrak{A} \times 1
		\arrow[rr,"1 \times I"]
		\arrow[rrrr,bend left=40,"\left(F{,}I\right)"{name=A}]
		\arrow[rrdd,bend right=15,"1"'{name=B}]
		&
		{}
		\arrow[to=B,Rightarrow,shorten=5,"\rho", near start]
		&
		\mathfrak{A}^2
		\arrow[rr,"F^2"]
		\arrow[dd,"\circ"']
		\arrow[to=A,Leftarrow,shorten=5,"\iota"',near start]
		&
		{}
		\arrow[dd,Rightarrow,shorten=15,"\chi"]
		&
		\mathfrak{B}^2
		\arrow[dd,"\circ"]
		\\
		\\
		&&
		\mathfrak{A}
		\arrow[rr,"F"']
		&
		{}
		&
		\mathfrak{B}
	\end{tikzcd}

\noindent We refer the reader to Definition 4.10 of \cite{Gurski Coherence in Three Dimensional Category Theory} for the axioms that these data are required to satisfy.
\end{definition}

\noindent Definition \ref{Trihomomorphism special cases} specifies a special kind of trihomomorphisms which we will use in Section \ref{Path Objects Section} to study $\left(3, k\right)$-transfors via path objects of $\mathbf{Gray}$-categories.

\begin{definition}\label{Trihomomorphism special cases}
	A trihomomorphism $F:\mathfrak{A} \rightsquigarrow \mathfrak{B}$ will be called a \emph{$3$-pseudofunctor} if $\mathfrak{A}$ and $\mathfrak{B}$ are $\mathbf{Gray}$-categories, all $F_{X, Y}$ are $2$-functors, both $\iota$ and $\chi$ are invertible $2$-natural transformations, and all thee modifications $\gamma$, $\delta$ and $\omega$ are identities.
\end{definition}

\subsection{Definitions of $\left(3, k\right)$-transfors between $\mathbf{Gray}$-functors}\label{subsection weak higher maps between Gray functors}

\noindent We refer the reader to Section 4.3 of \cite{Gurski Coherence in Three Dimensional Category Theory} for the general definitions of the various higher dimensional maps between tricategories. However, the data and axioms involved in describing the higher cells simplify considerably if the tricategories being considered are $\mathbf{Gray}$-categories and if the trihomomorphisms being considered are $\mathbf{Gray}$-functors. The goal of this subsection is to recall these definitions in this special setting.

\begin{definition}\label{definition trinatural transformation}
	Let $F, G: \mathfrak{A} \rightarrow \mathfrak{B}$ be $\mathbf{Gray}$-functors. A \emph{trinatural transformation} $p: F \Rightarrow G$ consists of \begin{itemize}
		\item For every $X \in \mathfrak{A}$, an arrow $p_{X}: FX \rightarrow GX$ in $\mathfrak{B}$,
		\item For every arrow $f: X \rightarrow Y$ in $\mathfrak{A}$, an adjoint equivalence $p_{f} \dashv p_{f}^{*}$ in $\mathfrak{B}\left(FX, GY\right)$ with left adjoint $p_{f}: p_{Y}.Ff \Rightarrow Gf.p_{X}$.
		\item For every $2$-cell $\phi: f \Rightarrow g: X \rightarrow Y$ in $\mathfrak{A}$, an invertible $3$-cell $p_\phi$ called the \emph{local pseudonaturality constraint}, as depicted below left,
		\item For every object $X \in \mathfrak{A}$, an invertible $3$-cell $p^X$ called the \emph{unitor}, as depicted below centre,
		\item For every \begin{tikzcd}
			X \arrow[r, "f"] & Y \arrow[r, "g"] & Z
		\end{tikzcd} in $\mathfrak{A}$, an invertible $3$-cell $p_{g, f}$ called the \emph{compositor}, as depicted below right.
	\end{itemize}
	
	$$\begin{tikzcd}[column sep = 12, font=\fontsize{9}{6}]
		p_Y.Ff
		\arrow[rr,"1.F\phi"]
		\arrow[dd,"p_f"']
		&
		{}
		\arrow[dd,Rightarrow,shorten=15,"p_\phi"]
		&
		p_Y.Fg
		\arrow[dd,"p_g"]
		\\
		\\
		Gf.p_X
		\arrow[rr,"G\phi.1"']
		&
		{}
		&
		Gg.p_X &{}
	\end{tikzcd}\begin{tikzcd}[column sep = 4, font=\fontsize{9}{6}]
		p_X.F1_X
		\arrow[rr,equal]
		\arrow[dd,"p_{1_X}"']
		&
		{}
		\arrow[dd,Rightarrow,shorten=15,"p^X"]
		&
		p_X.1_{FX}
		\arrow[dd,"1_{p_X}"]
		\\
		\\
		G1_X.p_X
		\arrow[rr, equal]
		&
		{}
		&
		{1_{GX}}.p_X&{}
	\end{tikzcd}\begin{tikzcd}[column sep = 1, font=\fontsize{9}{6}]
		p_Z.Fg.Ff
		\arrow[rr,"p_{gf}"]
		\arrow[rdd,"p_g.1"']
		&
		{}
		&
		Gg.Gf.p_X
		\\
		\\
		&
		Gg.p_Y.Ff
		\arrow[ruu,"1.p_f"']
		\arrow[uu,Rightarrow,shorten=10,"p_{g{,}f}"']
	\end{tikzcd}$$
	
	\noindent These data are subject to the following axioms.
	
	\begin{itemize}
		\item The assignation $\left(f, \phi\right) \mapsto \left(p_{f}, p_{\phi}\right)$ is a pseudonatural transformation. For $\left(f, \phi\right) \mapsto \left(p_{f}, p_{\phi}\right)$ this says that $p_{1_{f}} = 1_{p_{f}}$ and that the following equations hold for every $\Omega: \phi \Rrightarrow \phi'$ and every \begin{tikzcd}
			f \arrow[r, Rightarrow, "\phi"] & g \arrow[r, Rightarrow, "\psi"] & h
		\end{tikzcd} in $\mathfrak{A}$.

		$$\begin{tikzcd}[column sep = 18, row sep = 20, font=\fontsize{9}{6}]
			p_Y.Ff
			\arrow[rr,bend left = 45, "1.F\phi"]
			\arrow[dd,"p_f"']
			&
			{}
			\arrow[d,Rightarrow, shorten = 5,"p_\phi"]
			&
			p_Y.Fg
			\arrow[dd,"p_g"]
			\\
			&{}&&=
			\\
			Gf.p_X
			\arrow[rr,bend left = 45,"G\phi.1" {name = C}]
			\arrow[rr, bend right = 45, "G\phi'.1"' {name = D}]
			&
			{}
			&
			Gg.p_X
			\arrow[from =C, to =D, Rightarrow, shorten = 10, shift right = 5, "G\Omega .1"]
		\end{tikzcd}	\begin{tikzcd}[column sep = 18, row sep = 20, font=\fontsize{9}{6}]
			p_Y.Ff
			\arrow[rr,bend left = 45, "1.F\phi" {name = A}]
			\arrow[rr, bend right = 45, "1.F\phi'"' {name = B}]
			\arrow[dd,"p_f"']
			&
			&
			p_Y.Fg
			\arrow[dd,"p_g"]
			\\
			&{}
			\arrow[d,Rightarrow, shorten = 5,"p_\phi"]
			\\
			Gf.p_X\
			\arrow[rr, bend right = 45, "G\phi'.1"']
			&
			{}
			&
			Gg.p_X &{}
			\arrow[from =A, to =B, Rightarrow, shorten = 10, shift right = 5, "1.F\Omega"]
		\end{tikzcd}$$ 
		
		$$\begin{tikzcd}[column sep = 18, font=\fontsize{9}{6}]
			p_Y.Ff
			\arrow[rr,"1.F\phi"]
			\arrow[dd,"p_f"']
			&
			{}
			\arrow[dd,Rightarrow,shorten=15,"p_\phi"]
			&
			p_Y.Fg
			\arrow[dd,"p_g"] \arrow[rr, "1.F\psi"]&{}\arrow[dd, Rightarrow, shorten = 15, "p_{\psi}"]& p_{Y}.Fh \arrow[dd, "p_{h}"]
			\\
			&&&&&=
			\\
			Gf.p_X
			\arrow[rr,"G\phi.1"']
			&
			{}
			&
			Gg.p_X\arrow[rr, "G\psi.1"'] &{}&Gh.p_{X}
		\end{tikzcd}
		\begin{tikzcd}[column sep = 18, row sep = 24, font=\fontsize{9}{6}]
			p_Y.Ff
			\arrow[rr,"1.F\left(\psi\phi\right)"]
			\arrow[dd,"p_f"']
			&
			{}
			\arrow[dd,Rightarrow,shorten=15,"p_{\psi\phi}"]
			&
			p_Y.Fh
			\arrow[dd,"p_h"]
			\\
			\\
			Gf.p_X
			\arrow[rr,"G\left(\psi.\phi\right).1"']
			&
			{}
			&Gh.p_{X}
		\end{tikzcd}$$
		
		\item The assignation $\left(g, f\right) \mapsto p_{g, f}$ given in the compositor is a modification. This says that \begin{itemize}
			\item For every $\phi: f \Rightarrow f': X \rightarrow Y$  and $g: Y \rightarrow Z$ the first equation depicted below holds,
			\item For every $f: X \rightarrow Y$ and $\psi: g \Rightarrow g': Y \Rightarrow Z$ the second equation depicted below holds.
		\end{itemize} 
		We will respectively refer to these conditions as the \emph{right whiskering law} and \emph{left whiskering law} for $p$. They are specified by $2$-cells in $\mathfrak{A}$ being whiskered by $1$-cells on the left or on the right respectively.
	\end{itemize}

	$$\begin{tikzcd}[row sep = 20, column sep = 15, font=\fontsize{9}{6}]
		&Gg.p_{Y}.Ff\arrow[rd, "1.p_{f}"]
		\arrow[d, Rightarrow, shorten = 5, "p_{g{,}f}"]
		\\
		p_{Z}.Fg.Ff
		\arrow[ddd, "1.1_{Fg}F\phi" description]
		\arrow[ru, "p_{g}.1"]
		\arrow[rr, "p_{gf}"']
		&{}\arrow[ddd, Rightarrow, shorten = 18, "p_{g\phi}"]
		& Gg.Gf.p_{X}\arrow[ddd, "1.G\phi.1" description]
		\\
		&&&=
		\\
		\\
		p_{Z}.Fg.Ff' \arrow[rr, "p_{gf'}"']
		&{}
		& Gg.Gf'.p_{X}
	\end{tikzcd}\begin{tikzcd}[column sep = 15, font=\fontsize{9}{6}]
		&Gg.p_{Y}.Ff\arrow[rd, "1.p_{f}"]
		\arrow[ddd, shorten = 10, "1.1.F\phi" description]
		\\
		p_{Z}.Fg.Ff
		\arrow[dd, Rightarrow, shorten = 10, "{\left(p_{g}\right)}_{\left(F\phi\right)}", shift left = 10]
		\arrow[ddd, "1.1.F\phi"']
		\arrow[ru, "p_{g}.1"]
		&& Gg.Gf.p_{X}
		\arrow[dd, Rightarrow, shorten = 10, "1.p_{\phi}"', shift right = 10]
		\arrow[ddd, "1.G\phi.1"]
		\\
		\\
		{}&Gg.p_{Y}.Ff'
		\arrow[rd, "1.p_{f'}"]
		\arrow[d, Rightarrow, shorten = 5, "p_{g{,}f'}"]
		&{}
		\\
		p_{Z}.Fg.Ff'
		\arrow[ru, "p_{g}.1"]
		\arrow[rr, "p_{gf'}"']
		&{}
		& Gg.Gf'.p_{X}
	\end{tikzcd}$$
	\\
	$$\begin{tikzcd}[row sep = 20, column sep = 15, font=\fontsize{9}{6}]
		&Gg.p_{Y}.Ff\arrow[rd, "1.p_{f}"]
		\arrow[d, Rightarrow, shorten = 5, "p_{g{,}f}"]
		\\
		p_{Z}.Fg.Ff
		\arrow[ddd, "1.F\psi.1_{Ff}" description]
		\arrow[ru, "p_{g}.1"]
		\arrow[rr, "p_{gf}"']
		&{}\arrow[ddd, Rightarrow, shorten = 18, "p_{\psi f}"]
		& Gg.Gf.p_{X}\arrow[ddd, "G\psi.1.1" description]
		\\
		&&&=
		\\
		\\
		p_{Z}.Fg'.Ff \arrow[rr, "p_{g'f}"']
		&{}
		& Gg'.Gf.p_{X}
	\end{tikzcd}\begin{tikzcd}[column sep = 15, font=\fontsize{9}{6}]
		&Gg.p_{Y}.Ff\arrow[rd, "1.p_{f}"]
		\arrow[ddd, shorten = 10, "G\psi.1.1" description]
		\\
		p_{Z}.Fg.Ff
		\arrow[dd, Rightarrow, shorten = 10, "p_{\psi}.1", shift left = 10]
		\arrow[ddd, "1_{p_{Z}}.F\psi.1"']
		\arrow[ru, "p_{g}.1"]
		&& Gg.Gf.p_{X}
		\arrow[dd, Rightarrow, shorten = 10, "{\left(G\psi\right)}_{\left(p_{f}\right)}"', shift right = 10]
		\arrow[ddd, "G\psi.1.1"]
		\\
		\\
		{}&Gg'.p_{Y}.Ff
		\arrow[rd, "1.p_{f}"]
		\arrow[d, Rightarrow, shorten = 5, "p_{g'{,}f}"]
		&{}
		\\
		p_{Z}.Fg'.Ff
		\arrow[ru, "p_{g'}.1"]
		\arrow[rr, "p_{g'f}"']
		&{}
		& Gg'.Gf.p_{X}
	\end{tikzcd}$$
	\begin{itemize}
		\item For every composable triple \begin{tikzcd}
			W \arrow[r, "e"] & X \arrow[r, "f"] &Y \arrow[r, "g"] & Z 
		\end{tikzcd} in $\mathfrak{A}$, the following equation, called the \emph{associativity coherence} holds in the hom-$2$-category $\mathfrak{B}\left(W, Z\right)$. 
		\\
		\\
		\begin{tikzcd}[font=\fontsize{9}{6}]
			&
			Gg.p_Y.Ff.Fe
			\arrow[rr,"1.p_f.1"]
			\arrow[rrrdd,"1.p_{fe}"']
			\arrow[dd,Rightarrow,shorten=15,"p_{{g}{,}{fe}}"]
			&&
			Gg.Gf.p_X.Fe
			\arrow[rdd,"1.1.p_e"]
			\arrow[d,Rightarrow,shorten=4,"1.p_{{f}{,}{e}}"']
			\\
			&&&
			{}
			\\
			p_Z.Fg.Ff.Fe
			\arrow[ruu,"p_g.1.1"]
			\arrow[rrrr,"p_{gfe}"']
			&
			{}
			&&&
			Gg.Gf.Ge.p_W
		\end{tikzcd}=
		\\
		\\
		\\
		\begin{tikzcd}[font=\fontsize{9}{6}]
			&
			Gg.p_Y.Ff.Fe
			\arrow[rr,"1.p_f.1"]
			\arrow[d,Rightarrow,shorten=4,"p_{{g}{,}{f}}.1"]
			&&
			Gg.Gf.p_X.Fe
			\arrow[rdd,"1.1.p_e"]
			\arrow[dd,Rightarrow,shorten=15,"p_{{gf}{,}{e}}"']
			\\
			&
			{}
			\\
			p_Z.Fg.Ff.Fe
			\arrow[ruu,"p_g.1.1"]
			\arrow[rrrr,"p_{gfe}"']
			\arrow[rrruu,"p_{gf}.1"']
			&&&
			{}
			&
			Gg.Gf.Ge.p_W
		\end{tikzcd}
		\item For every $f: X \rightarrow Y$, the following pastings of $3$-cells in $\mathfrak{B}$ are both equal to the identity on $p_{f}$. These equations are respectively called the \emph{left and right unit laws}.
	\end{itemize}
	\noindent 
	\begin{tikzcd}[column sep = 15, font=\fontsize{9}{6}]
		&
		Gf.p_X.F1_X
		\arrow[rdd,bend right,"{1_{Gf}}.{p_{1_X}}"'{name=B}]
		\arrow[rdd,bend left=40,"{1_{Gf}}.{1_{p_X}}"{name=A}]
		\arrow[dd,Rightarrow,shorten=8,shift right=8,"p_{{f}{,}{1_X}}"']
		\\
		\\
		p_Y.Ff.F1_X
		\arrow[ruu,"{p_f}.{1_{F1_X}}",bend left=30]
		\arrow[rr,"p_f"']
		&
		{}
		&
		Gf.G1_X.p_X
		\arrow[from=A,to=B,Rightarrow,shorten=10,"1.{p^X}"',shift left=2]
		&
		{}
	\end{tikzcd}  \begin{tikzcd}[column sep = 15, font=\fontsize{9}{6}]
		&
		G1_Y.p_Y.Ff
		\arrow[rdd,bend left=30,"{1_{G1_Y}}.{p_f}"]
		\arrow[dd,Rightarrow,shorten=8,"p_{{1_Y}{,}{f}}",shift left=8]
		\\
		\\
		p_Y.F1_Y.Ff
		\arrow[ruu,bend left=40,"{1_{p_Y}}.{1_{Ff}}"{name=A}]
		\arrow[ruu,bend right,"{p_{1_Y}}.{1_{Ff}}"'{name=B}]
		\arrow[rr,"p_f"']
		&
		{}
		&
		G1_Y.Gf.p_X
		\arrow[from=A,to=B,Rightarrow,shorten=10,"{p^Y}.1",shift right=2]
	\end{tikzcd}
	
\end{definition}

\begin{remark}
	The unit and counit of the adjoint equivalence $p_{f} \dashv p_{f}^{*}$ will be denoted $p_{\eta}^{f}$ and $p_{\varepsilon}^f$ respectively. There are also $3$-cell components $p_{\phi}^{*}$ depicted below left which are the mates of the $3$-cell components $p_{\phi}$. Finally, there are modification conditions for $p_{f}^{\eta}$ and $p_{f}^{\varepsilon}$, with the equation for $p_{f}^{\eta}$ depicted below. These modification conditions are derivable using the triangle identities of the adjunction $p_{f} \dashv p_{f}^{*}$, and the definition of $p_{\phi}^{*}$ via mates.  
	\\

	\noindent \begin{tikzcd}[column sep = 12, font=\fontsize{9}{6}]
		p_Y.Ff
		\arrow[rr,"1.F\phi"]
		\arrow[dd, leftarrow, "p_{f}^{*}"']
		&
		{}
		\arrow[dd,Rightarrow,shorten=15,"p_{\phi}^{*}"]
		&
		p_Y.Fg
		\arrow[dd,"p_{g}^{*}", leftarrow]
		\\
		\\
		Gf.p_X
		\arrow[rr,"G\phi.1"']
		&
		{}
		&
		Gg.p_X &{}
	\end{tikzcd} \begin{tikzcd}[row sep = 20, column sep = 15, font=\fontsize{9}{6}]
		&Gf.p_{X}\arrow[rd, "p^*_f"]
		\arrow[d, Leftarrow, shorten = 5, "p^{\eta}_f"]
		\\
		p_{Y}.Ff
		\arrow[ddd, "1.F\phi" description]
		\arrow[ru, "p_{f}"]
		\arrow[rr, "1"']
		&{}
		& p_Y .Ff\arrow[ddd, "1.F\phi" description]
		\\
		&=&&=
		\\
		\\
		p_{Y}.Fg \arrow[rr, "1"']
		&{}
		& p_{Y}.Fg
	\end{tikzcd}\begin{tikzcd}[column sep = 15, font=\fontsize{9}{6}]
		&Gf.p_{X}\arrow[rd, "p_{f}^{*}"]
		\arrow[ddd, "G\phi .1" description]
		\\
		p_{Y}.Ff
		\arrow[dd, Leftarrow, shorten = 10, "p_{\phi}", shift left = 10]
		\arrow[ddd, "1.F\phi"']
		\arrow[ru, "p_{f}"]
		&& p_{Y}.Ff
		\arrow[dd, Leftarrow, shorten = 10, "p_{\phi}^{*}"', shift right = 10]
		\arrow[ddd, "1.F\phi"]
		\\
		\\
		{}&Gg.p_{X}
		\arrow[rd, "p_{g}^{*}"]
		\arrow[d, Leftarrow, shorten = 5, "p_{g}^{\eta}"]
		&{}
		\\
		p_{Y}.Fg
		\arrow[ru, "p_{g}"]
		\arrow[rr, "1"']
		&{}
		& p_{Y}.Fg
	\end{tikzcd}
\end{remark}

\noindent It will be useful to have terminology for trinatural transformations of different degrees of strictness. This terminology will be explained in Remark \ref{terminology for higher cells between tricats}.

\begin{definition}\label{degrees of strictness trinatural transformations}
	Call a trinatural transformation $p: F \rightarrow G$
	
	\begin{itemize}
		\item A \emph{pseudo-icon equivalence} if $F$ and $G$ agree on objects and its $1$-cell components are identities.
		\item An \emph{invertible ico-icon} if it is a pseudo-icon equivalence, $F$ and $G$ agree on arrows and the adjoint equivalences in $\mathfrak{B}\left(FX, GY\right)$ are identities.
		\item \emph{locally strict} if its component pseudonatural equivalences $\left(f, \phi\right) \mapsto \left(p_{f}, p_{\phi}\right)$ are $2$-natural isomorphisms, hence if each $p_{\phi}$, $p_{\phi}^{*}$, $p_{f}^{\eta}$ and $p_{f}^{\varepsilon}$ are identities.
		\item \emph{strict} if all components are identities, except perhaps $p_{X}$ for objects $X \in \mathfrak{A}$.
		\item \emph{unital} if its unitors are identities.
		\item \emph{compositional} if its compositors are identities.
		\item \emph{semi-strict} if it is unital and compositional.
		\item A \emph{$3$-pseudonatural transformation} if it is semi-strict and locally strict.
	\end{itemize}
\end{definition}

\noindent The reader may have noticed that the compositors and unitors of a trinatural transformation, and the axioms they must satisfy, resemble the data and axioms of a pseudofunctor. We will make this analogy precise in Proposition \ref{trinatural transformations as 2-Cat pseudofunctors}.

\begin{definition}\label{trimodification definition}
	Let $F, G: \mathfrak{A} \rightarrow \mathfrak{B}$ be $\mathbf{Gray}$-functors and let $p, q: F \rightarrow G$ be trinatural transformations. A \emph{trimodification} $\sigma: p \rightarrow q$ consists of
	
	\begin{itemize}
		\item For every $X \in \mathfrak{A}$, a $2$-cell $\sigma_{X}: p_{X} \Rightarrow q_{X}$
		\item For every $f: X \rightarrow Y$, an invertible $3$-cell \begin{tikzcd}
			p_{Y}.Ff
			\arrow[rr, "p_{f}"]
			\arrow[dd, "\sigma_{Y}.1_{Ff}"']
			& {}\arrow[dd, Rightarrow, shorten = 10, "\sigma_{f}"]
			&
			Gf.p_{X}
			\arrow[dd, "1_{Gf}.\sigma_{X}"]
			\\
			\\
			q_{Y}.Ff \arrow[rr, "q_{f}"']
			&{}&
			Gf.q_{X} 
		\end{tikzcd} 
	\end{itemize}
	\noindent These data are subject to the following axioms.
	
	\begin{itemize}
		\item For every $X, Y \in \mathfrak{A}$, the assignation $f \mapsto \sigma_{f}$ defines an invertible modification, from the pseudonatural transformation below left to the pseudonatural transformation below right. Here ${\left(\Phi\right)}^{*}$ and ${\left(\Phi\right)}_{*}$ denote restriction and extension along $\Phi$ respectively. 
		
		$$\begin{tikzcd}[column sep = 10, row sep = 15, font=\fontsize{9}{6}]
			&&
			\mathfrak{A}\left(X{,}Y\right)
			\arrow[rrdd,"G"]
			\arrow[lldd, "F"']
			\\
			\\
			\mathfrak{B}\left(FX{,}FY\right)
			\arrow[rrdd,bend right,"{p_{Y}}_{*}"']
			&&&&
			\mathfrak{B}\left(GX{,} GY\right)
			\arrow[lldd,"p_{X}^{*}"'name=A]
			\arrow[lldd,bend left,"q_{X}^{*}" name=B]
			\arrow[llll,shorten=25,Leftarrow,"{p}^{X{,}Y}"]
			\\
			\\
			&&
			\mathfrak{B}\left(FX{,}GY\right)
			\arrow[from=A,to=B,"\sigma_{X}^{*}",Rightarrow,shorten=8]
		\end{tikzcd}\begin{tikzcd}[column sep = 10, row sep = 15, font=\fontsize{9}{6}]
			&&
			\mathfrak{A}\left(X{,}Y\right)
			\arrow[rrdd,"G"]
			\arrow[lldd, "F"']
			\\
			\\
			\mathfrak{B}\left(FX{,}FY\right)
			\arrow[rrdd, "{q_{Y}}_{*}"name = A]
			\arrow[rrdd,bend right,"{p_{Y}}_{*}"'name=B]
			&&&&
			\mathfrak{B}\left(GX{,} GY\right)
			\arrow[lldd,bend left,"q_{X}^{*}"]
			\arrow[llll,shorten=25,Leftarrow,"{p}^{X{,}Y}"]
			\\
			\\
			&&
			\mathfrak{B}\left(FX{,}GY\right)
			\arrow[from=A,to=B,"{\sigma_{Y}}_{*}",Leftarrow,shorten=8]
		\end{tikzcd}$$
		
		\noindent  This says that for every $\phi: f \Rightarrow g: X \rightarrow Y$ in $\mathfrak{A}$, the following equation holds in $\mathfrak{B}$. We will refer to this as the \emph{local modification condition} for $\sigma$.
	\end{itemize}
	
	\noindent\begin{tikzcd}[column sep = 12, font=\fontsize{9}{6}]
		&&p_{Y}.Ff
		\arrow[lldd, "1.F\phi"']
		\arrow[rrdd, "\sigma_{Y}.1" description]
		\arrow[dddd, Rightarrow, "{\left(\sigma_{Y}\right)}_{\left(F\phi\right)}", shorten = 30, shift right = 6]
		\arrow[rr, "p_f"]
		&&Gf.p_{X}
		\arrow[dd, Rightarrow, "\sigma_{f}", shorten = 10]
		\arrow[rrdd, "1.\sigma_{X}"]
		\\
		\\
		p_{Y}.Fg\arrow[rrdd, "\sigma_{Y}.1"']
		&&&&q_{Y}.Ff
		\arrow[lldd, "1.F\phi" description]
		\arrow[rr, "q_{f}"]
		\arrow[dd, Rightarrow, "q_{\phi}", shorten = 10]
		&&Gf.q_{X}\arrow[lldd, "G\phi.1_{q_{X}}"]&=
		\\
		\\
		&&q_{Y}.Fg
		\arrow[rr, "q_{g}"']
		&&Gg.q_{X}
	\end{tikzcd}\begin{tikzcd}[column sep = 12, font=\fontsize{9}{6}]
		&&p_{Y}.Ff
		\arrow[lldd, "1.F\phi"']
		\arrow[dd, Rightarrow, "p_{\phi}"', shorten = 10]
		\arrow[rr, "p_f"]
		&&Gf.p_{X}
		\arrow[lldd, "G\phi.1" description]
		\arrow[dddd, Rightarrow, "{\left(G\phi\right)}_{\left(\sigma_{X}\right)}"', shorten = 30, shift left = 6]
		\arrow[rrdd, "1.\sigma_{X}"]
		\\
		\\
		p_{Y}.Fg
		\arrow[rrdd, "\sigma_{Y}.1"']
		\arrow[rr, "p_{g}"]
		&&Gg.p_{X}
		\arrow[dd, Rightarrow, "\sigma_{g}"', shorten = 10]
		\arrow[rrdd, "1.\sigma_{X}" description]
		&&
		&&Gg.p_{X}
		\arrow[lldd, "G\phi.1"]
		\\
		\\
		&&Gf.q_{X}
		\arrow[rr, "q_{g}"']
		&&Gg.q_{X}
	\end{tikzcd}
	\begin{itemize}
		\item For every $X \in \mathfrak{A}$, the following equation, called the \emph{unit law}, holds in $\mathfrak{B}$.

		\begin{tikzcd}[column sep = 18, row sep = 15, font=\fontsize{9}{6}]
			p_Y.F1_{X}
			\arrow[rr,bend left = 45, "1_{p_X}"]
			\arrow[dd,"\sigma_{X}.1"']
			&
			=
			&
			G1_{X}.p_{X}
			\arrow[dd,"1.\sigma_{X}"]
			\\
			&{}&&=
			\\
			q_{X}.F1_{X}
			\arrow[rr,bend left = 45,"1_{q_X}" {name = C}]
			\arrow[rr, bend right = 45, "q_{1_{X}}"' {name = D}]
			&
			{}
			&
			G1_{X}.q_X
			\arrow[from =C, to =D, Rightarrow, shorten = 10, "q^{X}"]
		\end{tikzcd}	\begin{tikzcd}[column sep = 18, row sep = 15, font=\fontsize{9}{6}]
			p_X.F1_{X}
			\arrow[rr,bend left = 45, "1_{p_X}" {name = A}]
			\arrow[rr, bend right = 45, "p_{1_{X}}"' {name = B}]
			\arrow[dd,"\sigma_{X}.1_{F1_{X}}"']
			&
			&
			G1_{X}.p_{X}
			\arrow[dd,"1_{G1_{X}}.\sigma_{X}"]
			\\
			&{}
			\arrow[d,Rightarrow, shorten = 5,"\sigma_{1_{X}}"]
			\\
			q_{X}.F1_{X}\
			\arrow[rr, bend right = 45, "q_{1_{X}}"']
			&
			{}
			&
			G1_{X}.q_X &{}
			\arrow[from =A, to =B, Rightarrow, shorten = 10, "p^X"]
		\end{tikzcd} 
		
		\item For every \begin{tikzcd}
			X \arrow[r, "f"] & Y \arrow[r, "g"] & Z
		\end{tikzcd} in $\mathfrak{A}$, the following equation, called the \emph{composition law}, holds in $\mathfrak{B}$.
		
		\begin{tikzcd}[row sep = 20, column sep = 15, font=\fontsize{9}{6}]
			&Gg.p_{Y}.Ff\arrow[rd, "1.p_{f}"]
			\arrow[d, Rightarrow, shorten = 5, "p_{g{,}f}"]
			\\
			p_{Z}.Fg.Ff
			\arrow[ddd, "\sigma_{Z}.1.1" description]
			\arrow[ru, "p_{g}.1"]
			\arrow[rr, "p_{gf}"']
			&{}\arrow[ddd, Rightarrow, shorten = 18, "\sigma_{gf}"]
			& Gg.Gf.p_{X}\arrow[ddd, "1.1.\sigma_{X}" description]
			\\
			&&&=
			\\
			\\
			q_{Z}.Fg.Ff \arrow[rr, "q_{gf}"']
			&{}
			& Gg.Gf.q_{X}
		\end{tikzcd}\begin{tikzcd}[column sep = 15, font=\fontsize{9}{6}]
			&Gg.p_{Y}.Ff\arrow[rd, "1.p_{f}"]
			\arrow[ddd, shorten = 10, "1_{Gg}.\sigma_{Y}.1" description]
			\\
			p_{Z}.Fg.Ff
			\arrow[dd, Rightarrow, shorten = 10, "\sigma_{g}.1", shift left = 10]
			\arrow[ddd, "\sigma_{Z}.1.1"']
			\arrow[ru, "p_{g}.1"]
			&& Gg.Gf.p_{X}
			\arrow[dd, Rightarrow, shorten = 10, "1.\sigma_{f}"', shift right = 10]
			\arrow[ddd, "1.1.\sigma_{X}"]
			\\
			\\
			{}&Gg.q_{Y}.Ff
			\arrow[rd, "1.q_{f}"]
			\arrow[d, Rightarrow, shorten = 5, "q_{g{,}f}"]
			&{}
			\\
			q_{Z}.Fg.Ff
			\arrow[ru, "q_{g}.1"]
			\arrow[rr, "q_{gf}"']
			&{}
			& Gg.Gf.q_{X}
		\end{tikzcd}
	\end{itemize}
	
	\noindent Call a trimodification
	
	\begin{itemize}
		\item \emph{strict} if its $3$-cell components are identities.
		\item \emph{costrict} if its $2$-cell components are identities.
	\end{itemize}
\end{definition}

\begin{remark}
	The $3$-cell component $\sigma_{f}$ of a trimodification $\sigma$ has a uniquely determined mate given by a $3$-cell $\sigma_{f}^{*}$ which can be produced by pasting $\sigma_{f}$ along $p_{f}^{\varepsilon}$ and $q_{f}^{\varepsilon}$. By mateship, $\sigma_{f}^{*}$ also satisfies conditions which are dual to those specified for $\sigma_{f}$.
	
	$$\begin{tikzcd}[font=\fontsize{9}{6}]
		p_{Y}.Ff
		\arrow[rr, leftarrow, "p_{f}^{*}"]
		\arrow[dd, "\sigma_{Y}.1_{Ff}"']
		& {}\arrow[dd, Rightarrow, shorten = 10, "\sigma_{f}^{*}"]
		&
		Gf.p_{X}
		\arrow[dd, "1_{Gf}.\sigma_{X}"]
		\\
		\\
		q_{Y}.Ff \arrow[rr, leftarrow, "q_{f}^{*}"']
		&{}&
		Gf.q_{X} 
	\end{tikzcd} $$
	
\end{remark}
\begin{definition}
	Let $\sigma: p \rightarrow q$ be a trimodification as in Definition \ref{trimodification definition} and let $\tau: p \rightarrow q$ be another trimodification. A \emph{perturbation} $\sigma \rightarrow \tau$ consists of the assignation to every object $X \in \mathfrak{A}$, a $3$-cell $\Omega: \sigma_{X} \Rrightarrow \tau_{X}$ in $\mathfrak{B}$ such that the following equation holds for every $f: X \rightarrow Y$ in $\mathfrak{A}$.

	$$\begin{tikzcd}[column sep = 18, row sep = 20, font=\fontsize{9}{6}]
		p_Y.Ff
		\arrow[rr,bend left = 45, "\sigma_{Y}.1"]
		\arrow[dd,"p_f"']
		&
		{}
		\arrow[d,Rightarrow, shorten = 5,"\sigma_{f}"]
		&
		q_Y.Ff
		\arrow[dd,"q_{f}"]
		\\
		&{}&&=
		\\
		Gf.p_X
		\arrow[rr,bend left = 45,"1.\sigma_{X}" {name = C}]
		\arrow[rr, bend right = 45, "1.\tau_{X}"' {name = D}]
		&
		{}
		&
		Gf.q_X
		\arrow[from =C, to =D, Rightarrow, shorten = 10, shift right = 5, "1.\Omega_{X}"]
	\end{tikzcd}	\begin{tikzcd}[column sep = 18, row sep = 20]
		p_Y.Ff
		\arrow[rr,bend left = 45, "\sigma_{Y}.1" {name = A}]
		\arrow[rr, bend right = 45, "\tau_{Y}.1"' {name = B}]
		\arrow[dd,"p_f"']
		&
		&
		q_Y.Ff
		\arrow[dd,"q_{f}"]
		\\
		&{}
		\arrow[d,Rightarrow, shorten = 5,"\tau_f"]
		\\
		Gf.p_{X}\
		\arrow[rr, bend right = 45, "1.\tau_{X}"']
		&
		{}
		&
		Gf.q_{X} &{}
		\arrow[from =A, to =B, Rightarrow, shorten = 10, shift right = 5, "\Omega_{Y}.1"]
	\end{tikzcd} $$
	
\end{definition}

\subsection{$\mathbf{Gray}$-categories of $\left(3, k\right)$-transfors}\label{Subsection Gray categories of 3 k transfors}

\begin{remark}\label{Gray Functor Gray category with weak higher cells}
	It is known (Theorem 7.2.2 of \cite{Gurski PhD}) that when $\mathfrak{B}$ is a $\mathbf{Gray}$-category then for any tricategory $\mathfrak{A}$, there is a $\mathbf{Gray}$-category $\mathbf{Tricat}\left(\mathfrak{A}, \mathfrak{B}\right)$ whose objects are trihomomorphisms, morphisms are trinatural transformations, $2$-cells are trimodifications and $3$-cells are perturbations. We review the structure of the full-sub-$\mathbf{Gray}$-category on $\mathbf{Gray}$-functors in the special setting where $\mathfrak{A}$ is also a $\mathbf{Gray}$-category. We refer the reader to chapter 7 of \cite{Gurski PhD} for proofs.
	
	\begin{enumerate}
		\item Given $\mathbf{Gray}$-functors $F, G: \mathfrak{A} \rightarrow \mathfrak{B}$ and trinatural transformations $p, q: F \rightarrow G$, the hom-category $\mathbf{Tricat}_{s}\left(\mathfrak{A}, \mathfrak{B}\right)\left(F, G\right)\left(p, q\right)$ has objects given by trimodifications $\sigma: p \Rightarrow q$ and morphisms given by perturbations $\Omega: \sigma \Rrightarrow \tau$. Composition in this category is inherited component-wise from the hom-categories of $\mathfrak{B}$, and the identity perturbation on $\sigma$ has $3$-cell component on $X \in \mathfrak{A}$ given by the identity $1_{\sigma_{X}}$.
		\item Given $\mathbf{Gray}$-functors $F, G: \mathfrak{A} \rightarrow \mathfrak{B}$, we describe the hom-$2$-category $\mathbf{Tricat}_{s}\left(\mathfrak{A}, \mathfrak{B}\right)\left(F, G\right)$. \begin{itemize}
			\item Its objects are trinatural transformations $p: F \rightarrow G$, with the identity trimodification on $p$ having identity $2$-cell and $3$-cell components. 
			\item Its hom-categories between trinatural transformations $p$ and $q$ are given as described in Part (1).
			\item Given trimodifications \begin{tikzcd}
				p \arrow[r, "\sigma"] &q \arrow[r, "\tau"] & r
			\end{tikzcd}, their composite is the trimodification whose $2$-cell component on $X$ is given by \begin{tikzcd}
				{p}_{X} \arrow[r, "{\sigma}_{X}"] &{q}_{X} \arrow[r, "{\tau}_{X}"] & {r}_{X}
			\end{tikzcd} 
			and whose $3$-cell component on $f: X \rightarrow Y$ is given by the following pasting in the hom-$2$-category $\mathfrak{B}\left(FX, GY\right)$.
			
			$$\begin{tikzcd}[font=\fontsize{9}{6}]
				p_{Y}.Ff
				\arrow[dd, "p_{f}"']
				\arrow[rr, "\sigma_{Y}.1"]
				& {}\arrow[dd, Leftarrow, shorten = 10, "\sigma_{f}"]
				&
				q_{Y}.Ff \arrow[dd, "q_{f}"] \arrow[rr, "\tau_{Y}.{1}"]
				&{}\arrow[dd, shorten = 10, "\tau_{f}", Leftarrow]
				&r_{Y}.Ff\arrow[dd, "r_{f}"]
				\\
				\\
				Gf.p_{X} 
				\arrow[rr, "1.\sigma_{X}"']
				&{}&
				Gf.q_{X}\arrow[rr, "1.\tau_{X}"']
				&{}&Gf.r_{X}
			\end{tikzcd}$$ 
			\item Given horizontally composable perturbations 
			\begin{tikzcd}
				p
				\arrow[rr, bend left, "\sigma" name = A]
				\arrow[rr, bend right, "\sigma '"' name = B]
				&& q
				\arrow[rr,bend left, "\tau" name = C]
				\arrow[rr, bend right, "\tau '"' name = D]
				&& r
				\arrow[from = A, to = B, Rightarrow, shorten = 5, "\Sigma"]
				\arrow[from = C, to = D, Rightarrow, shorten =5, "\Omega"]
			\end{tikzcd} their horizontal composite $\Omega * \Sigma: \tau. \sigma \Rrightarrow \tau'. \sigma'$ has component on $X$ given by the horizontal composite of their components $\Omega_{X}*\Sigma_{X}$ in the hom-$2$-category $\mathfrak{B}\left(FX, GX\right)$.
		\end{itemize} 
		\item The identity trinatural transformation on a $\mathbf{Gray}$-functor $F$ is the strict trinatural transformation whose $1$-cell components on objects are also identities.
		\item Given $\mathbf{Gray}$-functors $F, G, H: \mathfrak{A} \rightarrow \mathfrak{B}$, the composition $2$-functor
		
		$$\circ: \mathbf{Tricat}_{s}\left(\mathfrak{A}, \mathfrak{B}\right)\left(G, H\right) \otimes \mathbf{Tricat}_{s}\left(\mathfrak{A}, \mathfrak{B}\right)\left(F, G\right) \rightarrow \mathbf{Tricat}_{s}\left(\mathfrak{A}, \mathfrak{B}\right)\left(F, H\right)$$
		
		\noindent is defined in the following way. \begin{enumerate}
			\item Given a pair of composable trinatural transformations \begin{tikzcd}
				F \arrow[r, "p"] & G \arrow[r, "q"] & H
			\end{tikzcd}, their composite will have \begin{itemize}
				\item $1$-cell component on $X$ given by ${\left(q \circ p\right)}_{X}:=$ \begin{tikzcd}
					FX \arrow[r, "p_{X}"] & GX \arrow[r, "q_{X}"] & HX
				\end{tikzcd}
				\item $2$-cell component on $f: X \rightarrow Y$ given by the $2$-cell depicted below, with the rest of the adjoint equivalence constructed similarly.
				$${\left(q \circ p\right)}_{f}:= \begin{tikzcd}
					q_{Y}.p_{Y}.Ff \arrow[r, "1.p_{f}"] & q_{Y}.Gf.p_{X} \arrow[r, "q_{f}.1"] & Hf.q_{X}.p_{X}
				\end{tikzcd}$$
				\item $3$-cell component ${\left(q \circ p\right)}_{\phi}$ on $\phi: f \Rightarrow g: X \rightarrow Y$ given by the following pasting in the hom-$2$-category $\mathfrak{B}\left(FX, HY\right)$, with the $3$-cell component ${\left(q\circ p\right)}_{\phi}^{*}$ defined similarly.
				
				$$\begin{tikzcd}[column sep = 12, font=\fontsize{9}{6}]
					q_{Y}.p_{Y}.Ff
					\arrow[dd,"1.1.F\phi"']
					\arrow[rr,"1.p_f"]
					&
					{}
					\arrow[dd,Leftarrow,shorten=15,"1.p_\phi"]
					&
					q_{Y}.Gf.p_{X} \arrow[rr, "q_{f}.1"]
					\arrow[dd,"1.G\phi .1"]
					&{}\arrow[dd, Leftarrow, shorten = 15, "q_{\phi}.1"]
					& Hf.q_{X}.p_{X}\arrow[dd, "H\phi .1.1"]
					\\
					\\
					q_{Y}p_{Y}.Fg
					\arrow[rr,"1.p_{g}"']
					&
					{}
					&
					q_{Y}.Gg.p_X \arrow[rr, "q_{g}.1"']
					&{}
					& Hg.q_{X}.p_{X}
				\end{tikzcd}$$
				\item Unitor at $X$ given by the following horizontal composite in the hom-$2$-category $\mathfrak{B}\left(FX, HX\right)$. 
				
				$$\begin{tikzcd}[font=\fontsize{9}{6}]
					q_{X}p_{X}
					\arrow[rr, bend left, "1.1_{p_{X}}" name = A]
					\arrow[rr, bend right, "1.p_{1_{X}}"' name = B]
					&& q_{X}.p_{X}
					\arrow[rr,bend left, "1_{q_{X}}.1" name = C]
					\arrow[rr, bend right, "q_{1_{X}}.1"' name = D]
					&& q_{X}.p_{X}
					\arrow[from = A, to = B, Rightarrow, shorten = 5, "1.p^{X}"]
					\arrow[from = C, to = D, Rightarrow, shorten =5, "q^{X}.1"]
				\end{tikzcd}$$
				\item Compositor at \begin{tikzcd}
					X \arrow[r, "f"] & Y \arrow[r, "g"] & Z
				\end{tikzcd} given by the following pasting in the hom-$2$-category $\mathfrak{B}\left(FX, HZ\right)$.

				$$\begin{tikzcd}[font=\fontsize{9}{6}, column sep = 30]
					&q_{Z}.p_{Z}.Fg.Ff
					\arrow[d, "1.p_{g}.1"']
					\arrow[ddr, "1.p_{gf}", bend left, shift left = 3]
					\\
					&q_{Z}.Gg.p_{Y}.Ff
					\arrow[r, shorten = 15, Rightarrow, "1.p_{g{,}f}"]
					\arrow[ld, "q_{g}.1.1"']
					\arrow[rd, "1.1.p_{f}"'] &{}
					\\
					Hg.q_{Y}.p_{Y}.Ff
					\arrow[rr, Rightarrow, shorten = 40, "{\left(q_{g}\right)}_{\left(p_{f}\right)}"]
					\arrow[rd, "1.1.p_{f}"'] &{}& q_{Z}.Gg.Gf.p_{X}
					\arrow[ldd, "q_{gf}.1", bend left, shift left = 3]
					\arrow[ld, "q_{g}.1.1"']
					\\
					&Hg.q_{Y}.Gf.p_{X}
					\arrow[r, Rightarrow, shorten = 15, "q_{g{,}f}.1"]
					\arrow[d, "1.q_{f}.1"']
					&{}
					\\
					& Hg.Hf.q_{X}.p_{X}
				\end{tikzcd}$$
			\end{itemize}
			\item Consider \begin{tikzcd}
				F\arrow[rr, "p"]
				&& G
				\arrow[rr,bend left, "q" name = C]
				\arrow[rr, bend right, "q '"' name = D]
				&& H
				\arrow[from = C, to = D, Rightarrow, shorten =5, "\tau"]
			\end{tikzcd}, where $p, q$ and $q'$ are trinatural transformations and $\tau$ is a trimodification. Then the whiskering is given by the trimodification whose $2$-cell component on $X$ is given by the whiskering of $\tau_{X}$ with $p_{X}$, and whose $3$-cell component on $f: X \rightarrow Y$ is given by the pasting in the hom-$2$-category $\mathfrak{B}\left(FX, HX\right)$ depicted below left.
			\item Consider \begin{tikzcd}
				F
				\arrow[rr, bend left, "p" name = A]
				\arrow[rr, bend right, "p '"' name = B]
				&& G\arrow[rr, "q"]
				&& H
				\arrow[from = A, to = B, Rightarrow, shorten = 5, "\sigma"]
			\end{tikzcd}, where $p, p'$ and $q$ are trinatural transformations and $\sigma$ is a trimodification. Then the whiskering is given by the trimodification whose $2$-cell component on $X$ is given by the whiskering of $\sigma_{X}$ with $q_{X}$, and whose $3$-cell component on $f: X \rightarrow Y$ is given by the pasting in the hom-$2$-category $\mathfrak{B}\left(FX, HX\right)$ depicted below right.
			\item The whiskering of perturbations $\Omega$ with trinatural transformations $p$ on either side is given by the component-wise whiskering of the $3$-cell $\Omega_{X}$ in $\mathfrak{B}$ with the $1$-cell $p_{X}$ in $\mathfrak{B}$.
			\item Let \begin{tikzcd}
				F
				\arrow[rr, bend left, "p" name = A]
				\arrow[rr, bend right, "p '"' name = B]
				&& G
				\arrow[rr,bend left, "q" name = C]
				\arrow[rr, bend right, "q '"' name = D]
				&& H
				\arrow[from = A, to = B, Rightarrow, shorten = 5, "\sigma"]
				\arrow[from = C, to = D, Rightarrow, shorten =5, "\tau"]
			\end{tikzcd} be a pair of interchangeable trimodifications. Then their interchanger $\tau_{\sigma}$ is the perturbation whose component at an object $X$ is given by the interchanger in $\mathfrak{B}$ of $\tau_{X}$ and $\sigma_{X}$.
			
		\end{enumerate}
	\end{enumerate}
	
	$$\begin{tikzcd}[column sep = 8, font=\fontsize{9}{6}]
		q_{Y}.p_{Y}.Ff
		\arrow[dd,"1.\sigma_{Y}.1"description]
		\arrow[rr,"1.p_f"]
		&
		{}
		\arrow[dd,Rightarrow,shorten=15,"1.\sigma_f", shift right = 5]
		&
		q_{Y}.Gf.p_{X} \arrow[rr, "q_{f}.1"]
		\arrow[dd,"1.1.\sigma_{X}"description]
		&{}\arrow[dd, Rightarrow, shorten = 15, "{\left(q_{f}\right)}_{\left(\sigma_{X}\right)}", shift right = 5]
		& Hf.q_{X}.p_{X}\arrow[dd, "1.1.\sigma_{X}"description]
		\\
		\\
		q_{Y}{p'}_{Y}.Ff
		\arrow[rr,"1.{p'}_{f}"']
		&
		{}
		&
		q_{Y}.Gf.{p'}_{X} \arrow[rr, "q_{f}.1"']
		&{}
		& Hf.q_{X}.{p'}_{X}&{}
	\end{tikzcd}\begin{tikzcd}[column sep = 8, font=\fontsize{9}{6}]
		q_{Y}.p_{Y}.Ff
		\arrow[dd,"\tau_{Y}.1.1"description]
		\arrow[rr,"1.p_f"]
		&
		{}
		\arrow[dd,Rightarrow,shorten=15,shift right = 4, "{\left(\tau_{Y}\right)}_{\left(p_{f}\right)}"]
		&
		q_{Y}.Gf.p_{X} \arrow[rr, "q_{f}.1"]
		\arrow[dd,"\tau_{Y}.1.1"description]
		&{}\arrow[dd, Rightarrow, shorten = 15, "\tau_{f}.1"]
		& Hf.q_{X}.p_{X}\arrow[dd, "1.\tau_{X}.1"description]
		\\
		\\
		{q'}_{Y}{p}_{Y}.Ff
		\arrow[rr,"1.{p}_{f}"']
		&
		{}
		&
		{q'}_{Y}.Gf.{p}_{X} \arrow[rr, "{q'}_{f}.1"']
		&{}
		& Hf.{q'}_{X}.{p}_{X}
	\end{tikzcd}$$
	
\end{remark}

\begin{remark}\label{terminology for higher cells between tricats}
	What we call an `invertible ico-icon' is precisely an invertible version of an ico-icon as in \cite{Low dimensional structures formed by tricategories}, and what we call a `pseudo-icon equivalence' is exactly the data of a pseudo-icon $p$ as defined in \cite{Low dimensional structures formed by tricategories} and a specified adjoint equivalence in $\mathbf{Tricat}_{3}\left(\mathfrak{A}, \mathfrak{B}\right)$ of which $p$ is the left adjoint. Strict trinatural transformations are precisely $\mathbf{Gray}$-natural transformations. The adjective `semi-strict' will be justified by Proposition \ref{Trinatural transformations out of cofibrant Gray categories} and Theorem \ref{Gr on weak higher cells between tricategories} part (1). They have been considered in Definition 15.1 of \cite{Gurski Coherence in Three Dimensional Category Theory} where they were called `$1$-strict', and a lax version was considered in \cite{crans tensor of gray categories} where they were called `lax-$1$-transfors'. We also note that in the semi-strict case the left and right whiskering laws simplify somewhat and the associativity and left and right unit laws hold automatically. Strict trimodifications between $\mathbf{Gray}$-natural transformations are precisely the `$\mathbf{Gray}$-modifications' that appear as $2$-cells in the internal homs on $\mathcal{V}$-$\mathbf{Cat}$ which are inherited from $\mathcal{V} = \mathbf{Gray}$. Costrict trimodifications are precisely the costrict $2$-cells, in the sense described in section 5 of \cite{Icons}, in the underlying sesquicategory of $\mathbf{Tricat}_{s}\left(\mathfrak{A}, \mathfrak{B}\right)$.
	\\
	\\
	\noindent With some care, $3$-pseudofunctors and $3$-pseudonatural transformations can be seen as certain $\mathcal{V}$-enriched pseudofunctors and pseudonatural transformations in the sense of \cite{Garner Shulman Enriched Categories as a Free Cocompletion}, for an appropriate monoidal $2$-category $\mathcal{V}$. The underlying monoidal category of $\mathcal{V}$ is given by $2$-$\mathbf{Cat}$ with the $\mathbf{Gray}$-tensor product. Now, $2$-natural transformations are not closed under the $\mathbf{Gray}$-tensor product, but $2$-cells in $\mathbf{V}$ are given by pseudonatural transformations which can be built from $2$-natural transformations via composition and the $\mathbf{Gray}$-tensor product. Then $3$-pseudofunctors and $3$-pseudonatural transformations are $\mathcal{V}$-enriched pseudofunctors and pseudonatural transformations, in which the $2$-cell components are actually $2$-natural rather than built out of $2$-natural transformations in this way. This perspective will not be needed for the main results of this paper.
\end{remark}

\begin{remark}\label{Remark free vs operational three dimensions}
	We distinguish between free and operational coherence data in $\left(3, k\right)$-transfors out of a $\mathbf{Gray}$-category $\mathfrak{A}$ into some other $\mathbf{Gray}$-category $\mathfrak{B}$. This extends the notions introduced in subsection \ref{subsection free versus derivable coherences} to the three dimensional setting.
	\begin{itemize}
		\item ($k = 0$) All coherence data of trihomomorphisms are operational, mediating preservation of the sesquicategorical operations in the underlying two-dimensional structure of the domain. In particular,
		\begin{itemize}
			\item There are pseudofunctors between hom-$2$-categories, which as we have already discussed in Section \ref{subsection free versus derivable coherences} have only operational data, 
			\item The $1$-cell components of $\chi$ and $\iota$ in the hom-$2$-categories of $\mathfrak{B}$ mediate preservation of composition and identities of morphisms in $\mathfrak{A}$,
			\item The $2$-cell components of $\chi$ in the hom-$2$-categories of $\mathfrak{B}$ mediate whiskering of $1$-cells by $2$-cells in $\mathfrak{A}$,
			\item The $2$-cell components of the invertible modifications $\gamma$, $\delta$ and $\omega$ mediate respect for left and right unit laws, and the associativity law.
		\end{itemize}
		\item ($k = 1$) The compositors and unitors are operational, mediating the axioms of a pseudonatural transformation pertaining to respect for composition and identities. On the other hand, the following coherence data are free.
		\begin{itemize}
			\item For $f: X \rightarrow Y$, the adjoint equivalences $p_{f} \dashv p_{f}^{*}$ mediates varying an object along a morphism.
			\item For $\phi: f \Rightarrow g$, the invertible $3$-cell $p_{\phi}$ mediates varying a morphism along a $2$-cell.
		\end{itemize}
		\item ($k = 2$) The only coherence datum is the $3$-cell $\sigma_{f}$, which mediates varying an object along a morphism. It is a free coherence.
		\item ($k = 3$) Perturbations do not have any coherence data.
	\end{itemize}
	\noindent In the $n = 2$ dimensional setting we observed that when the underlying $\left(n-1\right)$-dimensional categorical structure (i.e. underlying category) of the domain is free on a graph, then it is possible to replace a $\left(2, k\right)$-transfor with a suitably equivalent one in which operational coherences are identities, but it is not in general possible to also make free coherences identities. This observation will be extended to the three dimensional setting in Subsection \ref{Subsection Extending strictification to 3 k transfors}.
\end{remark}

\noindent By the discussion in Remark \ref{Remark free vs operational three dimensions}, a trinatural transformation is semi-strict if and only if its operational coherences are identities. The strictification of a trinatural transformation will have this property. Semi-strict trinatural transformations have the undesirable property of not being closed under composition, an observation that we record in Lemma \ref{tautological condition for failure of semi-stricts to be closed under composition}, to follow. Section \ref{Section Solutios to the Failure of Semi-strict trinatural transformations to be closed under composition} is devoted to solutions to this failure of closure under composition.

\begin{lemma}\label{tautological condition for failure of semi-stricts to be closed under composition}
	Suppose \begin{tikzcd}
		F \arrow[r, "p"] & G \arrow[r, "q"] & H
	\end{tikzcd} is a composable pair of semi-strict trinatural transformations between $\mathbf{Gray}$-functors from $\mathfrak{A}$ to $\mathfrak{B}$.
	
	\begin{enumerate}
		\item The composite $q \circ p$ is unital.
		\item The composite $q \circ p$ is compositional if and only if the following whiskered interchanger is the identity. Note that the unlabelled regions strictly commute by compositionality of $p$ and $q$.

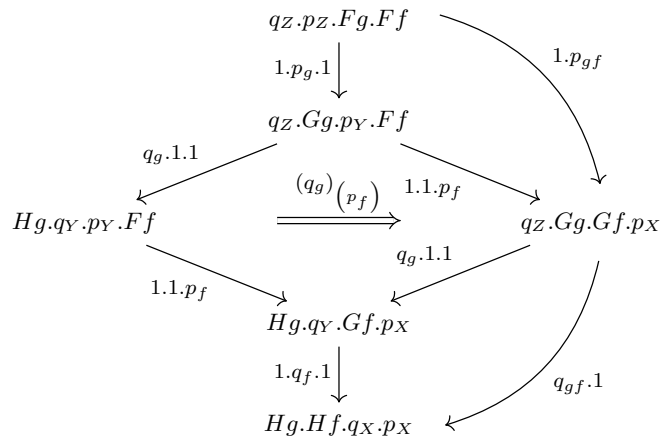
		$$\begin{tikzcd}[font=\fontsize{9}{6}, column sep = 30]
			&q_{Z}.p_{Z}.Fg.Ff
			\arrow[d, "1.p_{g}.1"']
			\arrow[ddr, "1.p_{gf}", bend left, shift left = 3]
			\\
			&q_{Z}.Gg.p_{Y}.Ff
			\arrow[ld, "q_{g}.1.1"']
			\arrow[rd, "1.1.p_{f}"'] &{}
			\\
			Hg.q_{Y}.p_{Y}.Ff
			\arrow[rr, Rightarrow, shorten = 40, "{\left(q_{g}\right)}_{\left(p_{f}\right)}"]
			\arrow[rd, "1.1.p_{f}"'] &{}& q_{Z}.Gg.Gf.p_{X}
			\arrow[ldd, "q_{gf}.1", bend left, shift left = 3]
			\arrow[ld, "q_{g}.1.1"']
			\\
			&Hg.q_{Y}.Gf.p_{X}
			\arrow[d, "1.q_{f}.1"']
			&{}
			\\
			& Hg.Hf.q_{X}.p_{X}
		\end{tikzcd}$$
		\item $q\circ p$ is not semi-strict in general.
	\end{enumerate}
\end{lemma}

\section{Three-dimensional strictification}\label{Section Extending strictification to 3, k transfors}

\subsection{Gurski's $\mathbf{Gr}$}\label{Subsection Gurski's Gr}

\noindent We describe the underlying $2$-computad of a tricategory $\mathfrak{A}$. The free sesquicategory on this $2$-computad underlies the $\mathbf{Gray}$-category $\mathbf{Gr}\left(\mathfrak{A}\right)$. As such, strictification of tricategories can be chosen so that the resulting structure is free in codimension one.
\\
\\
\noindent Given a path of morphisms $f$ in $\mathfrak{A}$ and an association $\alpha$ of $f$, let ${f}_{\alpha}$ denote the morphism in $\mathfrak{A}$ produced by evaluating the path $f$ according to the association $\alpha$. There is a $2$-computad $C_{2}\left(\mathfrak{A}\right)$ given by the following data.
	
	\begin{itemize}
		\item Objects are given by those of $\mathfrak{A}$.
		\item Generating edges $f: X \rightarrow Y$ are given by morphisms $f: X \rightarrow Y$ in $\mathfrak{A}$.
		\item $2$-cells $\left(\alpha_{0}, \alpha_{1}, \phi\right): f \Rightarrow g: X \rightarrow Y$ are given by the data of
		
		\begin{itemize}
			\item An association ${\alpha_{0}}$ of $f$.
			\item An association ${\alpha_{1}}$ of $g$.
			\item A $2$-cell $\phi: f_{\alpha_{0}} \Rightarrow g_{\alpha_{1}}$ in $\mathfrak{A}$.
		\end{itemize}
	\end{itemize}

\noindent The underlying sesquicategory of the $\mathbf{Gray}$-category $\mathbf{Gr}\left(\mathfrak{A}\right)$ described in Section 10.4 of \cite{Gurski Coherence in Three Dimensional Category Theory} is free on the $2$-computad that we have just described. The $3$-cells in $\mathbf{Gr}\left(\mathfrak{A}\right)$ are defined with respect to some arbitrary choice of evaluations of their source and target $2$-cells, into actual $2$-cells in the tricategory $\mathfrak{A}$. This $\mathbf{Gray}$-category is indeed well-defined (Theorem 10.8, \cite{Gurski Coherence in Three Dimensional Category Theory}), and is triequivalent to the original tricategory $\mathfrak{A}$ (Theorem 10.9, \cite{Gurski Coherence in Three Dimensional Category Theory}).
\\
\\
\noindent We refer to $\mathbf{Gray}$-categories whose underlying sesquicategories are free on a $2$-computad as \emph{cofibrant}, since this condition characterises cofibrancy in the Lack model structure on the category $\mathbf{Gray}\text{-}\mathbf{Cat}$ \cite{Quillen Gray Cat}. The significance of cofibrancy of $\mathbf{Gr}\left(\mathfrak{A}\right)$ for strictification of $\left(3, k\right)$-transfors in analogous to what we discussed in Remark \ref{pseudofunctors out of cofibrant 2-categories}. Proposition \ref{trihomomorphism out of cofibrant Gray category}, to follow, records this precisely. The existence of a trinaturally biequivalent $\mathbf{Gray}$-functor $\overline{F}$ has been shown in Corollary 4.1 of \cite{Gray categories model algebraic tricategories}, however we give an explicit construction which demonstrates how freeness is used. 

\begin{proposition}\label{trihomomorphism out of cofibrant Gray category}
	Let $\mathfrak{A}$ and $\mathfrak{B}$ be $\mathbf{Gray}$-category and suppose the underlying sesquicategory of $\mathfrak{A}$ is free on a $2$-computad $\mathbb{A}$. Let $F: \mathfrak{A} \rightsquigarrow \mathfrak{B}$ be a trihomomorphism. Then there is \begin{enumerate}
		\item a $\mathbf{Gray}$-functor $\bar{F}: \mathfrak{A} \rightarrow \mathfrak{B}$ which agrees with $F$ on data in $\mathbb{A}$ and on $3$-cells.
		\item a unital pseudo-icon equivalence $e: \bar{F} \rightarrow F$ whose components on identities are given by the unitors of $F$, components on generating morphisms are given by identities, and components on a composite \begin{tikzcd}
			X \arrow[r, "g"] & Y \arrow[r, "h"] &Z
		\end{tikzcd} are given by \begin{itemize}
			\item The pasting depicted below if $g$ is generating and $h$ is not an identity. In this case the compositor $e_{h, g}$ is the identity.
		\end{itemize}
		$$\begin{tikzcd}[row sep = 20, column sep = 30, font=\fontsize{9}{6}]
			FX \arrow[r, "\bar{F}g"]\arrow[d, "1_{FX}"']
			&FY \arrow[r, "\bar{F}h"] \arrow[d, "1_{FY}"]
			&FZ\arrow[d, "1_{FZ}"]
			\arrow[ld, Rightarrow, shorten = 5, "e_{h}"]
			\\
			FX \arrow[rr, bend right = 60, "F\left(hg\right)"'] \arrow[r, "Fg"']
			\arrow[r, equal, shorten = 5, shift left = 8]
			& FY \arrow[r, "Fh"']\arrow[d, Rightarrow, "\chi_{h{,}g}", shorten = 2]
			& FZ
			\\
			&{}
		\end{tikzcd}$$
		\begin{itemize}
			\item The identity if either $h$ or $g$ is the identity and the other is generating, in which case the compositor $e_{h, g}$ is given by $\gamma_{g}$ or $\delta_{h}$ respectively.
			\item The pasting depicted below if $g$ is generating and $h = h_{2}.h_{1}$ with $h_{1}$ generating and $h_{2}$ non-identity. In this case the compositor $e_{h, g}$ is given by a whiskering of the associator $\omega_{h_{2}, h_{1}, g}$ of $F$ by the $2$-cell $e_{h_{2}}$.
			
			$$\begin{tikzcd}[row sep = 20, column sep = 35, font=\fontsize{9}{6}]
				FX \arrow[r, "\bar{F}g"]
				\arrow[d, "1_{FX}"']
				&FY \arrow[r, "\bar{F}h_{1}"]
				\arrow[d, "1_{FY}"]
				&FW
				\arrow[d, "1_{FW}"]
				\arrow[r, "\bar{F}h_{2}"]
				&FZ
				\arrow[ld, Rightarrow, shorten=2,"e_{h_{2}}"]
				\arrow[d, "1_{FZ}"]
				\\
				FX 
				\arrow[rrr, bend right = 80,near end, "F\left(h_{2}h_{1}g\right)"']
				\arrow[rr, bend right = 60, "F\left(h_{1}g\right)"']
				\arrow[r, "Fg"']
				\arrow[r, equal, shorten = 5, shift left = 8]
				& FY 
				\arrow[r, equal, shorten = 5, shift left = 8]
				\arrow[r, "Fh_{1}"']\arrow[d, Rightarrow, "\chi_{h{,}g}", shorten = 2]
				& FW
				\arrow[dd, Rightarrow, shorten = 10, "\chi_{h_{2}{,}h_{1}g}"]
				\arrow[r, "Fh_{2}"']
				&FZ
				\\
				&{}
				\\
				&&{}
			\end{tikzcd}$$
		\end{itemize}
	\end{enumerate}		
	
\end{proposition}

\begin{proof}
	For part (1), $\mathbf{Gray}$-functoriality requires checking preservation of relations on $3$-cells, such as their identities, composition and whiskering by $1$-cells and $2$-cells. This follows by coherence for the trihomomorphism $F$ (Theorem 10.13 of \cite{Gurski Coherence in Three Dimensional Category Theory}). For part (2), we describe the remaining data of the pseudo-icon equivalence $e: \bar{F} \rightarrow F$. Its $2$-cell components when $g$ is not generating are determined by the data described above given that $e$ needs to satisfy the associativity axiom, using the pentagon equation for $F: \mathfrak{A} \rightsquigarrow \mathfrak{B}$. The $3$-cell component at a generating $2$-cell of the form $\phi: f_{m}.f_{m-1}... f_{0} \Rightarrow g_{n}.g_{n-1}...g_{0}$ is given in terms of the pseudonaturality constraint of the compositor of $F$ at $\phi$. This is extended to general $2$-cells using the coherence data of the pseudofunctors between hom $2$-categories of $F: \mathfrak{A} \rightsquigarrow \mathfrak{B}$. The whiskering laws then follow from the modification axiom for $\omega$, while naturality in $3$-cells follows from the analogous property for $\chi$.
\end{proof}

\subsection{$3$-pseudofunctors and path objects in $\mathbf{Gray}\text{-}\mathbf{Cat}$}\label{Path Objects Section}

\noindent The aim of this Section is to describe generalised path objects for a $\mathbf{Gray}$-category $\mathfrak{B}$. These will be sub-$\mathbf{Gray}$-categories of $\mathbf{Tricat}\left(\mathfrak{P}_{m}, \mathfrak{B}\right)$, where $\mathfrak{P}_{m}$ denotes the free-living $k$-cell. In particular, the sub-$\mathbf{Gray}$-categories of $\mathbf{Tricat}\left(\mathfrak{P}_{m}, \mathfrak{B}\right)$ in which we will be interested has objects given by $\mathbf{Gray}$-functors $\mathfrak{P}_{m} \rightarrow \mathfrak{B}$, morphisms given by semi-strict trinatural transformations, and higher cells given by arbitrary $\left(3, k\right)$-transfors. By Lemma \ref{Gray Functor Gray category with weak higher cells}, semi-strict trinatural transformations out of $\mathfrak{P}_{m}$ \emph{are} closed under composition. Sufficient conditions for this closure, which apply to $\mathfrak{P}_{n}$, will be given in Proposition \ref{Gray categories with only trivial composites}.

\begin{example}\label{coherent hom from 2 to B}
	Let $\mathfrak{B}$ be a $\mathbf{Gray}$-category. We use Proposition \ref{Gray categories with only trivial composites} to give an explicit description of the $\mathbf{Gray}$-category $[\mathfrak{P}_{1}, \mathfrak{B}]$. Its objects are $\mathbf{Gray}$-functors $\mathfrak{P}_{1} \rightarrow \mathfrak{B}$. These consist of just the data of a morphism $x: X_{0} \rightarrow X_{1}$ in $\mathfrak{B}$, but we will write an object of $[\mathfrak{P}_{1}, \mathfrak{B}]$ as a triple $\left(x, X_{0}, X_{1}\right)$. A morphism $\left(f \dashv f^{*}, f_{0}, f_{1}\right): \left(x, X_{0}, X_{1}\right) \rightarrow \left(y, Y_{0}, Y_{1}\right)$ is a a semi-strict trinatural transformation, which consists of precisely the data of a pair of morphisms $f_{0}:X_{0} \rightarrow Y_{0}$ and $f_{1}: X_{1} \rightarrow Y_{1}$, and an adjoint equivalence in $\mathfrak{B}\left(X_{0}, Y_{1}\right)$ whose left adjoint is written as $f: y.f_{0} \rightarrow f_{1}.x$ and whose right adjoint is written as $f^{*}: f_{1}.x \rightarrow y.f_{0}$. A $2$-cell $\left(\phi, \phi_{0}, \phi_{1}\right): \left(f \dashv f^{*}, f_{0}, f_{1}\right) \Rightarrow \left(g, g^{*}, g_{0}, g_{1}\right)$ is a trimodification. This consists of precisely the data of a pair of $2$-cells $\phi_{0}: f_{0} \Rightarrow g_{0}$ and $\phi_{1}: f_{1} \Rightarrow g_{1}$, as well as an invertible $3$-cell as depicted below left. Finally, a $3$-cell $\left(\Gamma_{0}, \Gamma_{1}\right): \left(\phi, \phi_{0}, \phi_{1}\right) \Rrightarrow \left(\psi, \psi_{0}, \psi_{1}\right)$ is a perturbation, which consists of a pair of $3$-cells $\Gamma_{0}: \phi_{0} \Rrightarrow \psi_{0}$ and $\Gamma_{1}: \phi_{1} \Rrightarrow \psi_{1}$ satisfying the equation depicted below right. 
	\\
	\begin{tikzcd}[font=\fontsize{9}{6}]
		y.f_{0} \arrow[rr, "1.\phi_{0}"]\arrow[dd, "f"'] &{}\arrow[dd, Rightarrow, shorten = 10, "\phi"]& y.g_{0}\arrow[dd, "g"]
		\\
		\\
		f_{1}.x \arrow[rr, "\phi_{1}.1"'] &{}& g_{1}.x &{}
	\end{tikzcd}\begin{tikzcd}[font=\fontsize{9}{6}]
		y.f_{0} \arrow[rr, bend right, "1.\psi_{0}"' name = B]\arrow[rr, bend left, "\phi_{0}" name = A]
		\arrow[dd, "f"'] &{}& y.g_{0}\arrow[dd, "g"]
		\\
		&&&=
		\\
		f_{1}.x \arrow[rr, bend right, "\psi_{1}.1"' name = C] &{}& g_{1}.x 
		\arrow[from= A, to =B, Rightarrow, "\Gamma_{0}", shorten = 10]\arrow[from =B, to =C, Rightarrow, shorten = 15, "\psi"]
	\end{tikzcd}\begin{tikzcd}[font=\fontsize{9}{6}]
		y.f_{0} \arrow[rr, bend left, "\phi_{0}" name = A]
		\arrow[dd, "f"'] &{}& y.g_{0}\arrow[dd, "g"]
		\\
		\\
		f_{1}.x \arrow[rr, bend right, "\psi_{1}.1"' name = C] \arrow[rr, bend left, "\phi_{1}.1" name = D]&{}& g_{1}.x &{}
		\arrow[from= A, to =D, Rightarrow, shorten = 12, "\phi"]\arrow[from =D, to =C, Rightarrow, shorten = 10, "\Gamma_{1}"]
	\end{tikzcd}
	\\
	\\
	\noindent The $\mathbf{Gray}$-category structure on $[\mathfrak{P}_{1}, \mathfrak{B}]$ can also be described directly in terms of pastings of cells in $\mathfrak{B}$. Note that the full-sub-$\mathbf{Gray}$-category of $[\mathfrak{P}_{1}, \mathfrak{B}]$ on those objects $x: X_{0} \rightarrow X_{1}$ which are part of a biequivalence internal to $\mathfrak{B}$ is precisely the path object constructed in the proof of Proposition 4.1 in \cite{Quillen Gray Cat}.
\end{example}

\begin{remark}\label{general trinatural transformations via 2-cat pseudofunctors}
	Consider the $\mathbf{Gray}$-category $[\mathfrak{P}_{1}, \mathfrak{B}]$ just defined in Example \ref{coherent hom from 2 to B}. There are evident $\mathbf{Gray}$-functors $s_{0}, t_{0}: [\mathfrak{P}_{1}, \mathfrak{B}] \rightarrow \mathfrak{B}$ which take sources and targets of all data in the following way.
	\begin{align*}
		s_{0}\left(x, X_{0}, X_{1}\right) = X_{0} && t_{0}\left(x, X_{0}, X_{1}\right) = X_{1}
		\\
		s_{0}\left(f \dashv f^{*}, f_{0}, f_{1}\right) = f_{0} && t_{0}\left(f \dashv f^{*}, f_{0}, f_{1}\right) = f_{1}
		\\
		s_{0}\left(\phi, \phi_{0}, \phi_{1}\right) = \phi_{0} && t_{0}\left(\phi, \phi_{0}, \phi_{1}\right) = \phi_{1}
		\\
		s_{0}\left(\Gamma_{0}, \Gamma_{1}\right) = \Gamma_{0} && t_{0}\left(\Gamma_{0}, \Gamma_{1}\right) = \Gamma_{1}
	\end{align*}
	\noindent Given a pair of $\mathbf{Gray}$-functors $\left(F, G\right): \mathfrak{A} \rightarrow \mathfrak{B} \times \mathfrak{B}$, a semi-strict trinatural transformation $p: F \rightarrow G$ is precisely the data of a $\mathbf{Gray}$-functor $p': \mathfrak{A} \rightarrow [\mathfrak{P}_{1}, \mathfrak{B}]$ giving a lifting of $\left(F, G\right)$ along $\left(s_{0}, t_{0}\right): [\mathfrak{P}_{1}, \mathfrak{B}] \rightarrow \mathfrak{B} \times \mathfrak{B}$. On the other hand, a general trinatural transformation $p: F \rightarrow G$ can also be described as a lifting of $\left(F, G\right): \mathfrak{A} \rightarrow \mathfrak{B} \times \mathfrak{B}$ along $\left(s_{0}, t_{0}\right): [\mathfrak{P}_{1}, \mathfrak{B}]$, this time not by a $\mathbf{Gray}$-functor but by somewhat weaker map $p': \mathfrak{A} \rightsquigarrow [\mathfrak{P}_{1}, \mathfrak{B}]$ as we will now make precise.
\end{remark}

\begin{proposition}\label{trinatural transformations as 2-Cat pseudofunctors}
	Let $F, G: \mathfrak{A} \rightarrow \mathfrak{B}$ be $\mathbf{Gray}$-functors.
	\begin{enumerate}
		\item There is a bijection between \begin{itemize}
			\item the set of trinatural transformations $p: F \rightarrow G$,
			\item the set of $3$-pseudofunctors $p': \mathfrak{A} \rightarrow [\mathbf{2}, \mathfrak{B}]$ satisfying $s_{0}.p' = F$ and $t_{0}.p' = G$.
		\end{itemize} 
		\noindent Moreover, this restricts to a bijection between semi-strict trinatural transformations and $\mathbf{Gray}$-functors.
		\item Let $p, q: F \rightsquigarrow G$ be trinatural transformations and let $p', q': \mathfrak{A} \rightarrow [\mathbf{2}, \mathfrak{B}]$ be the corresponding $3$-pseudofunctors from part (1). There is a bijection between \begin{itemize}
			\item the set of trimodificational adjoint equivalences $\sigma: p \Rightarrow q$
			\item the set of $3$-pseudonatural transformations $\sigma': p' \Rightarrow q'$ satisfying $s_{0}.\sigma' = 1_{F}$ and $t_{0}.\sigma' = 1_{G}$
		\end{itemize} 
	\end{enumerate}
\end{proposition}

\begin{proof}
	For part (1), we claim that the $3$-pseudofunctor $p'$ will \begin{itemize}
		\item Send an object $X$ to the object $p_{X}: FX \rightarrow GX$ in $[\mathfrak{P}_{1}, \mathfrak{B}]$,
		\item Have actions on hom-$2$-categories given by $2$-functors ${p'}_{X, Y}: \mathfrak{A}\left(X, Y\right) \rightarrow [\mathfrak{P}_{1}, \mathfrak{B}]\left(p_{X}, p_{Y}\right)$ which \begin{itemize}
			\item On objects will send $f$ to the object of the $2$-category $[\mathfrak{P}_{1}, \mathfrak{B}]\left(p_{X}, p_{Y}\right)$ determined by $\left( p_{f}\dashv p_{f}^{*}, Ff, Gf\right)$.
			\item On morphisms will send $\phi: f \rightarrow g$ to the morphism $\left(p_{\phi}, F\phi, G\phi\right)$ from $\left(p_{f}\dashv p_{f}^{*}, Ff, Gf\right)$ to $\left( p_{g}\dashv p_{g}^{*}, Fg, Gg\right)$ in the $2$-category $[\mathfrak{P}_{1}, \mathfrak{B}]\left(p_{X}, p_{Y}\right)$.
			\item On $2$-cells, will send $\Gamma: \phi \Rightarrow \psi$ to the $2$-cell $\left(F\Gamma, G\Gamma\right)$.
		\end{itemize}
		\item Have unitors $1_{p'X} \Rightarrow p'\left(1_{X}\right)$ given by the $2$-cells $\left(p^{X}, 1_{F1_{X}}, 1_{G1_{X}}\right)$ in $[\mathfrak{P}_{1},\mathfrak{B}]$.
		\item Have compositors $p'\left(g\right).p'\left(f\right) \Rightarrow p'\left(gf\right)$ given by the $2$-cells $\left(p_{g, f}, 1_{Fgf}, 1_{Ggf}\right)$ in $[\mathfrak{P}_{1},\mathfrak{B}]$.
	\end{itemize}

	\noindent We need to check that the hom-$2$-functors described above are well-defined.
	
	\begin{itemize}
		\item The fact that $\left(F\Gamma, G\Gamma\right)$ is well-defined as a $2$-cell in $[\mathfrak{P}_{1}, \mathfrak{B}]\left(p_{X}, p_{Y}\right)$ is precisely local pseudonaturality in $3$-cells for $p$, and the fact that this assignation is functorial is precisely functoriality between hom-categories of both $F$ and $G$.
		\item The fact that  $\phi \mapsto \left(p_{\phi}, F\phi, G\phi\right)$ is functorial between underlying categories is precisely due to the unit and composition preservation aspects of the pseudonatural transformation $\left(f,\phi\right) \mapsto \left(p_{f}, p_{\phi}\right)$.
		\item The fact that the $2$-functor respects horizontal composition of $2$-cells is precisely the same property for the corresponding hom-$2$-functors $F_{X, Y}:\mathfrak{A}\left(X, Y\right) \rightarrow \mathfrak{B}\left(FX, FY\right)$ and $G_{X, Y}: \mathfrak{A}\left(X, Y\right) \rightarrow \mathfrak{B}\left(GX, GY\right)$.
	\end{itemize}
	
	\noindent We need to check that we have described a well-defined $3$-pseudofunctor.
	
	\begin{itemize}
		\item 
		Naturality in $1$-cells of the compositor of the $3$-pseudofunctor $p'$ corresponds precisely to the left and right whiskering laws for the trinatural transformation $p$. In other words, it corresponds precisely to the modification axiom for the compositors of the trinatural transformation $p$. Naturality in $2$-cells $\Gamma: \phi \Rightarrow \psi$ follows from local pseudonaturality of $p$ in $\Gamma$, and uses the fact that $F$ and $G$ preserve interchangers in the case where $\Gamma$ is an interchanger $2$-cell in $\mathfrak{A}\left(Y, Z\right) \otimes \mathfrak{A}\left(X, Y\right)$.
		\item The fact that the left and right unit laws for the $3$-pseudofunctor $p'$ hold on the nose rather than up to invertible modifications correspond precisely to the left and right unit laws for the trinatural transformation $p$.
		\item The fact that the associativity condition for the $3$-pseudofunctor $p'$ holds on the nose rather than up to an invertible modification is precisely the associativity condition for  trinatural transformation $p$.
	\end{itemize}
	
	\noindent Conversely, observe that for any $3$-pseudofunctor $q': \mathfrak{A}\rightsquigarrow [\mathfrak{P}_{1}, \mathfrak{B}]$ the property of satisfying the lifting condition $\left(s_{0}, t_{0}\right). q' = \left(F, G\right)$ is tantamount to all of the following conditions holding.
	
	\begin{itemize}
		\item The object $q'\left(X\right)\in [\mathfrak{P}_{1}, \mathfrak{B}]$ being determined by a morphism in $\mathfrak{B}$ from $FX$ to $GX$.
		\item $q'\left(f: X \rightarrow Y\right)_{0} = Ff: FX \rightarrow FY$ and $q'\left(f: X \rightarrow Y\right)_{1} = Gf: GX \rightarrow GY$.
		\item $q'\left(\phi: f \Rightarrow g\right)_{0} = F\phi: Ff \Rightarrow Fg$ and $q'\left(\phi: f \Rightarrow g\right)_{1} = G\phi: Gf \Rightarrow Gg$.
		\item $q'\left(\Gamma: \phi \Rrightarrow \psi\right)_{0} = F\Gamma: F\phi \Rrightarrow F\psi$ and $q'\left(\Gamma: \phi \Rrightarrow \psi\right)_{1} = G\Gamma: G\phi \Rrightarrow G\psi$.
		\item The unitors and compositors of $q'$ being $2$-cells $\left(\alpha, \alpha_{0}, \alpha_{1}\right)$ in $[\mathfrak{P}_{1}, \mathfrak{B}]$ satisfying the condition that $\alpha_{0}$ and $\alpha_{1}$ are identity $2$-cells in $\mathfrak{B}$.
	\end{itemize}
	
	\noindent Notice that the $3$-pseudofunctor $p'$ which we defined indeed satisfies all of these conditions, and that any $3$-pseudofunctor $q': \mathfrak{A} \rightsquigarrow [\mathfrak{P}_{1}, \mathfrak{B}]$ which also satisfies these conditions gives rise to a trinatural transformation $q: F \rightarrow G$ since we have described a one to one correspondence between the data and axioms of trinatural transformations and $3$-pseudofunctors satisfying the lifting condition. Observe also that the unitors and compositors for the trinatural transformation $p$ are identities if and only if the corresponding data for the $3$-pseudofunctor $p'$ are identities. This completes the proof of part (1).
	\\
	\\
	\noindent Now let $q: F \rightarrow G$ be another trinatural transformation giving rise to another $q': \mathfrak{A} \rightarrow[\mathbf{2}, \mathfrak{B}]$ satisfying $\left(s_{0}, t_{0}\right).q' = \left(F, G\right)$. For part (2) we need to describe a bijection between $3$-pseudonatural transformations $\sigma': p' \Rightarrow q'$ satisfying $\left(s_{0}, t_{0}\right).\sigma' = 1_{\left(F, G\right)}$ and trimodificational adjoint equivalences $\sigma: p \Rightarrow q$. But starting with  $\big(\sigma, \left(\sigma_{X} \dashv \sigma_{X}^{*}\right)_{X \in \mathfrak{A}}\big)$, we claim that the $3$-pseudonatural transformation  $\big(\sigma, \left(\sigma_{X} \dashv \sigma_{X}^{*}\right)_{X \in \mathfrak{A}}\big)'$ will have
	
	\begin{itemize}
		\item $1$-cell component on $X$ given by the morphism $\left(\sigma_{X}\dashv \sigma_{X}^{*}, 1_{FX}, 1_{GX}\right)$ from $p_{X}$ to $q_{X}$ in $[\mathbf{2}, \mathfrak{B}]$.
		\item $2$-cell component on $f: X \rightarrow Y$ given by the $2$-cell $\left( \sigma_{f}, 1_{1_{FX}}, 1_{1_{GX}}\right)$ from $\left(\sigma_{X}\dashv \sigma_{X}^{*}, 1_{FX}, 1_{GX}\right)$ to $\left(\sigma_{Y}\dashv \sigma_{Y}^{*}, 1_{FY}, 1_{GY}\right)$ in $[\mathbf{2}, \mathfrak{B}]$.
	\end{itemize} 
	
	\noindent Now observe that naturality in $1$-cells of $f \mapsto  \left( \sigma_{f}, 1_{1_{FX}}, 1_{1_{FX}}\right)$ is precisely the local modification condition for $\sigma$, while $2$-naturality follows from local pseudonaturality of $p$ and $q$. The fact that the unit and composition conditions for $3$-pseudonaturality hold on the nose corresponds precisely to the unit and composition laws for the trimodification $\sigma$. Next, observe that for a $3$-pseudonatural transformation $\tau': p' \Rightarrow q'$, satisfaction of the lifting condition $\left(s_{0}, t_{0}\right).\tau' = 1_{\left(F, G\right)}$ is tantamount to the conditions \begin{itemize}
		\item $\tau'\left(X\right)_{0} = 1_{FX}$ and $\tau'\left(X\right)_{1} = 1_{GX}$
		\item $\tau'\left(f: X \rightarrow Y\right)_{0}$ and $\tau'\left(f: X \rightarrow Y\right)_{1}$ being identity $2$-cells in $\mathfrak{B}$.
	\end{itemize} 
	
	\noindent Notice that the $3$-pseudonatural transformation $\big(\sigma,  \left(\sigma_{X} \dashv \sigma_{X}^{*}\right)_{X \in \mathfrak{A}}\big)': p' \Rightarrow q'$ that we have described indeed satisfies all of these conditions. Conversely, any $3$-pseudonatural transformation $\tau': p' \Rightarrow q'$ which also satisfies these conditions gives rise to a trimodification $\tau: p \Rightarrow q$ by taking $\tau_{X}$ to be the $2$-cell component of ${\tau'}_{X}$ and $\tau_{f}$ to be the $3$-cell component of ${\tau'}_{f}$. Moreover, observe that by the description of morphisms in $[\mathbf{2}, \mathfrak{B}]$ given in Example \ref{coherent hom from 2 to B}, each ${\tau'}_{X}$ is the left adjoint part of a specified adjoint equivalence. We have hence described a one to one correspondence between the data and axioms of trimodificational equivalences and $3$-pseudonatural transformations satisfying the lifting conditions. This completes the proof part (2).
\end{proof}

\noindent In Theorem \ref{trinatural transformations as 2-Cat pseudofunctors} part (1) we have shown that $\mathbf{Gray}$-categories of the form $[\mathfrak{P}_{1}, \mathfrak{B}]$ can be used to model general trinatural transformations via semi-strict ones, using $3$-pseudofunctors. This clarifies the similarities between the data (i.e. compositors and unitors) and axioms (i.e. left and right whiskering, left and right unit, and associativity laws) for trinatural transformations and for pseudofunctors.
\\
\\
\noindent We will now describe the $\mathbf{Gray}$-category $[\mathfrak{P}_{2}, \mathfrak{B}]$ and the associated $\mathbf{Gray}$-functors $s_{1}, t_{1}: [\mathfrak{P}_{2}, \mathfrak{B}] \rightarrow [\mathfrak{P}_{1}, \mathfrak{B}]$. This facilitates a description of trimodifications $\sigma: p \Rightarrow q$ via  $3$-pseudofunctors $\sigma': \mathfrak{A} \rightarrow [\mathfrak{P}_{2}, \mathfrak{B}]$ satisfying $\left(s_{1}, t_{1}\right).\sigma' = \left(p', q'\right)$. 

\begin{example}\label{coherent hom from free 2-cell to B}
	The $\mathbf{Gray}$-category $[\mathfrak{P}_{2}, \mathfrak{B}]$ admits the following description. \begin{itemize}
		\item Its objects $\left(\mathbf{X}, x_{0}, x_{1}, X_{0}, X_{1}\right)$ are given by $2$-cells $\mathbf{X}: x_{0} \Rightarrow x_{1}: X_{0} \rightarrow X_{1}$ in $\mathfrak{B}$.
		\item Its morphisms $\left(\mathbf{f},  \bar{f}_{0} \dashv \bar{f}_{0}^{*}, \bar{f}_{1} \dashv \bar{f}_{1}^{*}, f_{0}, f_{1}\right)$ consist of \begin{itemize}
			\item For $i \in \{0, 1\}$ a morphism $f_{i}: X_{i} \rightarrow Y_{i}$ in $\mathfrak{B}$,
			\item For $i \in \{0, 1\}$, an adjoint equivalence in $\mathfrak{B}\left(X_{0}, Y_{1}\right)$ with left adjoint $\bar{f}_{i}: x_{1}.f_{i} \rightarrow f_{i}.x_{0}$, and right adjoint $f_{i}^{*}$.
			\item An invertible $3$-cell $\bar{f}$ as depicted below left.
		\end{itemize}
		
		\item Its $2$-cells $\left(\Phi_{0}, \Phi_{1}, \phi_{0}, \phi_{1} \right):\left(\mathbf{f}, \bar{f}_{0} \dashv \bar{f}_{0}^{*}, \bar{f}_{1} \dashv \bar{f}_{1}^{*}, f_{0}, f_{1}\right) \Rightarrow \left(\mathbf{g}, \bar{g}_{0} \dashv \bar{g}_{0}^{*}, \bar{g}_{1} \dashv \bar{g}_{1}^{*}, g_{0}, g_{1}\right)$ consist of \begin{itemize}
			\item For $i \in \{0, 1\}$, a $2$-cell $\phi_{i}: f_{i} \Rightarrow g_{i}$, 
			\item For $i \in \{0, 1\}$, an invertible $3$-cell $\Phi_{i}$ as depicted below right.
		\end{itemize}
		
		$$\begin{tikzcd}[font=\fontsize{9}{6}]
			y_{0}.f_{0} \arrow[rr, "\mathbf{Y}.1"]\arrow[dd, "\bar{f}_{0}"'] &{}\arrow[dd, Rightarrow, shorten = 10, "\mathbf{f}"]& y_{1}.f_{0}\arrow[dd, "\bar{f}_{1}"]
			\\
			\\
			f_{1}.x_{0} \arrow[rr, "1.\mathbf{X}"'] &{}& f_{1}.x_{1} &{}
		\end{tikzcd}\begin{tikzcd}[font=\fontsize{9}{6}]
			f_{1}.x_{i} \arrow[rr, "\phi_{1}.1"]\arrow[dd, "\bar{f}_{i}"'] &{}\arrow[dd, Rightarrow, shorten = 10, "\Phi_{i}"]& g_{1}.x_{i}\arrow[dd, "\bar{g}_{i}"]
			\\
			\\
			y_{i}.f_{0} \arrow[rr, "1.\phi_{0}"'] &{}& y_{i}.g_{0} &{}
		\end{tikzcd}$$
		
		\noindent Such that the equation below holds in the hom-$2$-category $\mathfrak{B}\left(X_{0}, Y_{1}\right)$.

		$$\begin{tikzcd}[column sep = 12, font=\fontsize{9}{6}]
			&&f_{1}.x_{1}
			\arrow[lldd,leftarrow, "\phi_{1}.1"']
			\arrow[rrdd, "1.\mathbf{X}" description]
			\arrow[dddd, Rightarrow, "{\left(\phi_{1}\right)}_{\left(\mathbf{X}\right)}", shorten = 30, shift right = 6]
			\arrow[rr, "\bar{f}_{1}"]
			&&y_{1}.f_{0}
			\arrow[dd, Rightarrow, "{\Phi}_{1}", shorten = 10]
			\arrow[rrdd, "1.\phi_{0}"]
			\\
			\\
			f_{1}.x_{0}\arrow[rrdd, "\phi_{1}.1"']
			&&&&g_{1}.x_{1}
			\arrow[lldd, leftarrow, "1.\mathbf{X}" description]
			\arrow[rr, "\bar{g}_{1}"]
			\arrow[dd, Rightarrow, "\mathbf{g}", shorten = 10]
			&&y_{1}.g_{0}\arrow[lldd, leftarrow, "\mathbf{Y}.1"]&=
			\\
			\\
			&&g_{1}.x_{0}
			\arrow[rr, "\bar{g}_{0}"']
			&&y_{0}.g_{0}
		\end{tikzcd}\begin{tikzcd}[column sep = 12, font=\fontsize{9}{6}]
			&&f_{1}.x_{1}
			\arrow[lldd,leftarrow, "\phi_{1}.1"']
			\arrow[dd, Rightarrow, "\mathbf{f}"', shorten = 10]
			\arrow[rr, "\bar{f}_{1}"]
			&&y_{1}.f_{0}
			\arrow[lldd, leftarrow, "\mathbf{Y}.1" description]
			\arrow[dddd, Rightarrow, "{\left(\mathbf{Y}\right)}_{\left(\phi_{0}\right)}"', shorten = 30, shift left = 6]
			\arrow[rrdd, "1.\phi_{0}"]
			\\
			\\
			f_{1}.x_{0}
			\arrow[rrdd, "1.\mathbf{X}"']
			\arrow[rr, "\bar{f}_{0}"]
			&&y_{0}.f_{0}
			\arrow[dd, Rightarrow, "\Phi_{0}"', shorten = 10]
			\arrow[rrdd, "1.\phi_{0}" description]
			&&
			&&y_{1}.g_{0}
			\arrow[lldd, leftarrow, "\mathbf{Y}.1"]
			\\
			\\
			&&g_{1}.x_{0}
			\arrow[rr, "\bar{g}_{0}"']
			&&y_{0}.g_{0}
		\end{tikzcd}$$
		
		\item Its $3$-cells $\left(\Gamma_{0}, \Gamma_{1}\right): \left(\Phi_{0}, \Phi_{1}, \phi_{0}, \phi_{1}\right) \Rrightarrow \left(\Psi_{0}, \Psi_{1}, \psi_{0}, \psi_{1}\right)$ consist of two $3$-cells $\Gamma_{0}:\phi_{0} \Rrightarrow \psi_{0}$ and $\Gamma_{1}: \phi_{1} \Rrightarrow \psi_{1}$ such that the equation depicted below holds in $\mathfrak{B}\left(X_{0}, Y_{1}\right)$ for $i \in \{0, 1\}$.

		$$\begin{tikzcd}[font=\fontsize{9}{6}]
			y.f_{0} \arrow[rr, bend right, "1.\psi_{0}"' name = B]\arrow[rr, bend left, "\phi_{0}" name = A]
			\arrow[dd, "\bar{f}_{i}"'] &{}& y.g_{0}\arrow[dd, "\bar{g}_{i}"]
			\\
			&&&=
			\\
			f_{1}.x \arrow[rr, bend right, "\psi_{1}.1"' name = C] &{}& g_{1}.x 
			\arrow[from= A, to =B, Rightarrow, "\Gamma_{0}", shorten = 10]\arrow[from =B, to =C, Rightarrow, shorten = 15, "\Psi_{i}"]
		\end{tikzcd}\begin{tikzcd}[row sep = 25, font=\fontsize{9}{6}]
			y.f_{0} \arrow[rr, bend left, "\phi_{0}" name = A]
			\arrow[dd, "\bar{f}_{i}"'] &{}& y.g_{0}\arrow[dd, "\bar{g}_{i}"]
			\\
			\\
			f_{1}.x \arrow[rr, bend right, "\psi_{1}.1"' name = C] \arrow[rr, bend left, "\phi_{1}.1" name = D]&{}& g_{1}.x &{}
			\arrow[from= A, to =D, Rightarrow, shorten = 12, "\Phi_{i}"]\arrow[from =D, to =C, Rightarrow, shorten = 10, "\Gamma_{1}"]
		\end{tikzcd}$$
		
	\end{itemize}
	\noindent There are $\mathbf{Gray}$-functors $s_{1}, t_{1}: [\mathfrak{P}_{2}, \mathfrak{B}] \rightarrow [\mathfrak{P}_{1}, \mathfrak{B}]$ which take sources and targets of all data in the following way.
	
	\begin{align*}
		s_{1}\left(\mathbf{X}: x_{0} \Rightarrow x_{1}: X_{0} \rightarrow X_{1}\right) = \left(x_{0}, X_{0}, X_{1}\right) && t_{1}\left(\mathbf{X}: x_{0} \Rightarrow x_{1}: X_{0} \rightarrow X_{1}\right)  = \left(x_{1}, X_{0}, X_{1}\right)
		\\
		s_{1}\left(\mathbf{f},  \bar{f}_{0} \dashv \bar{f}_{0}^{*}, \bar{f}_{1} \dashv \bar{f}_{1}^{*}, f_{0}, f_{1}\right)  = \left(\bar{f}_{0} \dashv \bar{f}_{0}^{*}, f_{0}, f_{1}\right) && 	t_{1}\left(\mathbf{f},  \bar{f}_{0} \dashv \bar{f}_{0}^{*}, \bar{f}_{1} \dashv \bar{f}_{1}^{*}, f_{0}, f_{1}\right)  = \left(\bar{f}_{1} \dashv \bar{f}_{1}^{*}, f_{0}, f_{1}\right)
		\\
		s_{1}\left(\Phi_{0},\Phi_{1}, \phi_{0}, \phi_{1}\right) = \left(\Phi_{0}, \phi_{0}, \phi_{1}\right) && 
		t_{1}\left(\Phi_{0},\Phi_{1}, \phi_{0}, \phi_{1}\right) = \left(\Phi_{1}, \phi_{0}, \phi_{1}\right) 
		\\
		s_{1}\left(\Gamma_{0}, \Gamma_{1}\right) =\left(\Gamma_{0}, \Gamma_{1}\right) && t_{1}\left(\Gamma_{0}, \Gamma_{1}\right) =\left(\Gamma_{0}, \Gamma_{1}\right)
	\end{align*}
\end{example}

\begin{proposition}\label{trimodifications via 3-pseudofunctors}
	Let $p, q: F \rightarrow G: \mathfrak{A} \rightarrow \mathfrak{B}$ be trinatural transformations between $\mathbf{Gray}$-functors and let $p', q': \mathfrak{A} \rightarrow [\mathfrak{P}_{1}, \mathfrak{B}]$ be the corresponding $3$-pseudofunctors. There is a bijection between \begin{itemize}
			\item the set of trimodifications $\sigma: p \Rightarrow q$,
			\item the set of $3$-pseudofunctors $\sigma': \mathfrak{A} \rightarrow [\mathfrak{P}_{2}, \mathfrak{B}]$ which satisfy $\left(s_{1}, t_{1}\right).\sigma' = \left(p', q'\right)$.
		\end{itemize}
\end{proposition}

\begin{proof}
	Let $\sigma: p \Rightarrow q$ be a trimodification and define the $3$-pseudofunctor $\sigma': \mathfrak{A} \rightarrow [\mathfrak{P}_{2}, \mathfrak{B}]$ in the following way.
	
	\begin{itemize}
		\item It sends an object $X \in \mathfrak{A}$ to the $2$-cell $\sigma_{X}: p_{X} \Rightarrow q_{X}$ in $\mathfrak{B}$,
		\item Between hom-$2$-categories, it sends \begin{itemize}
			\item $f: X \rightarrow Y$ in $\mathfrak{A}$ to the morphism $\left(\sigma_{f}, p_{f}\dashv p_{f}^{*}, q_{f}\dashv q_{f}^{*}, Ff, Gf\right)$ from $\sigma_{X}$ to $\sigma_{Y}$ in $[\mathfrak{P}_{2}, \mathfrak{B}]$.
			\item $\phi: f \Rightarrow g: X \rightarrow Y$ in $\mathfrak{A}$ to the $2$-cell $\left(q_{\phi}, p_{\phi}, F\phi, G\phi\right)$ from $\left(\sigma_{f}, p_{f}\dashv p_{f}^{*}, q_{f}\dashv q_{f}^{*}, Ff, Gf\right)$ to $\left(\sigma_{g}, p_{g}\dashv p_{g}^{*}, q_{g}\dashv q_{g}^{*}, Fg, Gg\right)$ in $[\mathfrak{P}_{2}, \mathfrak{B}]$.
			\item $\Gamma: \phi \Rrightarrow \psi: f \Rightarrow g: X \rightarrow Y$ in $\mathfrak{A}$ to the $3$-cell $\left(G\Gamma, F\Gamma\right)$ from  $\left(q_{\phi}, p_{\phi}, F\phi, G\phi\right)$ to  $\left(q_{\psi}, p_{\psi}, F\psi, G\psi\right)$ in $[\mathfrak{P}_{2}, \mathfrak{B}]$.
		\end{itemize}
		\item Its unitor's component on $X$ given by the $2$-cell $\left(p^{X}, q^{X}, 1_{1_{FX}}, 1_{1_{GX}}\right)$ in $[\mathfrak{P}_{2}, \mathfrak{B}]$ from the identity on $\sigma_{X}$ to $\left(\sigma_{1_{X}}, p_{1_{X}}, q_{1_{X}}, F1_{X}, G1_{X}\right)$.
		\item It has compositor component on \begin{tikzcd}[column sep = 16]
			X \arrow[r, "f"] & Y \arrow[r, "g"] & Z
		\end{tikzcd} given by the $2$-cell $\left(p_{g, f}, q_{g, f}, 1_{Fgf}, 1_{Ggf}\right)$.
	\end{itemize}
	\noindent We need to check that the hom-$2$-functors $\mathfrak{A}\left(X, Y\right) \rightarrow [\mathfrak{P}_{2}, \mathfrak{B}]\left(\sigma_{X}, \sigma_{Y}\right)$ are well-defined. The equation for $2$-cells is precisely the local modification axiom for $\sigma$, while the equation for $3$-cells is precisely local pseudonaturality of $p$. Functoriality between hom-categories is precisely that same property for $F$ and $G$. Functoriality between underlying categories of $\mathfrak{A}\left(X, Y\right) \rightarrow [\mathfrak{P}_{2}, \mathfrak{B}]\left(\sigma_{X}, \sigma_{Y}\right)$ is precisely the unit and composition laws for $\left(f, \phi\right)\mapsto \left(p_{f}, p_{\phi}\right)$ and $\left(f, \phi\right)\mapsto \left(q_{f}, q_{\phi}\right)$. Preservation of horizontal composition follows from that same property for the actions between hom-$2$-categories of $F$ and $G$. Thus the hom-$2$-functors are well-defined. We now check that we have described a well-defined $3$-pseudofunctor. Naturality in $1$-cells of the compositor of $\sigma'$ is precisely the left and right whiskering laws for the compositors of $p$ and $q$, while $2$-naturality on $\Gamma: \phi \Rightarrow \psi$ follows from local pseudonaturality of $p$ and $q$ in the corresponding $3$-cell $\Gamma$, again using $F$ and $G$'s preservation of interchangers when $\Gamma$ is an interchanger $2$-cell in $\mathfrak{A}\left(Y, Z\right) \otimes \mathfrak{A}\left(X, Y\right)$. Similarly, the left and right unit and associativity axioms for the $3$-pseudofunctor coincide with corresponding the trinaturality axioms for $p$ and $q$.
	\\
	\\
	\noindent Conversely, a $3$-pseudofunctor $\tau': \mathfrak{A} \rightsquigarrow [\mathfrak{P}_{2}, \mathfrak{B}]$ satisfies $\left(s_{1}, t_{1}\right).\tau' = \left(p', q'\right)$ precisely if the following conditions hold
	
	\begin{itemize}
		\item For all $X \in \mathfrak{A}$, the $2$-cell $\tau'\left(X\right)$ has source and target $1$-cells given by $p_{X}$ and $q_{X}$ respectively.
		\item For all morphisms $f: X \rightarrow Y$ in $\mathfrak{A}$, the morphism $\tau'\left(f\right) = \left(\mathbf{f}, \bar{f}_{1} \dashv \bar{f}_{1}^{*},  \bar{f}_{0} \dashv \bar{f}_{0}^{*}, f_{0}, f_{1}\right)$ is such that $f_{0}= Ff$, $f_{1}= Gf$, $\bar{f}_{0}\dashv\bar{f}_{0}^{*} = p_{f}\dashv p_{f}^{*}$ and  $\bar{f}_{1}\dashv\bar{f}_{1}^{*} = q_{f}\dashv q_{f}^{*}$.
		\item For all $2$-cells $\phi: f \Rightarrow g$ in $\mathfrak{A}$, $\tau'\left(\phi\right) = \left(p_{\phi}, q_{\phi}, F\phi, G\phi\right)$.
		\item For all $3$-cells $\Gamma: \phi \Rrightarrow \psi$ in $\mathfrak{A}$, $\tau'\left(\Gamma\right) = \left(F\Gamma, G\Gamma\right)$.
		\item For all objects $X \in \mathfrak{A}$, the unitor $\left(\Phi_{0}, \Phi_{1}, \phi_{0}, \phi_{1}\right): 1_{\tau'\left(X\right)} \Rightarrow \tau'\left(1_{X}\right)$ is such that $\phi_{0}$ and $\phi_{1}$ are identities, $\Phi_{0} = p^{X}$ and $\Phi_{1} = q^{X}$.
		\item For all \begin{tikzcd}
			X \arrow[r, "f"] & Y \arrow[r, "g"] & Z
		\end{tikzcd} in $\mathfrak{A}$, the compositor $\left(\Phi_{0}, \Phi_{1}, \phi_{0}, \phi_{1}\right): \tau'\left(g\right).\tau'\left(f\right) \Rightarrow \tau'\left(gf\right)$ is such that $\phi_{0}$ and $\phi_{1}$ are identities, $\Phi_{0} = p_{g, f}$ and $\Phi_{1} = q_{g, f}$.
	\end{itemize}
	\noindent And in this case $\left(X, f\right) \mapsto \left(\tau_{X}, \tau_{f}\right)$ is well-defined as a trimodification due to the corresponding properties of the $3$-pseudofunctor, with naturality in $1$-cells of the compositor for $\tau'$ corresponding to the local modification axiom for $\tau$ and $2$-naturality in $\Gamma: \phi: \Rightarrow \psi$ following from local pseudonaturality of $p$ and $q$ in $\Gamma$, and using $\mathbf{Gray}$-functoriality of $F$ and $G$ when $\Gamma$ is an interchanger in $\mathfrak{A}\left(Y, Z\right)\otimes \mathfrak{A}\left(X,Y\right)$. It is also clear that these constructions are inverse to one another. This completes the proof.  
\end{proof}

\noindent An explicit description of the $\mathbf{Gray}$-category $[\mathfrak{P}_{3}, \mathfrak{B}]$ can also be given, similar to that in Examples \ref{coherent hom from 2 to B} and \ref{coherent hom from free 2-cell to B}. Objects consist of a $3$-cell in $\mathfrak{B}$, morphisms have data subject to one axiom, while $2$-cells and $3$-cells have data subject to two axioms in both cases. There are again source and target $\mathbf{Gray}$-functors  $s_{2}, t_{2}: [\mathfrak{P}_{3}, \mathfrak{B}] \rightarrow [\mathfrak{P}_{2}, \mathfrak{B}]$. Then, given $\sigma, \tau: p \Rightarrow q: F \rightarrow G$ there is a bijection between $3$-pseudofunctors $\Omega': \mathfrak{A} \rightsquigarrow [\mathfrak{P}_{3}, \mathfrak{B}]$ satisfying $\left(s_{2}, t_{2}\right).\Omega'$ and perturbations $\Omega: \sigma \Rrightarrow \tau$. The proof of this claim is similar to the proofs of Propositions \ref{trinatural transformations as 2-Cat pseudofunctors} and \ref{trimodifications via 3-pseudofunctors}, so we omit it. We collect the results of this section into Theorem \ref{2-globular set modelling weak higher cells via 3 pseudofunctors}. Proposition \ref{internal category with weak composition} extends the globular information of Theorem \ref{2-globular set modelling weak higher cells via 3 pseudofunctors} to include the underlying category structure of $\mathbf{Tricat}_{s}\left(\mathfrak{A}, \mathfrak{B}\right)$.

\begin{theorem}\label{2-globular set modelling weak higher cells via 3 pseudofunctors}
	Let $\mathfrak{B}$ be a $\mathbf{Gray}$-category and form the $3$-globular $\mathbf{Gray}$-category $\bar{\mathfrak{B}}$ displayed below. Let $F, G: \mathfrak{A} \rightarrow \mathfrak{B}$ be $\mathbf{Gray}$-functors and let $\left(F, G\right): \mathfrak{A} \rightarrow \mathfrak{B} \times \mathfrak{B}$ denote the corresponding $2$-functor induced by the universal property of the product. This $3$-globular $\mathbf{Gray}$-category classifies $n$-cells in $\mathbf{Tricat}_{s}[\mathfrak{A}, \mathfrak{B}]\left(F, G\right)$ via $3$-pseudofunctors $\mathfrak{A} \rightarrow \bar{\mathfrak{B}}_{n+1}$. 
	
	$$\begin{tikzcd}
		{[}\mathfrak{P}_{3}{,}\mathfrak{B}{]} \arrow[rr, shift left, "s"] \arrow[rr, shift right, "t"'] && {[}\mathfrak{P}_{2}{,}\mathfrak{B}{]} \arrow[rr, shift left, "s"] \arrow[rr, shift right, "t"'] &&  {[}\mathfrak{P}_{1}{,}\mathfrak{B}{]} \arrow[rr, shift left, "s"] \arrow[rr, shift right, "t"'] &&\mathfrak{B}
	\end{tikzcd}$$
	
\end{theorem}

\noindent In Proposition \ref{internal category with weak composition}, to follow, describes composites of semi-strict trinatural transformations via $3$-pseudofunctors. See Chapter 8 of \cite{Borceaux Vol 1} for background on internal categories.

\begin{proposition}\label{internal category with weak composition}
	\begin{enumerate}
		\item There is a $3$-pseudofunctor $m: [\mathbf{2}, \mathfrak{B}]\times_\mathfrak{B} [\mathbf{2}, \mathfrak{B}] \rightsquigarrow [\mathbf{3}, \mathfrak{B}]$ with identity unitor and with compositor on the pair of composable pairs depicted below given by the interchanger of ${g'}_{f}$ in $\mathfrak{B}$, whiskered with $g$ and $f'$.	
		
		$$\begin{tikzcd}[row sep = 5]
			\left(x_{0}{,} X_{0}{,} X_{1}\right)
			\arrow[rr, "\left(f \dashv f^{*}{,} f_{0}{,} f_{1}\right)"]
			&&
			\left(y_{0}{,} Y_{0}{,} Y_{1}\right)
			\arrow[rr, "\left(g \dashv g^{*}{,} g_{0}{,} g_{1}\right)"]
			&&
			\left(z_{0}{,} Z_{0}{,} Z_{1}\right)
			\\
			\left(x_{1}{,} X_{1}{,} X_{2}\right)
			\arrow[rr, "\left(f{'} \dashv f^{*}{'}{,} f_{0}{'}{,} f_{1}{'}\right)"']
			&&
			\left(y_{1}{,} Y_{1}{,} Y_{2}\right)
			\arrow[rr, "\left(g{'} \dashv g^{*}{'}{,} g_{0}{'}{,} g_{1}{'}\right)"']
			&&
			\left(z_{1}{,} Z_{1}{,} Z_{2}\right)
		\end{tikzcd}$$
		\item $m$ is biequivalence inverse in $\mathbf{Gray}$-$\mathbf{Cat}_{3}$ to the canonical $\mathbf{Gray}$-functor to the pullback.
		\item The structure given by $\Delta_{\leq 3}^\text{op}\ni\mathbf{n}\mapsto {[}\mathbf{n}{,}\mathfrak{B}{]}$ satisfies the axioms for a category internal to $\mathbf{Gray}$-$\mathbf{Cat}$, except for the fact that the composition map depicted below is a $3$-pseudofunctor rather than a $\mathbf{Gray}$-functor.
		
		$$\begin{tikzcd}
			{[}\mathbf{2}{,}\mathfrak{B}{]}\times_{\mathfrak{B}} {[}\mathbf{2}{,}\mathfrak{B}{]} \arrow[r, "m"] &{[}\mathbf{3}{,}\mathfrak{B}{]} \arrow[r, "{[}d_{1}{,}\mathfrak{B}{]}"]&{[}\mathbf{2}{,}\mathfrak{B}{]}
		\end{tikzcd}$$ 
	\end{enumerate}
	
\end{proposition}

\begin{proof}
	For part (1), $m$ sends the pair of $2$-cells consisting of \begin{tikzcd}
		\left(f \dashv f^{*}{,} f_{0}{,} f_{1}\right) \arrow[rr,"\left(\phi{,} \phi_{0}{,} \phi_{1}\right)"] && \left(g \dashv g^{*}{,} g_{0}{,} g_{1}\right) \end{tikzcd} and \begin{tikzcd}\left(g \dashv g^{*}{,} g_{0}{,} g_{1}\right)\arrow[rr, "\left(\phi'{,}\phi_{1}{,}\phi_{2}\right)"] && \left(h \dashv h^{*}{,} h_{0}{,} h_{1}\right)
	\end{tikzcd} to the semi-strict trinatural transformation whose components on generators are given by $\phi$ and $\phi'$. The behaviour on $3$-cells is similar, and the axioms for a $3$-pseudofunctor are clear from this description. Part (2) may be checked via direct calculation. For part (3), the unit and associativity laws follow from those in $\mathfrak{B}$.
\end{proof}

\subsection{Extending strictification to $\left(3, k\right)$-transfors}\label{Subsection Extending strictification to 3 k transfors}

\noindent We are now ready to extend $\mathbf{Gr}$ to $\left(3, k\right)$-transfors for positive $k$, and exhibit $\eta_\mathfrak{A}: \mathfrak{A} \rightsquigarrow \mathbf{Gr}\left(\mathfrak{A}\right)$ as a universal trihomomorphism into a $\mathbf{Gray}$-category. There are subtleties involved in forming a $\mathbf{Gray}$-category out of the resulting semi-strict $\left(3, k\right)$-transfors, which we will need to do in order to state the universal property. These subtleties will be addressed in Section \ref{Section Solutios to the Failure of Semi-strict trinatural transformations to be closed under composition}.
\\
\\
\noindent We first recall Lemma 15.4 of \cite{Gurski Coherence in Three Dimensional Category Theory} in part (1) of Lemma \ref{Gurski Lemma semi strict}, and extend it to trimodifications and perturbations in parts (2) and (3). This extends Proposition \ref{trihomomorphism out of cofibrant Gray category} to $\left(3, k\right)$-transfors for positive $k$.

\begin{lemma}\label{Gurski Lemma semi strict}
	Let $\mathfrak{A}$ be a cofibrant $\mathbf{Gray}$-category and let $F, G: \mathfrak{A} \rightarrow \mathfrak{B}$ be $\mathbf{Gray}$-functors.
	
	\begin{enumerate}
		\item 
		To (uniquely) define a semi-strict trinatural transformation $p: F \rightarrow G$ it suffices to \begin{itemize}
			\item define $1$-cell components on objects,
			\item define adjoint equivalence $2$-cell components on generating morphisms
			\item define $3$-cell components on generating $2$-cells
			\item show that the $3$-cell components satisfy local pseudonaturality with respect to the $3$-cells of $\mathfrak{A}$.
		\end{itemize}
		\item If $p, q: F \rightarrow G$ are semi-strict trinatural transformations, to (uniquely) define a trimodification $\sigma: p \Rightarrow q$ it suffices to \begin{itemize}
			\item define $2$-cell components on objects,
			\item define $3$-cell components on generating morphisms
			\item show that these data satisfy the local modification condition for generating $2$-cells in $\mathfrak{A}$.
		\end{itemize}
		\item If $\sigma, \tau: p \Rightarrow q$ are trimodifications, to (uniquely) define a perturbation $\Omega: \sigma \Rrightarrow \tau$ it suffices to define its $3$-cell components $\Omega_{X}: \sigma_{X}\Rrightarrow \tau_{X}$ and check the perturbation axiom for all generating morphisms in $\mathfrak{A}$.
	\end{enumerate}
\end{lemma}

\begin{proof}
	Part (1) is Lemma 15.4 of \cite{Gurski Coherence in Three Dimensional Category Theory}. It follows from Theorem \ref{trinatural transformations as 2-Cat pseudofunctors}, given that the underlying sesquicategory of $\mathfrak{A}$ is free on a $2$-computad. In particular, a $\mathbf{Gray}$-functor $p': \mathfrak{A} \rightarrow [\mathfrak{P}_{1}, \mathfrak{B}]$ is determined by its outputs on the generators in the underlying $2$-computad, as well as its outputs on $3$-cells. But the outputs on $3$-cells are already determined by the requirements that $s.p' = F$ and $t.p' = G$, so the only condition needing to be checked is that $p'\left(\Omega\right)$ is a well-defined $3$-cell in $[\mathfrak{P}_{1}, \mathfrak{B}]$. This corresponds to the fourth dot point. The proof for part (2) is analogous, using Proposition \ref{trimodifications via 3-pseudofunctors}. Once again, the behaviour on $3$-cells is determined by commutativity with $s$ and $t$, and well-definedness of $\sigma'\left(\phi\right)$ as a $2$-cell in $[\mathfrak{P}_{2}, \mathfrak{B}]$ precisely corresponds to the local modification condition for $\sigma$. The proof for part (3) is also analogous, and uses $[\mathfrak{P}_{3}, \mathfrak{B}]$.
\end{proof}

\begin{proposition}\label{Trinatural transformations out of cofibrant Gray categories}
	Let $F, G: \mathfrak{A} \rightarrow \mathfrak{B}$ be $\mathbf{Gray}$-functors with $\mathfrak{A}$ cofibrant. Let $p: F \rightarrow G$ be a trinatural transformation. Then there is \begin{enumerate}
		\item A unique semi-strict trinatural transformation $\bar{p}: F \rightarrow G$ agreeing with $p$ on objects and generating morphisms, and defined on $2$-cells using the unitors and compositors of $p$.
		\item A unique invertible costrict trimodification $i: \bar{p} \rightarrow p$ whose $3$-cell components at generating morphisms are identities. 
	\end{enumerate}
\end{proposition}

\begin{proof}
	 View the trinatural transformation $p: F \rightarrow G$ as a $3$-pseudofunctor $p': \mathfrak{A} \rightsquigarrow [\mathbf{2}, \mathfrak{B}]$ using Proposition \ref{trinatural transformations as 2-Cat pseudofunctors} part (1). Now use cofibrancy of $\mathfrak{A}$, the fact that $3$-pseudofunctors are certain trihomomorphisms, and Proposition \ref{trihomomorphism out of cofibrant Gray category} to construct a $\mathbf{Gray}$-functor $\underline{p'}: \mathfrak{A} \rightarrow [\mathbf{2}, \mathfrak{B}]$ and a pseudo-icon equivalence $\underline{e}: \underline{p'} \rightarrow p'$. But then observe that $\underline{e}$ is in fact a $3$-pseudonatural transformation. This is because as per Proposition \ref{trihomomorphism out of cofibrant Gray category}, both its operational and free $3$-cell components are given in terms of coherence $3$-cells of the trihomomorphism $p'$, and these $3$-cells are all identities in $[\mathbf{2}, \mathfrak{B}]$ since $p'$ is a $3$-pseudofunctor. Observe also that $s\underline{e} = {1}_{p'}$ and $t\underline{e} = {1}_{\underline{p'}}$, so that Proposition \ref{trinatural transformations as 2-Cat pseudofunctors} part (2) gives the desired trimodification $e: \bar{p}\Rightarrow p$, where the semi-strict trinatural transformation $\bar{p}: F \rightarrow G$ corresponds to the $\mathbf{Gray}$-functor $\underline{p'}: \mathfrak{A} \rightarrow [\mathbf{2}, \mathfrak{B}]$.
\end{proof}

\begin{theorem}\label{Gr on weak higher cells between tricategories}
	Let $F, G: \mathfrak{A} \rightarrow \mathfrak{B}$ be trihomomorphisms between tricategories. Let $\mathbf{Gr}\left(F\right), \mathbf{Gr}\left(G\right): \mathbf{Gr}\left(\mathfrak{A}\right) \rightarrow \mathbf{Gr}\left(\mathfrak{B}\right)$ be the corresponding $\mathbf{Gray}$-functors described in Section 10.6 of \cite{Gurski Coherence in Three Dimensional Category Theory}. Let $p: F \rightarrow G$ be a trinatural transformation.
	\begin{enumerate}
		\item There is a unique semi-strict trinatural transformation $\mathbf{Gr}\left(p\right): \mathbf{Gr}\left(F\right) \rightarrow \mathbf{Gr}\left(G\right)$ whose
		\begin{itemize}
			\item $1$-cell component on an object $X$ is given by $p_{X}$,
			\item $2$-cell component on a path of length one $f: X \rightarrow Y$ given by $p_{f}: \left(p_{Y}, Ff\right) \Rightarrow \left(Gf, p_{X}\right)$,
			\item $3$-cell component on a generating $2$-cell $\left(\mathfrak{p}_{0}: [f]_\text{min} \Rightarrow [f]_{\alpha_{0}}, \mathfrak{p}_{1}: [g]_{\alpha_{1}} \Rightarrow [g]_\text{min}\right): [f] \Rightarrow [g]$ is given by the pasting of $p_{\phi}$ with various compositors of $p$.
		\end{itemize} 
		\item Let $q: F \rightarrow G$ be another trinatural transformation and let $\sigma: p \Rightarrow q$ be a trimodification. Then there is a trimodification $\mathbf{Gr}\left(\sigma\right): \mathbf{Gr}\left(p\right) \Rightarrow \mathbf{Gr}\left(q\right)$ whose $2$-cell component on $X$ is given by $\sigma_{X}$ and whose $3$-cell component on a generating morphism $f: X \rightarrow Y$ is given by $\sigma_{f}$.
		\item Let $\tau: p \Rightarrow q$ be another trimodification and let $\Omega: \sigma \Rrightarrow \tau$ be a perturbation. Then there is a perturbation $\mathbf{Gr}\left(\Omega\right): \mathbf{Gr}\left(\sigma\right) \Rrightarrow \mathbf{Gr}\left(\tau\right)$ whose component on $X$ is given by $\Omega_{X}$.
	\end{enumerate} 
\end{theorem}

\begin{proof}
	\noindent Part (1) can be proved by combining Proposition \ref{Trinatural transformations out of cofibrant Gray categories} with Theorem 10.16 of \cite{Gurski Coherence in Three Dimensional Category Theory}. In particular, the desired $\mathbf{Gr}\left(p\right): \mathbf{Gr}\left(F\right) \rightarrow \mathbf{Gr}\left(G\right)$ is the one given by Proposition \ref{Trinatural transformations out of cofibrant Gray categories} from the trinatural transformation displayed below. Here $\psi$ is as in Theorem 10.16 of \cite{Gurski Coherence in Three Dimensional Category Theory}. Parts (2) and (3) follow similarly. Their proofs respectively use Proposition \ref{trimodifications via 3-pseudofunctors} part (1) and the analogous classifying property of $[\mathfrak{P}_{3}, \mathfrak{B}]$ with respect to perturbations, as part of Theorem \ref{2-globular set modelling weak higher cells via 3 pseudofunctors}.
	
	$$\begin{tikzcd}[column sep = 18]
		\mathbf{Gr}\left(F\right) \arrow[r, "1.\simeq"] & \mathbf{Gr}\left(F\right).\eta_\mathfrak{A}.\eta_{\mathfrak{A}}^{*}
		\arrow[r, "\psi_{F}.1"] & \left(\eta_\mathfrak{A}.F\right).{\eta}_{\mathfrak{A}}^{*}
		\arrow[r, "\left(1.p\right).1"] & 
		\left(\eta_\mathfrak{A}.G\right).{\eta}_\mathfrak{A}^{*}
		\arrow[r, "\psi_{G}^{-1}"] &  
		\mathbf{Gr}\left(G\right).\eta_\mathfrak{A}.\eta_{\mathfrak{A}}^{*}
		\arrow[r, "1.\simeq"] & 
		\mathbf{Gr}\left(G\right) 
	\end{tikzcd}$$
\end{proof}

\noindent We end this section by reviewing the low dimensional categorical structures formed by tricategories, which were introduced in \cite{Low dimensional structures formed by tricategories}. Later, Theorem \ref{trihomomorphism classifier} we will give an account of three dimensional strictification which will not be restricted to pseudo-icon equivalences, costrict trimodifications and identity perturbations. However, restricting one's attention in this way allows for a simpler formulation of semi-strictification in dimension three. In particular, Proposition \ref{tricategory of tricategories} gives an entirely three-dimensional account of three dimensional semi-strictification. This is in the same spirit as Remark \ref{Two Dimensional Statement of Two Dimensional Strictification}, although the statement is now tricategorical rather than $\mathbf{Gray}$-categorical.

\begin{proposition}\label{tricategory of tricategories}
	\begin{enumerate}
		\item There is a tricategory $\mathbf{Tricat}_{3}$ whose
		\begin{itemize}
			\item Objects are tricategories,
			\item Morphisms are trihomomorphisms,
			\item $2$-cells are pseudo-icon equivalences,
			\item $3$-cells are costrict trimodifications,
			\item Composition of $2$-cells and $3$-cells is defined so that they retain their respective identity components.
		\end{itemize} 
		\item The sub-tricategory of $\mathbf{Tricat}_{3}$ on $\mathbf{Gray}$-categories and $\mathbf{Gray}$-functors is a $\mathbf{Gray}$-category $\mathbf{Gray}$-$\mathbf{Cat}_{3}$, and in this case composition of $1$-cells and $2$-cells coincides with their usual composition as trinatural transformations and trimodifications.
		\item The triequivalence $\eta_\mathfrak{A}: \mathfrak{A} \rightarrow \mathbf{Gr}\left(\mathfrak{A}\right)$ defined in Theorem 10.11 of \cite{Gurski Coherence in Three Dimensional Category Theory} has a retraction $\eta_\mathfrak{A}^{*}$, and these morphisms are part of a biequivalence internal to $\mathbf{Tricat}_{3}$.
		\item The inclusion $\mathbf{Gray}$-$\mathbf{Cat}_{3} \rightarrow \mathbf{Tricat}_{3}$ has a left triadjoint given by $\mathbf{Gr}$, with the component of the unit of this triadjunction given at a tricategory $\mathfrak{A}$ by $\eta_\mathfrak{A}: \mathfrak{A} \rightarrow \mathbf{Gr}\left(\mathfrak{A}\right)$.
		\item If $\mathfrak{A}$ is a $\mathbf{Gray}$-category then the retraction $\eta_\mathfrak{A}^{*}$ is a $\mathbf{Gray}$-functor.
		\item If $\mathfrak{A}$ is a cofibrant $\mathbf{Gray}$-category then
		
		\begin{enumerate}
			\item the trihomomorphism $\eta_\mathfrak{A}$ is equivalent to a $\mathbf{Gray}$-functor $\bar{\eta}_\mathfrak{A}$ still satisfying  $\eta_\mathfrak{A}\bar{\eta}_\mathfrak{A} = 1_\mathfrak{A}$ via a pseudo-icon equivalence $e: \eta_\mathfrak{A} \rightarrow \bar{\eta}_\mathfrak{A}$.
			\item the pseudo-icon equivalence \begin{tikzcd}
				1_{\mathbf{Gr}\left(\mathfrak{A}\right)} \arrow[rr, "u"] && \eta_{\mathfrak{A}}\eta_{\mathfrak{A}}^{*} \arrow[rr, "e.1"] && \bar{\eta}_{\mathfrak{A}}{\eta}_{\mathfrak{A}}^{*}
			\end{tikzcd} is isomorphic to a semi-strict pseudo-icon equivalence $\bar{u}$ via a costrict trimodification.
			\item there is a biequivalence $\eta_\mathfrak{A}^{*} \dashv \bar{\eta}_\mathfrak{A}$ internal to the $\mathbf{Gray}$-category $\mathbf{Gray}$-$\mathbf{Cat}_{3}$ whose counit is the identity and whose unit is the semi-strict pseudo-icon $\bar{u}$.
		\end{enumerate} 
	\end{enumerate}
\end{proposition}

\begin{proof}
	See \cite{Low dimensional structures formed by tricategories} for parts (1) and (3). Part (2) follows from part (1) by observing that any weakness in the tricategory $\mathbf{Tricat}_{3}$ comes from its objects and morphisms being weak. In particular, $\mathbf{Gray}$-functors compose associatively and if $\mathfrak{B}$ is a $\mathbf{Gray}$-category then $\mathbf{Tricat}_{s}[\mathfrak{A}, \mathfrak{B}]$ is a $\mathbf{Gray}$-category whose interchanger on a pair of interchangable trimodifications $\alpha$ and $\beta$ is a perturbation whose components at $X \in \mathfrak{A}$ are given by the interchanger in $\mathfrak{B}$ of the $2$-cell components $\alpha_{X}$ and $\beta_{X}$. But by the interchange axioms of a $\mathbf{Gray}$-category, this will be the identity if $\alpha_{X}$ and $\beta_{X}$ are identities, as is the case for costrict trimodifications. For part (4), let $\mathfrak{B}$ be a $\mathbf{Gray}$-category and observe that in the composite depicted below, $I$ is the identity on hom $2$-categories by definition and essentially surjective on objects by Proposition \ref{trihomomorphism out of cofibrant Gray category}. Part (4) then follows by noticing that the restriction $\mathbf{Gray}\text{-}\mathbf{Cat}_{3}\left(\eta_{\mathfrak{A}}{,}\mathfrak{B}\right)$ is also a biequivalence since $\eta_{\mathfrak{A}}$ is a biequivalence.
	
	$$\begin{tikzcd}
		\mathbf{Gray}\text{-}\mathbf{Cat}_{3}\left(\mathbf{Gr}\left(\mathfrak{A}\right){,}\mathfrak{B}\right)
		\arrow[r, "I"] & 
		\mathbf{Tricat}_{3}\left(I\mathbf{Gr}\left(\mathfrak{A}\right){,}I\mathfrak{B}\right) \arrow[rrr, "
		\mathbf{Tricat}_{3}\left(\eta_{\mathfrak{A}}{,}\mathfrak{B}\right)"] &&& 
		\mathbf{Tricat}_{3}\left(\mathfrak{A}{,}I\mathfrak{B}\right) 
	\end{tikzcd}$$
	
	\noindent Part (5) is clear from the definition of $\eta_\mathfrak{A}^{*}$ given in \cite{Gurski Coherence in Three Dimensional Category Theory}. In part 6, part (a) follows from Proposition \ref{trihomomorphism out of cofibrant Gray category} and part (b) follows from Proposition \ref{Trinatural transformations out of cofibrant Gray categories}. Finally, part (6 c) follows from part (6 b) by Theorem 4.5 of \cite{Gurski biequivalence in tricategories}.
\end{proof}

\section{Solutions to the failure of semi-strict trinatural transformations to be closed under composition}\label{Section Solutios to the Failure of Semi-strict trinatural transformations to be closed under composition}

\noindent This Section explores alternative solutions to semi-strict trinatural transformations failing to be closed under composition, recorded in Lemma \ref{tautological condition for failure of semi-stricts to be closed under composition}. All but the last of these solutions involve making simplifying assumptions on either the domain or codomain.

\begin{itemize}
	\item Corollary \ref{3 -category suffices for semi-strictness} identifies an assumption on the codomain $\mathbf{Gray}$-category $\mathfrak{B}$ which suffices for semi-strict trinatural transformations to be closed under composition. This reiterates that failure of semi-stricts to be closed under composition is entirely to do with non-trivial interchangers in the codomain. 
	\item Proposition \ref{Gray categories with only trivial composites} identifies simplifying assumptions on the domain $\mathfrak{A}$ which suffice for semi-stricts to be closed under composition.
	\item Proposition \ref{alternative composition if domain is cofibrant resolves semi-strict's failure to be closed under composition} shows that if $\mathfrak{A}$ is cofibrant then it is possible to define a new composition structure on $\mathbf{Tricat}_{s}\left(\mathfrak{A}, \mathfrak{B}\right)$ under which semi-strict trinatural transformations will indeed be closed. This composition will moreover be shown to extend to a triequivalent $\mathbf{Gray}$-category structure.
	\item Finally, we suggest focusing on the weaker property of semi-strict \emph{decomposability} rather than semi-strictness itself. A precise definition will be given in Definition \ref{coherent hom def}. In a forthcoming paper \cite{Miranda semi-strictly generated closed structure on Gray-Cat} we will exhibit this as the internal hom fragment of a closed structure on $\mathbf{Gray}$-$\mathbf{Cat}$.
	\item The last two solutions will be applied in Theorem \ref{trihomomorphism classifier} to categorify the isomorphism between hom $2$-categories in the two-dimensional strictification triadjunction discussion in Subsection \ref{Subsection strictification of bicategories as a triadjoint}, to a triequivalence that might reasonably underlie a `semi-strictification tetra-adjunction'.
\end{itemize} 

\subsection{If the codomain is strict, or the domain has only trivial $1$-cell composites}

In this subsection, we consider \begin{tikzcd}
	F \arrow[r, "p"] & G \arrow[r, "q"] & H
\end{tikzcd} as in Lemma \ref{tautological condition for failure of semi-stricts to be closed under composition}. 

\begin{corollary}\label{3 -category suffices for semi-strictness}
	If $\mathfrak{B}$ is a $3$-category, hence if it has identity interchangers, then the composite $q \circ p$ is again semi-strict.
\end{corollary}

\begin{proof}
	In Lemma \ref{tautological condition for failure of semi-stricts to be closed under composition} part (2), the interchanger ${\left(q_{g}\right)}_{\left(p_{f}\right)}$ itself will be an identity.
\end{proof}

\noindent Proposition \ref{Gray categories with only trivial composites} identifies a special condition on $\mathfrak{A}$ instead of on $\mathfrak{B}$ under which semi-strict trinatural transformations are also closed under composition.

\begin{definition}\label{only trivial 1-cell composites definition}
	We say that a category $\mathcal{A}$ has \emph{only trivial composites} if for any diagram \begin{tikzcd}
		X \arrow[r, "f"] & Y \arrow[r, "g"] &Z
	\end{tikzcd} in $\mathcal{A}$, either $g$ is an identity or $f$ is an identity.. We say that a $\mathbf{Gray}$-category has \emph{only trivial $1$-cell composites} if its underlying category has only trivial composites.
\end{definition}

\begin{proposition}\label{Gray categories with only trivial composites}
	Suppose that a $\mathbf{Gray}$-category $\mathfrak{A}$ has only trivial $1$-cell composites. Let $\mathfrak{B}$ be any $\mathbf{Gray}$-category and let $F, G: \mathfrak{A} \rightarrow \mathfrak{B}$ be $\mathbf{Gray}$-functors.
	
	\begin{enumerate}
		\item A trinatural transformation $p: F \rightarrow G$ is semi-strict if and only if it has identity unitors.
		\item If $H: \mathfrak{A} \rightarrow \mathfrak{B}$ is another $\mathbf{Gray}$-functor and $p: F \rightarrow G$ and $q: G \rightarrow H$ are semi-strict, then their composite $q \circ p$ is also semi-strict.
	\end{enumerate} 
\end{proposition}

\begin{proof}
	For part (1) it suffices to show that the compositors $p_{g, f}$ for $p$ are also identities. Since $\mathfrak{A}$ has only trivial $1$-cell composites, this follows from the left and right unit laws for $p$. Part (2) follows from part (1), since as we have already observed above, the unitors for $q \circ p$ will be identities since those for $q$ and $p$ are by assumption.  
\end{proof}

\begin{example}
	The $\mathbf{Gray}$-categories $\mathfrak{P}_{n}$ of Subsection \ref{Path Objects Section} have only trivial $1$-cell composites.
\end{example}

\subsection{Redefining composition when the domain $\mathbf{Gray}$-category is cofibrant}
\noindent The following proposition identifies another situation in which semi-strict trinatural transformations do fail to be closed under their usual composition, but a different composition law can be defined under which they will be closed. This new composition law is also shown to extend to a triequivalent $\mathbf{Gray}$-category structure on $\mathbf{Tricat}_{s}\left(\mathfrak{A}, \mathfrak{B}\right)$. A $\mathbf{Gray}$-category structure on semi-strict trinatural transformations was assumed in Chapter 15 of \cite{Gurski Coherence in Three Dimensional Category Theory}, and in Remark \ref{salvaging Gurski Chapter 15} we record that their results remain true with composition defined in this alternative way.

\begin{proposition}\label{alternative composition if domain is cofibrant resolves semi-strict's failure to be closed under composition}
	Let $\mathfrak{A}$ be a cofibrant $\mathbf{Gray}$-category and for $\mathbf{Gray}$-functors $P, Q: \mathfrak{A} \rightarrow \mathfrak{B}$, let $\left(P, Q\right)$ denote the full sub-$2$-category of $\mathbf{Tricat}_{s}\left(\mathfrak{A}, \mathfrak{B}\right)\left(P, Q\right)$ on semi-strict trinatural transformations, and let $\left(P, Q\right)_{0}$ denote just the set of semi-strict trinatural transformations. For $\mathbf{Gray}$ functors $F, G, H : \mathfrak{A} \rightarrow \mathfrak{B}$, define a function $\bar{\left(-\right)\circ\left(-\right)}: {\left(G, H\right)}_{0} \times {\left(F, G\right)}_{0} \rightarrow {\left(F, H\right)}_{0}$ by setting $\bar{q \circ p}$ to be the semi-strict trinatural transformation constructed from the usual composite \begin{tikzcd}
		F \arrow[r, "p"] & G \arrow[r, "q"] &H
	\end{tikzcd} via Proposition \ref{Trinatural transformations out of cofibrant Gray categories}.  
	
	\begin{enumerate}
		\item $\bar{q \circ p}$ agrees with the usual composite $q \circ p$ on $1$-cell components for objects, $2$-cell components for generating morphisms and $3$-cell components for generating $2$-cells.
		\item The family of functions of the form $\bar{\left(-\right)\circ\left(-\right)}: {\left(G, H\right)}_{0} \times {\left(F, G\right)}_{0} \rightarrow {\left(F, H\right)}_{0}$ define the composition maps of a category structure on the graph whose objects are $\mathbf{Gray}$-functors from $\mathfrak{A}$ to $\mathfrak{B}$ and whose edges are semi-strict trinatural transformations.
		\item The functions of part (1) are the assignments on objects of $2$-functors  ${\left(G, H\right)} \otimes {\left(F, G\right)} \rightarrow {\left(F, H\right)}$.
		\item The $2$-functors of part (3) are the composition operations of a $\mathbf{Gray}$-category structure on the $2$-$\mathbf{Cat}$ enriched graph $[\mathfrak{A}, \mathfrak{B}]_{ss}$ whose objects are $\mathbf{Gray}$-functors from $\mathfrak{A}$ to $\mathfrak{B}$ and whose hom $2$-categories are of the form $\left(F, G\right)$.
		\item The inclusions $\left(F, G\right) \rightarrow \mathbf{Tricat}_{s}\left(\mathfrak{A}, \mathfrak{B}\right)\left(F, G\right)$ are the actions between hom-$2$-categories of an identity on objects trihomomorphism whose only weak components are their compositors on \begin{tikzcd}
			F \arrow[r, "p"] & G \arrow[r, "q"] &H
		\end{tikzcd}. 
		\item The trihomomorphisms of part (5) are triequivalences.
	\end{enumerate}
\end{proposition}

\begin{proof}
	Part (1) evidently follows from the construction described in Proposition \ref{Trinatural transformations out of cofibrant Gray categories}. For part (2), the left and right unit laws are immediate, since identities are semi-strict and the construction of Proposition \ref{Trinatural transformations out of cofibrant Gray categories} is inert on identities. Moreover observe that by part (1), the associativity condition holds on objects, generating morphisms and generating $2$-cells and by Lemma \ref{Gurski Lemma semi strict}, this data uniquely determines a semi-strict trinatural transformation.
	\\
	\\
	\noindent For part (3), we need to describe the rest of the $2$-functor $\left(G ,H\right) \otimes \left(F, G\right) \rightarrow \left(F, H\right)$ on generating morphisms and $2$-cells in the $\mathbf{Gray}$ tensor product. For a generator of the form $\left(1_{p}, \tau\right): \left(p, t\right) \Rightarrow \left(p, t'\right)$ with $\tau$ a trimodification, the trimodification $p\bar{\circ} \tau: p\bar{\circ}t \Rightarrow p\bar{\circ}t'$ is given as depicted below, where $i^{pt}$ and $i^{pt'}$ are the unique costrict trimodifications of Proposition \ref{Trinatural transformations out of cofibrant Gray categories} part (2).
	
	$$\begin{tikzcd}
		p \bar{\circ}t = \bar{pt} \arrow[r, "i^{pt}"] & pt \arrow[r, "p.\tau"]& pt' \arrow[r, "{\left(i^{pt'}\right)}^{-1}"] &\bar{pt'} =p\bar{\circ}t'
	\end{tikzcd}$$
\noindent Meanwhile, the $2$-functor $\left(-\right)\bar{\circ}\left(-\right)$ is defined on generating morphisms $\left(\sigma, 1_{q}\right): \left(s, q\right) \Rightarrow \left(s', q\right)$  and generating $2$-cells in a similar way, namely via whiskering with $i$. Finally, it is clear that these assignments respect the relations in the presentation of the $\mathbf{Gray}$-tensor product since the $i$ are costrict and have identity $3$-cell components on generating morphisms in $\mathfrak{A}$. The associativity and unit laws in the $\mathbf{Gray}$-category structure for part (4) follow for similar reasons.
	\\
	\\
	\noindent For part (5), the component on \begin{tikzcd}
		F \arrow[r, "p"] & G \arrow[r, "q"] &H
	\end{tikzcd} is given by $i^{qp}$. $2$-naturality of this compositor follows by costrictness of $i^{qp}$. Strict unitality and the usual pseudofunctor laws, and by extension the trihomomorphism axioms, are also clear for similar reasons. Finally, part (6) follows from part (5); the trihomomorphism is bijective on objects and an isomorphism between hom-categories by definition, and it is locally essentially surjective on objects by Proposition \ref{Trinatural transformations out of cofibrant Gray categories}.
\end{proof}

\begin{remark}\label{salvaging Gurski Chapter 15}
	The results in Chapter 15 of \cite{Gurski Coherence in Three Dimensional Category Theory} are true when the $\mathbf{Gray}$-category structure on $\mathbf{Codsc}\left(\mathfrak{K}\right)$ is interpreted as per Proposition \ref{alternative composition if domain is cofibrant resolves semi-strict's failure to be closed under composition}. The $\mathbf{Gray}$-category $\Delta_\text{ps}^{G}$ is defined in Definition 11.4. Observe that its underlying sesquicategory is free on a $2$-computad, with the only relations being between $3$-cells, and corresponding to equations needed for certain data to be adjoint equivalences. Proposition \ref{alternative composition if domain is cofibrant resolves semi-strict's failure to be closed under composition} therefore applies, and in particular the calculations in Chapter 15 of \cite{Gurski Coherence in Three Dimensional Category Theory} transfer along the triequivalence of Proposition \ref{alternative composition if domain is cofibrant resolves semi-strict's failure to be closed under composition} part (6).
\end{remark}

\subsection{Semi-strict decomposability and the semi-strictification `tetra-adjunction'}\label{subsection semi-strict decomposability}

\noindent The following definition addresses the issue that semi-strict trinatural transformations fail to be closed under composition without making any assumptions on the domain or codomain $\mathbf{Gray}$-categories.

\begin{definition}\label{coherent hom def}
	Let $\mathfrak{A}$ and $\mathfrak{B}$ be $\mathbf{Gray}$-categories. Call a trinatural transformation \emph{semi-strictly decomposable} if it admits a factorisation into a finite composite of semi-strict trinatural transformations. Let $[\mathfrak{A}, \mathfrak{B}]_\text{ssg}$ be the sub-$\mathbf{Gray}$-category of $\mathbf{Tricat}_{s}\left(\mathfrak{A}, \mathfrak{B}\right)$ on all objects and just those trinatural transformations which are semi-strictly decomposable. Call the $\mathbf{Gray}$-category $[\mathfrak{A}, \mathfrak{B}]_{\text{ssg}}$ the \emph{semi-strictly generated hom-$\mathbf{Gray}$-category} associated to $\mathfrak{A}$ and $\mathfrak{B}$.
\end{definition}

\begin{remark}\label{coherent hom via bosa-l2ff}
	It is important to note that semi-strict decomposability has been defined as a property of a trinatural transformation, rather than as structure which records any particular factorisation into a specific composite of semi-strict trinatural transformations. Another perspective on $[\mathfrak{A}, \mathfrak{B}]_\text{ssg}$ is that it appears in the factorisation depicted below.
	
	$$\begin{tikzcd}
		\mathcal{F}\left(\mathbf{G}\right) \arrow[r, "L"] & {[}\mathfrak{A}{,}\mathfrak{B}{]}_\text{ssg} \arrow[r, "R"] & \mathbf{Tricat}_{s}\left(\mathfrak{A}{,}\mathfrak{B}\right)
	\end{tikzcd}$$
	\\
	\noindent In this factorisation
	
	\begin{itemize}
		\item $\mathbf{G}$ is the reflexive graph whose vertices are $\mathbf{Gray}$-functors $\mathfrak{A} \rightarrow \mathfrak{B}$ and whose arrows are semi-strict trinatural transformations,
		\item  $\mathcal{F}: \mathbf{RGrph} \rightarrow \mathbf{Cat}$ denotes the functor which sends a reflexive graph to its free category,
		\item  $L$ is a bijective on objects and surjective on morphisms $\mathbf{Gray}$-functor,
		\item $R$ is a $\mathbf{Gray}$-functor that is injective on morphisms and fully faithful on $2$-cells and $3$-cells,
		\item $RL$ is transpose to the inclusion of $\mathbf{G}$ as a sub-reflexive graph of $\mathbf{Tricat}_{s}[\mathfrak{A}, \mathfrak{B}]$.		
	\end{itemize}   
	\noindent In this way, semi-strictly decomposable trinatural transformation can also be seen as equivalence classes of non-empty paths of semi-strict trinatural transformations, where two paths $\left(p_{1}, ... , p_{m}\right)$ and $\left(q_{1}, ..., q_{n}\right)$ are equivalent if they compose to give equal trinatural transformations.
\end{remark}

\noindent Definition \ref{coherent hom def} allows us to state a universal property of $\eta_{\mathfrak{A}}: \mathfrak{A} \rightsquigarrow \mathbf{Gr}\left(\mathfrak{A}\right)$ as a trihomomorphism classifier, which takes into account more general $\left(3, k\right)$-transfors than Proposition \ref{tricategory of tricategories}. This is achieved in Theorem \ref{trihomomorphism classifier}, before which we will need Lemma \ref{Lemma semi-strictification tetra adjunction}.

\begin{lemma}\label{Lemma semi-strictification tetra adjunction}
	Let $\mathfrak{A}$ be a tricategory and $\mathfrak{B}$ be a $\mathbf{Gray}$-category. Let $\eta_\mathfrak{A}: \mathfrak{A} \rightsquigarrow \mathbf{Gr}\left(\mathfrak{A}\right)$ be the triequivalence described in \cite{Gurski Coherence in Three Dimensional Category Theory} Theorem 10.11. Then consider the following composite, in which $R$ is as defined in Remark \ref{coherent hom via bosa-l2ff} and for any $\mathbf{Gray}$-category $\mathfrak{C}$, $R\mathfrak{C} := \mathfrak{C}$.
	\\
	$$\begin{tikzcd}
		{[}\mathbf{Gr}\left(\mathfrak{A}\right){,}\mathfrak{B}{]}_{\text{ssg}} \arrow[rr, "R"] && \mathbf{Tricat}\left(R\mathbf{Gr}\left(\mathfrak{A}\right){,}R\mathfrak{B}\right) \arrow[rr, "\mathbf{Tricat}\left(\eta_\mathfrak{A}{,}R\mathfrak{B}\right)"] && \mathbf{Tricat}\left(\mathfrak{A}{,}R\mathfrak{B}\right)
	\end{tikzcd}$$
	
	\begin{enumerate}
		\item $\mathbf{Tricat}\left(\eta_\mathfrak{A}, R\mathfrak{B}\right)$ is a $\mathbf{Gray}$-functor.
		\item  $\mathbf{Tricat}\left(\eta_\mathfrak{A}, R\mathfrak{B}\right)$ is a triequivalence.
		\item $R$ is biessentially surjective on objects.
		\item $R$ is locally essentially surjective on objects.
	\end{enumerate}
\end{lemma}

\begin{proof}
	Part (1) follows from Theorem A.6 of \cite{Buhne PhD}. Since $\eta_\mathfrak{A}$ is a triequivalence by Theorem 10.9 of \cite{Gurski Coherence in Three Dimensional Category Theory}, part (2) follows from Theorem A.9 of \cite{Buhne PhD}. Part (3) follows from Proposition \ref{trihomomorphism out of cofibrant Gray category} since $\mathbf{Gr}\left(\mathfrak{A}\right)$ is a cofibrant $\mathbf{Gray}$-category. For the same reason, part (4) follows from Proposition \ref{Trinatural transformations out of cofibrant Gray categories}.
\end{proof}

\begin{theorem}\label{trihomomorphism classifier}
	Consider the composable pair of $\mathbf{Gray}$-functors displayed in Lemma \ref{Lemma semi-strictification tetra adjunction}. \begin{enumerate}
		\item Their composite is a triequivalence.
		\item The further composition with the trihomomorphism $[\mathbf{Gr}\left(\mathfrak{A}\right), \mathfrak{B}]_\text{ss} \rightsquigarrow [\mathbf{Gr}\left(\mathfrak{A}\right), \mathfrak{B}]_\text{ssg}$ of Proposition \ref{alternative composition if domain is cofibrant resolves semi-strict's failure to be closed under composition} is also a triequivalence.
		\item Both of these triequivalences are moreover biessentially surjective on objects, locally essentially surjective on objects, and isomorphisms between hom-categories.
	\end{enumerate}
\end{theorem}

\begin{proof}
	Since $R$ is by definition fully faithful on $2$-cells and $3$-cells, part (1) follows from parts (3) and (4) of Lemma \ref{Lemma semi-strictification tetra adjunction}. Parts (2) and (3) then follow by closure under composition of the properties in question.
\end{proof}

\begin{remark}\label{Remark semi-strictification tetra-adjunction}
	Theorem \ref{trihomomorphism classifier} part (1) is suggestive of a `semi-strictification tetra-adjunction' between $\mathbf{Tricat}$ and its sub-tetracategory $\mathbf{Gray}\text{-}\mathbf{Cat}_{\text{ss}}$ on $\mathbf{Gray}$-categories, $\mathbf{Gray}$-functors, and semi-strictly decomposable trinatural transformations. The left tetra-adjoint would now be given by $\mathbf{Gr}\left(-\right)$ in place of $\mathbf{st}_{2}$. It classifies both trihomomorphisms via $\mathbf{Gray}$-functors, and trinatural transformations via semi-strict ones. Motivated by the universal property, one can present a $\mathbf{Gray}$-category  $\mathbf{st}_{3}\left(\mathfrak{A}\right)$ in a similar way to what we have described in Definition \ref{presentation of universal 2-category recipient of a pseudofunctor}. This $\mathbf{Gray}$-category $\mathbf{st}_{3}\left(\mathfrak{A}\right)$ is a biased version of $\mathbf{Gr}\left(\mathfrak{A}\right)$, analogously to what we discussed in Remark \ref{comparing pseudofunctor classifier to cofibrant replacement}. See Subsection 2.2.3 of \cite{Miranda PhD} for a detailed presentation of $\mathbf{st}_{3}\left(\mathfrak{A}\right)$.
	\\
	\\
	\noindent For reasons similar to those discussed in Remark \ref{tri-fully-faithful and triessential image}, $\mathbf{Gr}\left(-\right)$ should be tetra-fully faithful, and the restriction of $R$ to cofibrant $\mathbf{Gray}$-categories should be a tetra-equivalence. Moreover, by Theorem \ref{trihomomorphism classifier} part (2), the `sub-tetracategory on cofibrant $\mathbf{Gray}$-categories' should have an alternative composition structure in which all $2$-cells are semi-strict. Theorem \ref{trihomomorphism classifier} part (3) records the fact that the action between hom-$\mathbf{Gray}$-categories is stricter than what one might expect in a general tetra-adjunction. The analogous map in the two-dimensional setting was seen to be an isomorphism of $2$-categories rather than merely a biequivalence. In a forthcoming paper \cite{Miranda semi-strictly generated closed structure on Gray-Cat} we will form a four dimensional categorical structure $\mathbf{Gray}$-$\mathbf{Cat}_\text{ssg}$, however a tetracategory structure on the $4$-globular set $\mathbf{Tricat}$ remains open \cite{Gurski Coherence in Three Dimensional Category Theory} \cite{Trimble Tetracategory}.
\end{remark}

\begin{remark}\label{dimension three strict, Gray, semi-strict, fully weak}
	An analogue to the perspective offered in Remark \ref{dimension two strict, semi-strict, fully weak} can also be given in dimension three. \emph{Strict} is now taken to mean enriched over $2$-$\mathbf{Cat}$ with its cartesian closed structure, while enrichment over the $\mathbf{Gray}$ tensor product provides another intermediate level of strictness. This time, the entries in the semi-strict column do not combine to give a $\mathbf{Gray}$-category since semi-strict trinatural transformations are not closed under composition. However, closing them under composition results in $\mathbf{Gray}$-categories of the form $[\mathfrak{A}, \mathfrak{B}]_{ssg}$.
	
	\begin{center}
		\noindent\begin{tabularx}{1.05\textwidth} { 
				| >{\raggedright\arraybackslash}X 
				| >{\raggedright\arraybackslash}X 
				| >{\raggedright\arraybackslash}X
				| >{\centering\arraybackslash}X
				| >{\raggedleft\arraybackslash}X
				|}
			\hline
			Dimension & Strict & $\mathbf{Gray}$-enriched & Semi-strict & Fully-weak \\
			\hline
			\hline
			$n = 0$  & $3$-categories & $\mathbf{Gray}$-categories & $\mathbf{Gray}$-categories & tricategories  \\
			\hline
			$n = 1$ & $3$-functors & $\mathbf{Gray}$-functors & $\mathbf{Gray}$-functors & trihomomorphisms \\
			\hline
			$n = 2$ & $3$-natural transformations & $\mathbf{Gray}$-natural transformations & semi-strict trinatural transformations & trinatural transformations\\
			\hline
			$n = 3$ & $3$-modifications & $\mathbf{Gray}$-modifications & trimodifications & trimodifications\\
			\hline
			$n = 4$ & perturbations & perturbations & perturbations & perturbations\\
			\hline
		\end{tabularx}
	\end{center}
	
	\noindent In this setting, we see that
	
	\begin{enumerate}
		\item For $n >0$ the difference between the $2$-$\mathbf{Cat}$ enriched and $\mathbf{Gray}$-enriched notions is nominal, in the sense that if the objects are $\mathbf{Gray}$-categories with strict interchangers, i.e. $3$-categories, then there is no difference.
		\item The entries in the `semi-strict' column are as strict as entries in the `Fully weak' column become under three-dimensional strictification,
		\item The entries in the `semi-strict' column have operational coherences given by identities, but may have non-trivial free coherences.
	\end{enumerate}
\end{remark}

\section{Conclusion}

\noindent In this paper we have studied how $n \in\{2, 3\}$ dimensional strictification constructions in which the resulting structure is free in codimension one can be extended to strictify certain coherence data in $\left(n, k\right)$-transfors. We have identified the coherences which strictify in this way as those which mediate the operations of their domain $n$-dimensional categorical structure. On the other hand, coherences which mediate varying a $k$-cell along a $\left(k+1\right)$-cell, or \emph{free} coherences, remain weak after strictification is applied. Using the generalised path objects for $\mathbf{Gray}$-categories, introduced in Subsection \ref{Path Objects Section}, we showed how extending $\mathbf{Gr}$ to $\left(3, k\right)$-transfors for positive $k$ can be reduced to the strictification of $\mathbf{Gray}$-functors. We have examined various solutions to the failure of semi-strict trinatural transformations being closed under composition, in particular recovering the results in Chapter 15 of \cite{Gurski Coherence in Three Dimensional Category Theory}. One such solution, given in Definition \ref{coherent hom def}, makes no assumptions on the source and target $\mathbf{Gray}$-category. This is used to give certain triequivalences which may underlie a semi-strictification tetra-adjunction', in Theorem \ref{trihomomorphism classifier}. These results reinforce the connection between weak functor classifiers, strictification, and cofibrant replacement. In a forthcoming paper \cite{Miranda semi-strictly generated closed structure on Gray-Cat}, this semi-strictly generated hom of $\mathbf{Gray}$-categories ${[\mathfrak{A}, \mathfrak{B}]}_\text{ssg}$ will be shown to equip the category $\mathbf{Gray}$-$\mathbf{Cat}$ with a closed structure.

\end{document}